\RequirePackage{fix-cm}
\documentclass[smallextended]{svjour3}       
\smartqed  
\usepackage{graphicx}
\usepackage{amsmath,amsfonts}
\usepackage{amssymb}
\usepackage{latexsym}
\usepackage{mathrsfs}
\usepackage{colortbl}
\usepackage{amscd}
\usepackage{tikz-cd}
\usepackage{comment}
\usepackage{url}
\usepackage[subrefformat=parens]{subcaption}
\usepackage{mathtools}
\usepackage[section]{placeins}

\allowdisplaybreaks[3]

  \journalname{}

\def\diam{\mathop{\mathrm{diam}}}

\def\div{\mathop{\mathrm{div}}}
\def\divh{\mathop{\mathrm{div}_h}}
\def\rot{\mathop{\mathrm{rot}}}

\def\<{\mathop{\textless}}
\def\>{\mathop{\textgreater}}

\spnewtheorem{thr}{Theorem}{\bf}{\it}
\spnewtheorem{defi}{Definition}{\bf}{\it}
\spnewtheorem{lem}{Lemma}{\bf}{\it}
\spnewtheorem{coro}{Corollary}{\bf}{\it}
\spnewtheorem{assume}{Assumption}{\bf}{\it}
\spnewtheorem{ex}{Example}{\bf}{\it}
\spnewtheorem{Case}{Case}{\bf}{\it}
\spnewtheorem*{pf*}{Proof}{\bf}{\rm}
\spnewtheorem*{rem*}{Remark:}{\it}{\it}
\spnewtheorem*{ex*}{Example:}{\it}{\it}
\spnewtheorem{Cond}{Condition}{\bf}{\it}
\spnewtheorem{rem}{Remark}{\it}{\it}

\spnewtheorem*{lem1*}{Lemma 1}{\bf}{\rm}

\newcounter{sone}
\newcounter{stwo}
\newcounter{sthree}
\newcounter{sfour}
\newcounter{sfive}
\newcounter{ssix}
\newcounter{lone}
\newcounter{ltwo}
\newcounter{lthree}
\newcounter{lfour}
\newcounter{lfive}
\newcounter{lsix}
\makeatletter
    
    \@addtoreset{equation}{section}

\begin{document}

\title{Anisotropic weakly over-penalised symmetric interior penalty method for the Stokes equation 
}

\titlerunning{Anisotropic WOPSIP method for the Stokes equation}        

\author{Hiroki Ishizaka  
}


\institute{Hiroki Ishizaka \at
              Team fem, Matsuyama, Japan \\
              \email{h.ishizaka005@gmail.com}           
}

\date{Received: date / Accepted: date}

\maketitle

\begin{abstract}
In this study, we investigate an anisotropic weakly over-penalised symmetric interior penalty method for the Stokes equation {on convex domains}. Our approach is a simple discontinuous Galerkin method similar to the Crouzeix--Raviart finite element method. As our primary contribution, we show a new proof for the consistency term, which allows us to obtain an estimate of the anisotropic consistency error. The key idea of the proof is to apply the relation between the Raviart--Thomas finite element space and a discontinuous space. While inf-sup stable schemes of the discontinuous Galerkin method on shape-regular mesh partitions have been widely discussed, our results show that the Stokes element satisfies the inf-sup condition on anisotropic meshes. {Furthermore, we provide} an error estimate in an energy norm on anisotropic meshes. In numerical experiments, we compare calculation results for standard and anisotropic mesh partitions.

\end{abstract}

\section{Introduction} \label{intro}
In this study, we investigate a weakly over-penalised symmetric interior penalty (WOPSIP) method for the Stokes equations on anisotropic meshes. Brenner et al. first proposed a WOPSIP method \cite{BreOweSun08}, and several further works have considered similar techniques \cite{BarBre14,Bre15,BreOweSun12}. WOPSIP methods have two main advantages compared to standard symmetric interior penalty discontinuous Galerkin (dG) methods \cite{PieErn12,Riv08}. The first is that they are stable for any penalty parameter. Moreover, they work on nonconforming mesh partitions. Meanwhile, the drawback of the standard symmetric interior penalty discontinuous Galerkin method is that it requires tuning the penalty parameter for stability. Moreover, applying nonconforming meshes can be difficult in the classical Crouzeix--Raviart (CR) nonconforming finite element method (\cite{CroRav73}). In the present work, we explore an error analysis on conformal meshes for simplicity. We briefly consider the case of nonconforming meshes in Section \ref{sec=conc}. The WOPSIP method is similar to the classical CR finite element method and thus has similar features such as inf-sup stability on anisotropic meshes. However, several studies have imposed the condition of shape regularity to a family of meshes, i.e., triangles or tetrahedra cannot be overly flat in a shape-regular family of triangulations.

Meanwhile, anisotropic meshes are effective for problems in which the solution has anisotropic behaviour in some direction of the domain. The anisotropic meshes do not satisfy the shape-regular condition or include elements with large aspect ratios. Given this background, several anisotropic finite element methods have been developed in recent years. \cite{AcoApe10,AcoDur99,Ape99,ApeKemLin21,Ish22a,Ish21,IshKobTsu21a,IshKobTsu21b,IshKobTsu21c,IshKobSuzTsu21d}. These methods aim to obtain optimal error estimates under the \textit{semi-regular condition} defined in Assumption \ref{ass1} \cite{IshKobTsu21a} or the \textit{maximum-angle condition}, which allows us to use anisotropic meshes; see {Babu$\rm{\breve{s}}$ka et al}. \cite{BabAzi76} for two-dimensional and K$\rm{\check{r}}$\'{i}$\rm{\check{z}}$ek’s work \cite{Kri92} three-dimensional cases. In this study, we consider the anisotropic WOPSIP method for the Stokes equation and present optimal error estimates in the energy norm.

The WOPSIP method is nonconforming. Therefore, an error between the exact and WOPSIP finite element approximation solutions for the velocity with an energy norm and the pressure with the $L^2$-norm is divided into two parts. One part is an optimal approximation error in discontinuous finite element spaces, and the other is a consistency error term. For the former, the first-order CR interpolation errors (Theorem \ref{thr2}) and the error estimate of the $L^2$-projection (Theorem \ref{thr1}) are used. However, estimating the consistency error term on anisotropic meshes is challenging. Barker and Brenner \cite{BarBre14} apply the trace theorem when $d=2$. Therefore, the shape regularity condition on meshes may be unavoidable with their proposed technique \cite{BarBre14}. To overcome this difficulty, we use the relation between the lowest order Raviart--Thomas (RT) finite element interpolation and the discontinuous space; see Lemma \ref{lem3}. This relation derives an optimal error estimate of the consistency error (Lemma \ref{lem9}).

The remainder of this study is organised as follows. In Section 2, we introduce the (scaled) Stokes equation with the Dirichlet boundary condition, its weak formulation, and finite element settings for the WOPSIP method. In Section 3, we show the discrete Poincar\'e inequality. In Section 4, we present the proposed WOPSIP approximation for the continuous problem and discuss its stability and error estimates. In Section 5, we provide the results of a numerical evaluation. Finally, we conclude by noting some limitations of the present work and suggesting some possible avenues for future research in Section 6. Throughout, we denote by $c$ a constant independent of $h$ (defined later) and of the angles and aspect ratios of simplices unless specified otherwise, {and all constants $c$ are bounded if the maximum angle is bounded}. These values may change in each context. The notation $\mathbb{R}_+$ denotes the set of positive real numbers.

\section{Preliminaries}

\subsection{Weak formulation}
We consider the following problem. Let $\Omega \subset \mathbb{R}^d$, $d \in \{ 2 , 3 \}$ be a bounded polyhedral domain. Furthermore, we assume that $\Omega$ is convex if necessary. The (scaled) Stokes problem is to find $(u,p): \Omega \to \mathbb{R}^d \times \mathbb{R}$ such that
\begin{align}
\displaystyle
- \nu \varDelta u + \nabla p = f \quad \text{in $\Omega$}, \quad \div u  = 0 \quad \text{in $\Omega$}, \quad u = 0 \quad \text{on $\partial \Omega$}, \label{intro1}
\end{align}
where $\nu$ is a nonnegative parameter and $f:\Omega \to \mathbb{R}^d$ is a given function.

To deduce a weak form of the continuous problem \eqref{intro1}, we consider function spaces described as follows.
\begin{align*}
\displaystyle
V &:= H_0^1(\Omega)^d, \quad Q := L^2_0(\Omega) := \left \{ q \in L^2(\Omega); \ \int_{\Omega} q dx = 0 \right \},
\end{align*}
with norms:
\begin{align*}
\displaystyle
| \cdot |_V := | \cdot |_{H^1(\Omega)^d}, \quad \| \cdot \|_Q := \| \cdot \|_{L^2(\Omega)}.
\end{align*}

The variational formulation for the Stokes equations \eqref{intro1} is as follows. For any $f \in L^2(\Omega)^d$, find $(u,p) \in V \times Q$ such that
\begin{subequations} \label{stokes1}
\begin{align}
\displaystyle
\nu a(u,\varphi) + b(\varphi , p) &= \int_{\Omega} f \cdot \varphi dx \quad \forall \varphi \in V, \label{stokes1a} \\
b(u , q) &= 0 \quad \forall q \in Q. \label{stokes1b}
\end{align}
\end{subequations}
Here, $a: H^1(\Omega)^d \times H^1(\Omega)^d \to \mathbb{R}$ and $b:H^1(\Omega)^d \times L^2(\Omega) \to \mathbb{R}$ respectively denote bilinear forms defined by
\begin{align*}
\displaystyle
a(v,\psi) := \int_{\Omega} \nabla v : \nabla \psi dx = \sum_{i=1}^d \int_{\Omega}  \nabla v_i \cdot \nabla \psi_i dx, \quad   b(\psi , q) := - \int_{\Omega} \div \psi q dx. 
\end{align*}
Here, the colon denotes the scalar product of tensors.

Using the space of weakly divergence-free functions $V_{\sigma} := \{ v \in V  ; \ b(v,q) = 0 \quad \forall q \in Q \}$, the associated problem to \eqref{stokes1} is then to find $u \in V_{\sigma}$ such that
\begin{align}
\displaystyle
\nu a(u,\varphi)  &= \int_{\Omega} f \cdot \varphi dx \quad \forall \varphi \in V_{\sigma}. \label{stokes2}
\end{align}

The continuous inf-sup inequality
\begin{align}
\displaystyle
\inf_{q \in Q} \sup_{\psi \in V} \frac{b(\psi , q)}{ | \psi |_{V} \| q \|_{Q} }\geq \beta \>0 \label{stokes3}
\end{align}
has been shown to hold. Proofs can be found in John \cite[Theorem 3.46]{Joh16}, Ern and Guermond \cite[Lemma 53.9]{ErnGue21b}, and Girault and Raviart \cite[Lemma 4.1]{GirRav86}.

We set 
\begin{align*}
\displaystyle
H_*^1(\Omega) := H^1(\Omega) \cap L^2_0(\Omega), \quad \mathcal{H}
 := \{ v \in L^2(\Omega)^d: \ \div v = 0, \ v|_{\partial \Omega} \cdot n = 0 \},
 \end{align*}
 where $\div v = 0$ and $v|_{\partial \Omega} \cdot n = 0$ mean that $\int_{\Omega} (v \cdot \nabla ) q dx = 0$ for any $q \in H_*^1(\Omega) $. Then, the following $L^2$-orthogonal decomposition holds.
\begin{align}
\displaystyle
L^2(\Omega)^d = \mathcal{H} \oplus \nabla (H_*^1(\Omega) ), \label{helmdeco}
\end{align}
see \cite[Lemma 74.1]{ErnGue21c}. The $L^2$-orthogonal projection $P_{\mathcal{H}}: L^2(\Omega)^d \to \mathcal{H}$ resulting from this decomposition is often called the \textit{Leray projection}.

\begin{thr}[Stability] \label{thr5}
For any $f \in L^2(\Omega)^d$ or $f \in V^{\prime}$, the weak formulation \eqref{stokes1} of the Stokes problem is well-posed. Furthermore, if  $f \in V^{\prime}$, \begin{align}
\displaystyle
| u |_{V} &\leq \frac{1}{\nu} \| f \|_{V^{\prime}}, \quad
\| p \|_{Q} \leq \frac{2}{\beta} \| f \|_{V^{\prime}}. \label{stokes5}
\end{align}
If $f \in L^2(\Omega)^d$,
\begin{align}
\displaystyle
| u |_{V} &\leq \frac{c}{\nu}\| P_{\mathcal{H}} (f) \|_{L^2(\Omega)^d}, \quad \| p \|_{Q} \leq \frac{c}{\beta} \| f \|_{L^{2}(\Omega)^d}. \label{stokes6}
\end{align}
\end{thr}

\begin{pf*}
A proof was provided by John \cite[Theorem 4.6, Lemma 4.7]{Joh16}.
\qed
\end{pf*}

\subsection{Trace inequality}
The following trace inequality on anisotropic meshes is significant in this study. Some references can be found for the proof. Here, we follow Ern and Guermond \cite[Lemma 12.15]{ErnGue21a}, also see Andreas \cite[Lemma 2.3]{And15} and Kunnert \cite[Lemma 2.2]{Kun01}. We note that although Lemma 12.15 in \cite{ErnGue21a} imposes a shape-regular mesh condition, the condition is easily violated. For a simplex $T \subset \mathbb{R}^d$, let $\mathcal{F}_{T}$ be the collection of the faces of $T$. Let $|\cdot|_d$ denote the $d$-dimensional Hausdorff measure.

\begin{lem}[Trace inequality] \label{lem=trace}
Let  $T \subset \mathbb{R}^d$ be a simplex. There exists a positive constant $c$ such that for any $v = (v^{(1)}, \ldots,v^{(d)})^T \in H^{1}(T)^d$, $F \in \mathcal{F}_{T}$, and $h$,
\begin{align}
\displaystyle
\| v \|_{L^2(F)^d}
\leq c \ell_{T,F}^{- \frac{1}{2}} \left( \| v \|_{L^2(T)^d} + h_{T}^{\frac{1}{2}}  \| v \|_{L^2(T)^d}^{\frac{1}{2}} | v |_{H^1(T)^d}^{\frac{1}{2}} \right),\label{trace}
\end{align}
where $\ell_{T,F} := \frac{d! |T|_d}{|F|_{d-1}}$ denotes the distance of the vertex of $T$ opposite to $F$ to the face. 
\end{lem}

\begin{pf*}
Let $v = (v^{(1)},\ldots,v^{(d)})^T \in H^{1}(T)^d$. Let $P_{F}$ be the vertex of $T$ opposite to $F$. By the same argument for $p=2$ of \cite[Lemma 12.15]{ErnGue21a}, together with the fact that $|x - P_F| \leq h_T$ and $\frac{|F|_{d-1}}{|T|_d} = \frac{d !}{\ell_{T,F}}$, it holds that for $i \in \{ 1 , \ldots , d \}$,
\begin{align*}
\displaystyle
\| v^{(i)} \|^2_{L^2(F)}
&\leq \frac{d !}{\ell_{T,F}} \| v^{(i)} \|^2_{L^2(T)} + \frac{2 (d-1) !}{\ell_{T,F}} h_T \| v^{(i)} \|_{L^2(T)} \| \nabla v^{(i)} \|_{L^2(T)^d}.
\end{align*}
Using the Cauchy--Schwarz inequality, we obtain the target inequality together with the Jensen's inequality.
\qed
\end{pf*}

\begin{rem}
Because $|T|_d \approx h_T^d$ and $|F|_{d-1} \approx h_T^{d-1}$ on the shape-regular mesh, it holds that $\ell_{T,F} \approx h_T$. Then, the trace inequality \eqref{trace} is given as
\begin{align*}
\displaystyle
\| v \|_{L^2(F)^d}
\leq c h_T^{- \frac{1}{2}} \left( \| v \|_{L^2(T)^d} + h_{T}^{\frac{1}{2}}  \| v \|_{L^2(T)^d}^{\frac{1}{2}} | v |_{H^1(T)^d}^{\frac{1}{2}} \right).
\end{align*}
\end{rem}

\subsection{Meshes, mesh faces, averages and jumps} \label{regularmesh}
For simplicity, we consider conformal meshes. Let $\Omega \subset \mathbb{R}^d$ be a bounded polyhedral domain. Let $\mathbb{T}_h = \{ T \}$ be a simplicial mesh of $\overline{\Omega}$ made up of closed $d$-simplices such as
\begin{align*}
\displaystyle
\overline{\Omega} = \bigcup_{T \in \mathbb{T}_h} T,
\end{align*}
with $h := \max_{T \in \mathbb{T}_h} h_{T}$, where $ h_{T} := \diam(T)$. For simplicity, we assume that $\mathbb{T}_h$ is conformal: that is, $\mathbb{T}_h$ is a simplicial mesh of $\overline{\Omega}$ without hanging nodes.

Let $\mathcal{F}_h^i$ be the set of interior faces and $\mathcal{F}_h^{\partial}$ the set of the faces on the boundary $\partial \Omega$. We set $\mathcal{F}_h := \mathcal{F}_h^i \cup \mathcal{F}_h^{\partial}$. For any $F \in \mathcal{F}_h$, we define the unit normal $n_F$ to $F$ as follows. (\roman{sone}) If  $F \in \mathcal{F}_h^i$ with $F = T_1 \cap T_2$, $T_1,T_2 \in \mathbb{T}_h$, let $n_1$ and $n_2$ be the outward unit normals of $T_1$ and $T_2$, respectively. Then, $n_F$ is either of $\{ n_1 , n_2\}$; (\roman{stwo}) If $F \in \mathcal{F}_h^{\partial}$, $n_F$ is the unit outward normal $n$ to $\partial \Omega$.

Here, we consider $\mathbb{R}^q$-valued functions for some $q \in \mathbb{N}$. We define a broken (piecewise) Hilbert space as
\begin{align*}
\displaystyle
H^{1}(\mathbb{T}_h)^q := \{ v \in L^2(\Omega)^q: \ v|_{T} \in H^{1}(T)^q \quad \forall T \in \mathbb{T}_h  \}
\end{align*}
with the norms
\begin{align*}
\displaystyle
{| v |_{H^{1}(\mathbb{T}_h)^q}} &:= \left( \sum_{T \in \mathbb{T}_h} { | v |^2_{H^{1}(T)^q }} \right)^{\frac{1}{2}}.
\end{align*}
When $q=1$, we denote $H^{1}(\mathbb{T}_h) := H^{1}(\mathbb{T}_h)^1 $. Let $\varphi \in H^{1}(\mathbb{T}_h)$. Suppose that $F \in \mathcal{F}_h^i$ with $F = T_1 \cap T_2$, $T_1,T_2 \in \mathbb{T}_h$. Set $\varphi_1 := \varphi{|_{T_1}}$ and $\varphi_2 := \varphi{|_{T_2}}$. Set two nonnegative real numbers $\omega_{T_1,F}$ and $\omega_{T_2,F}$ such that
\begin{align*}
\displaystyle
\omega_{T_1,F} + \omega_{T_2,F} = 1.
\end{align*}
The jump and the skew-weighted average of $\varphi$ across $F$ is then defined as
\begin{align*}
\displaystyle
[\![\varphi]\!] := [\! [ \varphi ]\!]_F := \varphi_1 - \varphi_2, \quad  \{\! \{ \varphi\} \! \}_{\overline{\omega}} :=  \{\! \{ \varphi\} \! \}_{\overline{\omega},F} := \omega_{T_2,F} \varphi_1 + \omega_{T_1,F} \varphi_2.
\end{align*}
For a boundary face $F \in \mathcal{F}_h^{\partial}$ with $F = \partial T \cap \partial \Omega$, $[\![\varphi ]\!]_F := \varphi|_{T}$ and $\{\! \{ \varphi \} \!\}_{\overline{\omega}} := \varphi |_{T}$. For any $v \in H^{1}(\mathbb{T}_h)^d$, we use the notation
\begin{align*}
\displaystyle
&[\![v \cdot n]\!] := [\![ v \cdot n ]\!]_F := v_1 \cdot n_F - v_2 \cdot n_F,  \  \{\! \{ v\} \! \}_{\omega} :=  \{\! \{ v \} \! \}_{\omega,F} := \omega_{T_1,F} v_1 + \omega_{T_2,F} v_2,\\
&[\![v]\!] :=  [\![ v]\!]_F := v_1 - v_2,
\end{align*}
for the jump of the normal component of $v$, the weighted average of $v$, and the junp {of} $v$. 


We define a broken gradient operator as follows. For $\varphi \in H^1(\mathbb{T}_h)^d$, the broken gradient $\nabla_h:H^1(\mathbb{T}_h)^d \to L^2(\Omega)^{d \times d}$ is defined by
\begin{align*}
\displaystyle
(\nabla_h \varphi)|_T &:= \nabla (\varphi|_T) \quad \forall T \in \mathbb{T}_h.
\end{align*}

For all $T \in \mathbb{T}_h$, we define the broken $H(\div;T)$ space as
\begin{align*}
\displaystyle
H(\div;\mathbb{T}_h) := \left \{ v \in L^2(\Omega)^d; \ v |_T \in H(\div;T) \ \forall T \in \mathbb{T}_h  \right\},
\end{align*}
and the broken divergence operator $\divh : H(\div;\mathbb{T}_h) \to L^2(\Omega)$ such that for all $v \in H(\div;\mathbb{T}_h)$,
\begin{align*}
\displaystyle
(\divh v)|_T := \div (v |_T) \quad \forall T \in \mathbb{T}_h.
\end{align*}

Let $T \in \mathbb{T}_h$. For any $k \in \mathbb{N}_0$, let $\mathbb{P}^k(T)$ be the space of polynomials with degree at most $k$ in $T$. For any $F \in \mathcal{F}_h$, we define the $L^2$-projection $\Pi_F^{0}: L^2(F) \to \mathbb{P}^{0}(F)$ by
\begin{align*}
\displaystyle
\int_F (\Pi_F^{0} \varphi - \varphi)   ds = 0 \quad \forall \varphi \in L^2(F).
\end{align*}

\subsection{Penalty parameters and energy norms}
Deriving an appropriate penalty term is essential in discontinuous Galerkin methods (dG) on anisotropic meshes. The use of weighted averages gives robust dG schemes for various problems; see \cite{DonGeo22,PieErn12}. For any $v \in H^{1}(\mathbb{T}_h)^d$ and $\varphi \in H^{1}(\mathbb{T}_h)$,
\begin{align*}
\displaystyle
[\![ (v \varphi) \cdot n ]\!]_F
&=  \{\! \{ v \} \! \}_{\omega,F} \cdot n_F [\! [ \varphi ]\!]_F + [\![ v \cdot n ]\!]_F \{\! \{ \varphi\} \! \}_{\overline{\omega},F}.
\end{align*}
For example, if $u \in H^1_0(\Omega) \cap W^{2,1}(\Omega)$, setting $v := - \nabla u$, we have $ [\![ v \cdot n ]\!]_F = 0$ for all $F \in \mathcal{F}_h^i$, see \cite[Lemma 4.3]{PieErn12}. Using the trace (Lemma \ref{lem=trace}) and the H\"older inequalities, the weighted average gives the following estimate for the term $ \{\! \{ v \}\! \}_{\omega,F} \cdot n_F [\![ \varphi ]\!]_F$.
\begin{align}
\displaystyle
&\int_{F} \left|   \{\! \{ v \} \! \}_{\omega,F} \cdot n_F [\! [ \varphi ]\!]_F  \right| ds \nonumber\\
& \leq c \left(  \omega_{T_1,F}  \| v|_{T_1} \|_{H^1(T_1)^d}  \sqrt{ \ell_{T_1,F}}^{- 1} +  \omega_{T_2,F} \| v|_{T_2} \|_{H^1(T_2)^d}  \sqrt{ \ell_{T_2,F}}^{- 1}  \right)  \left \|  [\![ \varphi ]\!] \right\|_{L^2(F)} \nonumber\\
& \leq c  \left( \| v|_{T_1} \|_{H^1(T_1)^d}^2 + \| v|_{T_2} \|_{H^1(T_2)^d}^2 \right)^{\frac{1}{2}} \left( \omega_{T_1,F}^2 { \ell_{T_1,F}}^{-1} + \omega_{T_2,F}^2 {\ell_{T_2,F}}^{-1} \right)^{\frac{1}{2}} \left \|  [\![ \varphi ]\!] \right\|_{L^2(F)} \nonumber\\
&= c h^{\beta}  \left( \| v|_{T_1} \|_{H^1(T_1)^d}^2 + \| v|_{T_2} \|_{H^1(T_2)^d}^2 \right)^{\frac{1}{2}} \nonumber\\
&\quad \times \left( h^{- 2 \beta} \omega_{T_1,F}^2 { \ell_{T_1,F}}^{-1} + h^{- 2 \beta}\omega_{T_2,F}^2 {\ell_{T_2,F}}^{-1} \right)^{\frac{1}{2}} \left \|  [\![ \varphi ]\!] \right\|_{L^2(F)}, {\label{new=2=9}}
\end{align}
{where $\ell_{T_1,F}$ and $\ell_{T_2,F}$ are defined in the inequality \eqref{trace}. The weights $\omega_{T_1,F}$, $\omega_{T_2,F}$ and $\beta$ are nonnegative real numbers chosen latter on.} The associated penalty parameter when $ \omega_{T_1,F} =  \omega_{T_2,F} = \frac{1}{2}$ becomes 
\begin{align}
\displaystyle
\frac{1}{4}  h^{- 2 \beta} \left( \frac{1}{\ell_{T_1,F}} + \frac{1}{\ell_{T_2,F}} \right), \label{ave_pena}
\end{align}
also see \cite{KasTsu21}. Furthermore, \begin{align*}
\displaystyle
\frac{1}{4} h^{- 2 \beta} \max \left\{ \frac{1}{\ell_{T_1,F}}  , \frac{1}{\ell_{T_2,F}}  \right\} &\leq \frac{1}{4} h^{- 2 \beta} \left( \frac{1}{\ell_{T_1,F}} + \frac{1}{\ell_{T_2,F}} \right) \\
&\leq \frac{1}{2} h^{- 2 \beta}  \max \left\{ \frac{1}{\ell_{T_1,F}}  , \frac{1}{\ell_{T_2,F}}  \right\}.
\end{align*}
Another choice for the weighted parameters is such that for $F \in \mathcal{F}_h^i$ with $F = T_{1} \cap T_{2}$, $T_{1},T_{2} \in \mathbb{T}_h$,
\begin{align}
\displaystyle
 \omega_{T_i,F} :=  \frac{\sqrt{\ell_{T_i,F}}}{\sqrt{\ell_{T_1,F}} + \sqrt{\ell_{T_2,F}}}, \quad i=1,2. \label{weight}
\end{align}
Then, the associated penalty parameter is defined as
\begin{align}
\displaystyle
 \frac{2 h^{- 2 \beta}}{(\sqrt{ \ell_{T_1,F}} + \sqrt{\ell_{T_2,F}})^2}. \label{weight_pena}
\end{align}
This quantity is the special case of the parameter proposed in \cite{DonGeo22}. 
Given that
\begin{align*}
\displaystyle
\frac{2 h^{- 2 \beta}}{(\sqrt{ \ell_{T_1,F}} + \sqrt{\ell_{T_2,F}})^2}
&\leq \frac{2 h^{- 2 \beta}}{( 2 \min \{ \sqrt{ \ell_{T_1,F}} , \sqrt{ \ell_{T_2,F}} \})^2} = \frac{ h^{- 2 \beta}}{2 \min \{  \ell_{T_1,F} , \ell_{T_2,F} \}} \\
&\leq \frac{1}{2} h^{- 2 \beta} \left( \frac{1}{\ell_{T_1,F}} + \frac{1}{\ell_{T_2,F}} \right),
\end{align*}
the quantity \eqref{weight_pena} {makes \eqref{new=2=9} a shaper bound} on anisotropic meshes than the parameter \eqref{ave_pena} for sufficiently small $h$. See also Remark \ref{rem1}.

Therefore, we use the type \eqref{weight_pena} of the parameter. For any $F \in \mathcal{F}_h$, we set $h_F := \diam(F)$. Let $F \in \mathcal{F}_h^i$ with $F = T_1 \cap T_2$, $T_1,T_2 \in \mathbb{T}_h$ be an interior face and $F \in \mathcal{F}_h^{\partial}$ with $F = \partial T_{\partial} \cap \partial \Omega$, $T_{\partial} \in \mathbb{T}_h$ a boundary face. A new penalty parameter $\kappa_F$ for the WOPSIP method is defined as follows using \eqref{weight_pena} with $\beta =1$. 
\begin{align}
\displaystyle
\kappa_F :=
\begin{cases}
\displaystyle
h^{-2} \left( \sqrt{ \ell_{T_1,F}} + \sqrt{\ell_{T_2,F}} \right)^{-2} \quad \text{if $F \in \mathcal{F}_h^i$},\\
\displaystyle
h^{-2} \ell_{T_{\partial},F}^{-1} \quad \text{if $F \in \mathcal{F}_h^{\partial}$}.
\end{cases} \label{penalty0}
\end{align}
For the proof of the discrete Poincar\'e inequality ({Lemma \ref{lem6}}), we {will} use the following parameter.
\begin{align}
\displaystyle
\kappa_{F*} :=
\begin{cases}
\displaystyle
\left( \sqrt{ \ell_{T_1,F}} + \sqrt{\ell_{T_2,F}} \right)^{-2} \quad \text{if $F \in \mathcal{F}_h^i$},\\
\displaystyle
\ell_{T_{\partial},F}^{-1} \quad \text{if $F \in \mathcal{F}_h^{\partial}$}.
\end{cases} \label{penalty1}
\end{align}

\begin{rem} \label{rem1}
We set $s := \left( \frac{1}{2} \right)^{\varepsilon}$ for $\varepsilon \in \mathbb{R}$ and $\varepsilon \> 1$. Let $T_1$ be the triangle with vertices $(0,0)^T$, $(s,0)^T$, and $(s,1)^T$, and let $T_2$ be the triangle with vertices $(s,0)^T$, $(1,0)^T$, and $(s,1)^T$. Then, we have $F = \{ s\} \times (0,1)$, $\ell_{T_1,F} = s$, $\ell_{T_2,F} = 1-s$, ${h = \sqrt{(1-s)^2 + 1}}$, and
\begin{align*}
\displaystyle
\kappa_F &= {h^{-2}} \left( \sqrt{ \ell_{T_1,F}} + \sqrt{ \ell_{T_2,F} } \right)^{-2} \to {\frac{1}{2}} \quad \text{as $\varepsilon \to \infty$},\\
\kappa_{F*} &= \left( \sqrt{ \ell_{T_1,F}} + \sqrt{ \ell_{T_2,F} } \right)^{-2} \to 1 \quad \text{as $\varepsilon \to \infty$}.
\end{align*}
However, the quantity given in \eqref{ave_pena} with $\beta = 0$ diverges as follows.
\begin{align*}
\displaystyle
\frac{1}{4} \left( \frac{1}{\ell_{T_1,F}} + \frac{1}{\ell_{T_2,F}} \right)
\to \infty \quad \text{as $\varepsilon \to \infty$}.
\end{align*}
\end{rem}

\begin{rem} \label{rem=new3}
We set $s := \left( \frac{1}{2} \right)^{\varepsilon}$ for $\varepsilon \in \mathbb{R}$ and $\varepsilon \> 1$. Let $T_3$ be the triangle with vertices $(0,0)^T$, $(s,0)^T$, and $(s,1)^T$. Let $T_4$ be the triangle with vertices $(s,0)^T$, $(2s,0)^T$, and $(s,1)^T$. Then, we have $F = \{ s\} \times (0,1)$, $\ell_{T_3,F} = s$, $\ell_{T_4,F} = s$, ${h = \sqrt{s^2 + 1}}$, and
\begin{align*}
\displaystyle
\kappa_F &= {h^{-2}} \left( \sqrt{ \ell_{T_3,F}} + \sqrt{ \ell_{T_4,F} } \right)^{-2} \to \infty \quad \text{as $\varepsilon \to \infty$}, \\
\kappa_{F*} &= \left( \sqrt{ \ell_{T_1,F}} + \sqrt{ \ell_{T_2,F} } \right)^{-2} \to \infty \quad \text{as $\varepsilon \to \infty$}.
\end{align*}
This may not be avoided by triangulation. Therefore, the use of anisotropic meshes may {cause} an ill-conditioned linear system.
\end{rem}

\begin{rem}
To overcome this difficulty of Remark \ref{rem=new3}, we consider the following case. We set $s := \left( \frac{1}{2} \right)^{\varepsilon}$ for $\varepsilon \in \mathbb{R}$ and $\varepsilon \> 1$. Let $T_5$ be the rectangle with vertices $(0,0)^T$, $(s,0)^T$, $(s,1)^T$, and $(0,1)^T$. Let $T_6$ be the triangle with vertices $(s,0)^T$, $(1,0)^T$, and $(s,1)^T$. As in Remark \ref{rem1}, we have $F = \{ s\} \times (0,1)$, $\ell_{T_5,F} = s$, $\ell_{T_6,F} = 1-s$, and $\kappa_F,\kappa_{F*} \to 1$ as $\varepsilon \to \infty$. Therefore, the parameters may not be as large if quadrangles are used for elements adjacent to boundaries.

\end{rem}

We define the following norms for any $v \in H^1(\mathbb{T}_h)$.
\begin{align*}
\displaystyle
| v |_{h} := \left( | v |_{H^1(\mathbb{T}_h)}^2 + | v |_{J}^2 \right)^{\frac{1}{2}}
\end{align*}
with the jump seminorm
\begin{align*}
\displaystyle
 | v |_{J} := \left( \sum_{F \in \mathcal{F}_h} \kappa_F \| \Pi_F^{0} [\![ v ]\!] \|_{L^2(F)}^2  \right)^{\frac{1}{2}}
\end{align*}
and $\kappa_F$ defined as in \eqref{penalty0};
\begin{align*}
\displaystyle
| v |_{h*} := \left( | v |_{H^1(\mathbb{T}_h)}^2 + | v |_{J*}^2 \right)^{\frac{1}{2}}
\end{align*}
with
\begin{align*}
\displaystyle
 | v |_{J*} := \left( \sum_{F \in \mathcal{F}_h} \kappa_{F*} \| \Pi_F^{0} [\![ v ]\!] \|_{L^2(F)}^2  \right)^{\frac{1}{2}}
\end{align*}
and $\kappa_{F*}$ defined as in \eqref{penalty1}. For any $v \in H^1(\mathbb{T}_h)$, $ | v |_{J*} \leq  | v |_{J}$ for $h \leq 1$. The norm $| \cdot |_{h*}$ {will be used} to prove the discrete Poincar\'e inequality (Lemma \ref{lem5}).

\subsection{Edge characterisation on a simplex, a geometric parameter, and a condition} \label{element=cond}

\begin{Cond}[Case in which $d=2$] \label{cond1}
Let ${T} \in \mathbb{T}_h$ with the vertices ${p}_i$ ($i=1,\ldots,3$). We assume that $\overline{{p}_2 {p}_3}$ is the longest edge of ${T}$; i.e., $ h_{{T}} := |{p}_2 - {p}_ 3|$. We set  $h_1 = |{p}_1 - {p}_2|$ and $h_2 = |{p}_1 - {p}_3|$. We then assume that $h_2 \leq h_1$. Note that ${h_1 \approx h_T}$. 
\end{Cond}

\begin{Cond}[Case in which $d=3$] \label{cond2}
Let ${T} \in \mathbb{T}_h$ with the vertices ${p}_i$ ($i=1,\ldots,4$). Let ${L}_i$ ($1 \leq i \leq 6$) be the edges of ${T}$. We denote by ${L}_{\min}$  the edge of ${T}$ with the minimum length; i.e., $|{L}_{\min}| = \min_{1 \leq i \leq 6} |{L}_i|$. We set $h_2 := |{L}_{\min}|$ and assume that 
\begin{align*}
\displaystyle
&\text{the endpoints of ${L}_{\min}$ are either $\{ {p}_1 , {p}_3\}$ or $\{ {p}_2 , {p}_3\}$}.
\end{align*}
Among the four edges that share an endpoint with ${L}_{\min}$, we take the longest edge ${L}^{({\min})}_{\max}$. Let ${p}_1$ and ${p}_2$ be the endpoints of edge ${L}^{({\min})}_{\max}$. We thus have that
\begin{align*}
\displaystyle
h_1 = |{L}^{(\min)}_{\max}| = |{p}_1 - {p}_2|.
\end{align*}
We consider cutting $\mathbb{R}^3$ with the plane that contains the midpoint of edge ${L}^{(\min)}_{\max}$ and is perpendicular to the vector ${p}_1 - {p}_2$. Thus, we have two cases: 
\begin{description}
  \item[(Type \roman{sone})] ${p}_3$ and ${p}_4$  belong to the same half-space;
  \item[(Type \roman{stwo})] ${p}_3$ and ${p}_4$  belong to different half-spaces.
\end{description}
In each case, we set
\begin{description}
  \item[(Type \roman{sone})] ${p}_1$ and ${p}_3$ as the endpoints of ${L}_{\min}$, that is, $h_2 =  |{p}_1 - {p}_3| $;
  \item[(Type \roman{stwo})] ${p}_2$ and ${p}_3$ as the endpoints of ${L}_{\min}$, that is, $h_2 =  |{p}_2 - {p}_3| $.
\end{description}
Finally, we set $h_3 = |{p}_1 - {p}_4|$. Note that we implicitly assume that ${p}_1$ and ${p}_4$ belong to the same half-space. In addition, note that ${h_1 \approx h_T}$. 
\end{Cond}

We define vectors ${r}_n \in \mathbb{R}^d$, $n=1,\ldots,d$ as follows. If $d=2$,
\begin{align*}
\displaystyle
{r}_1 := \frac{p_2 - p_1}{|p_2 - p_1|}, \quad {r}_2 := \frac{p_3 - p_1}{|p_3 - p_1|},
\end{align*}
and if $d=3$,
\begin{align*}
\displaystyle
&{r}_1 := \frac{p_2 - p_1}{|p_2 - p_1|}, \quad {r}_3 := \frac{p_4 - p_1}{|p_4 - p_1|}, \quad
\begin{cases}
\displaystyle
{r}_2 := \frac{p_3 - p_1}{|p_3 - p_1|}, \quad \text{for case (\roman{sone})}, \\
\displaystyle
{r}_2 := \frac{p_3 - p_2}{|p_3 - p_2|} \quad \text{for case (\roman{stwo})}.
\end{cases}
\end{align*}
For a sufficiently smooth function $\varphi$ and vector function $v := (v_{1},\ldots,v_{d})^T$, we define the directional derivative as, for $i \in \{ 1 : d \}$,
\begin{align*}
\displaystyle
\frac{\partial \varphi}{\partial {r_i}} &:= ( {r}_i \cdot  \nabla_{x} ) \varphi = \sum_{i_0=1}^d ({r}_i)_{i_0} \frac{\partial \varphi}{\partial x_{i_0}^{}}, \\
\frac{\partial v}{\partial r_i} &:= \left(\frac{\partial v_{1}}{\partial r_i}, \ldots, \frac{\partial v_{d}}{\partial r_i} \right)^T 
= ( ({r}_i  \cdot \nabla_{x}) v_{1}, \ldots, ({r}_i  \cdot \nabla_{x} ) v_{d} )^T.
\end{align*}
For a multi-index $\beta = (\beta_1,\ldots,\beta_d) \in \mathbb{N}_0^d$, we use the notation
\begin{align*}
\displaystyle
\partial^{\beta}_{r} \varphi := \frac{\partial^{|\beta|} \varphi}{\partial r_1^{\beta_1} \ldots \partial r_d^{\beta_d}}, \quad h^{\beta} :=  h_{1}^{\beta_1} \cdots h_{d}^{\beta_d}.
\end{align*}

We proposed a geometric parameter $H_{T}$ in a prior work \cite{IshKobTsu21a}.
 \begin{defi} \label{defi1}
 The parameter $H_{{T}}$ is defined as
\begin{align*}
\displaystyle
H_{{T}} := \frac{\prod_{i=1}^d h_i}{|{T}|_d} h_{{T}}.
\end{align*}
\end{defi}
We introduce the geometric condition proposed in \cite{IshKobTsu21a}, which is equivalent to the maximum-angle condition \cite{IshKobSuzTsu21d}.

\begin{assume} \label{ass1}
A family of meshes $\{ \mathbb{T}_h\}$ has a semi-regular property if there exists $\gamma_0 \> 0$ such that
\begin{align}
\displaystyle
\frac{H_{T}}{h_{T}} \leq \gamma_0 \quad \forall \mathbb{T}_h \in \{ \mathbb{T}_h \}, \quad \forall T \in \mathbb{T}_h. \label{NewGeo}
\end{align}
\end{assume}

\subsection{Affine mappings and Piola transformations}
In anisotropic interpolation errors on anisotropic meshes, we follow a strategy proposed in several prior works \cite{Ish21,IshKobTsu21a,IshKobTsu21c}. Let $T \in \mathbb{T}_h$ with Condition \ref{cond1} when $d=2$ or Condition \ref{cond2} when $d=3$. We define an affine mapping $\Phi: \widehat{T} \to T$ as
\begin{align*}
\displaystyle
\Phi: \widehat{T} \ni \hat{x} \mapsto x := \Phi(\hat{x}) := {A}_{T} \hat{x} + b_{T} \in T,
\end{align*}
where ${A}_{{T}} \in \mathbb{R}^{d \times d}$ is an invertible matrix and $b_{T} \in \mathbb{R}^d$. See \cite[Section 2]{IshKobTsu21c}.

The Piola transformation $\Psi : L^1(\widehat{T})^d \to L^1({T})^d$ is defined as
\begin{align*}
\displaystyle
\Psi :  L^1(\widehat{T})^d  &\to  L^1({T})^d \\
\hat{v} &\mapsto v(x) :=  \Psi(\hat{v})(x) = \frac{1}{\det(A_T)} A_T \hat{v}(\hat{x}).
\end{align*}

\subsection{Finite element spaces and anisotropic interpolation error estimates}
For $k \in \mathbb{N}_0$, $\mathbb{P}^k(T)$ is spanned by the restriction to $T$ of polynomials in $\mathbb{P}^k$ where  $\mathbb{P}^k$ denotes the space of polynomials with degree at most $k$. Let $Ne$ be the number of elements included in the mesh $\mathbb{T}_h$. Thus, we write $\mathbb{T}_h = \{ T_j\}_{j=1}^{Ne}$.

\subsubsection{Discontinuous space and the $L^2$-orthogonal projection}
For $T_j \in \mathbb{T}_h$, $j \in \{ 1, \ldots , Ne \}$, let $\Pi_{\textcolor{blue}{T_j}}^0 : L^2(T_j) \to \mathbb{P}^0$ be the $L^2$-orthogonal projection defined as
\begin{align*}
\displaystyle
\Pi_{T_j}^0 \varphi := \frac{1}{|T_j|_d} \int_{\textcolor{blue}{T_j}} \varphi dx \quad \forall \varphi \in L^2(T_j).
\end{align*}

The following theorem gives an anisotropic error estimate of the projection $\Pi_{T_j}^0$. Obtaining this estimate is not a novel concept \cite{AcoDur99}. However, the settings for meshes used in previous works differ slightly from our settings in Section \ref{element=cond}. In our theory, the same estimate is obtained.

\begin{thr} \label{thr1}
For any $\hat{\varphi} \in H^{1}(\widehat{T})$ with ${\varphi} := \hat{\varphi} \circ {\Phi}^{-1}$,
\begin{align}
\displaystyle
\| \Pi_{T_j}^0 \varphi - \varphi \|_{L^2(T_j)} \leq c \sum_{i=1}^d h_i \left\| \frac{\partial \varphi}{\partial r_i} \right\|_{L^{2}(T_j)}. \label{L2ortho}
\end{align}
\end{thr}

\begin{pf*}
The scaling argument yields
\begin{align}
\displaystyle
\| \Pi_{T_j}^0 \varphi - \varphi \|_{L^2(T_j)}
\leq c |\det({A}_{T_j})|^{\frac{1}{2}} \| \Pi_{\widehat{T}}^0 \hat{\varphi} - \hat{\varphi} \|_{L^2(\widehat{T})}. \label{add=1}
\end{align}
For any $\hat{\eta} \in \mathbb{P}^{0}$, 
\begin{align*}
\displaystyle
\| \Pi_{\widehat{T}}^0 \hat{\varphi} - \hat{\varphi} \|_{L^2(\widehat{T})}
&\leq \| \Pi_{\widehat{T}}^0  (\hat{\varphi} - \hat{\eta}) \|_{L^2(\widehat{T})}  + \| \hat{\eta} - \hat{\varphi} \|_{L^2(\widehat{T})},
\end{align*}
because $\Pi_{\widehat{T}}^0 \hat{\eta} = \hat{\eta}$. The stability of the $L^2$-orthogonal projection yields
\begin{align*}
\displaystyle
\| \Pi_{\widehat{T}}^0  (\hat{\varphi} - \hat{\eta}) \|_{L^2(\widehat{T})} 
\leq c \| \hat{\varphi} - \hat{\eta} \|_{L^2(\widehat{T})}.
\end{align*}
Thus, 
\begin{align}
\displaystyle
\| \Pi_{\widehat{T}}^0 \hat{\varphi} - \hat{\varphi} \|_{L^2(\widehat{T})}
\leq c  \inf_{\hat{\eta} \in  \mathbb{P}^{0}} \| \hat{\varphi} - \hat{\eta} \|_{L^2(\widehat{T})}. \label{chap827}
\end{align}
From the Bramble--Hilbert-type lemma (e.g., see \cite[Lemma 4.3.8]{BreSco08}), there exists a constant $\hat{\eta}_{\beta} \in \mathbb{P}^{0}$ such that for any $\hat{\varphi} \in H^{1}(\widehat{T})$,
\begin{align}
\displaystyle
\| \hat{\varphi} - \hat{\eta}_{\beta} \|_{L^{2}(\widehat{T})} \leq C^{BH}(\widehat{T}) |\hat{\varphi}|_{H^{1}(\widehat{T})}. \label{chap829}
\end{align}
Using the inequality in \cite[Lemma 6]{IshKobTsu21c} with $m=0$ and $p=2$, the inequality \eqref{chap829} is estimated as
\begin{align}
\displaystyle
\| \hat{\varphi} - \hat{\eta}_{\beta} \|_{L^{2}(\widehat{T})}
&\leq c |\hat{\varphi}|_{H^{1}(\widehat{T})}
\leq c |\det({{A}_{T_j}})|^{-\frac{1}{2}} \sum_{i=1}^d h_i \left\| \frac{\partial \varphi}{\partial r_i} \right\|_{L^{2}(T_j)}.  \label{chap8210}
\end{align}
From \eqref{add=1}, \eqref{chap827}, and \eqref{chap8210}, we have the target inequality \eqref{L2ortho}.
\qed
\end{pf*}

\subsubsection{Discontinuous CR finite element space, an associated interpolation operator, and a Stokes element}
We introduce a discontinuous CR finite element space and an associated interpolation operator, as well as a Stokes element.

Let the points $\{ P_{T_j,1}, \ldots, P_{T_j,d+1} \}$ be the vertices of the simplex $T_j \in \mathbb{T}_h$ for $j \in \{1, \ldots , Ne \}$. Let $F_{T_j,i}$ be the face of $T_j$ opposite $P_{T_j,i}$ for $i \in \{ 1, \ldots , d+1\}$. We set $P := \mathbb{P}^1$, and take a set ${\Sigma}_{T_j} := \{ {\chi}^{CR}_{T_j,i} \}_{1 \leq i \leq d+1}$ of linear forms with its components such that for any $p \in \mathbb{P}^1$.
\begin{align}
\displaystyle
{\chi}^{CR}_{T_j,i}({p}) := \frac{1}{| {F}_{T_j,i} |_{d-1}} \int_{{F}_{T_j,i}} {p} d{s} \quad \forall i \in \{ 1, \ldots,d+1 \}. \label{CR1}
\end{align}
For each $j \in \{1, \ldots ,Ne \}$, the triple $\{ T_j ,  \mathbb{P}^1 , \Sigma_{T_j} \}$ is a finite element. Using the barycentric coordinates $ \{{\lambda}_{T_j,i} \}_{i=1}^{d+1}: \mathbb{R}^d \to \mathbb{R}$ on the reference element, the nodal basis functions associated with the degrees of freedom by \eqref{CR1} are defined as
\begin{align}
\displaystyle
{\theta}^{CR}_{T_j,i}({x}) := d \left( \frac{1}{d} - {\lambda}_{T_j,i} ({x}) \right) \quad \forall i \in \{ 1, \ldots ,d+1 \}. \label{CR2}
\end{align}
For $j \in \{1, \ldots ,Ne \}$ and $i \in \{ 1, \ldots , d+1\}$, we define the function $\phi_{j(i)}$ as
\begin{align}
\displaystyle
\phi_{j(i)}(x) :=
\begin{cases}
\theta_{T_j,i}^{CR}(x), \quad \text{$x \in T_j$}, \\
0, \quad \text{$x \notin T_j$}.
\end{cases} \label{CR5}
\end{align}
We define a discontinuous finite element space as
\begin{align}
\displaystyle
X_{dc,h}^{CR} &:= \left\{ \sum_{j=1}^{Ne} \sum_{i=1}^{d+1} c_{j(i)} \phi_{j(i)}; \  c_{j(i)} \in \mathbb{R}, \ \forall i,j \right\} \subset P_{dc,h}^1. \label{CR6}
\end{align}
For $s \in \mathbb{N}_0$, we define a discontinuous finite element space as
\begin{align*}
\displaystyle
P_{dc,h}^{s} &:= \left\{ p_h \in L^2(\Omega); \ p_h|_{T_j} \circ {\Phi} \in \mathbb{P}^{s}(\widehat{T}) \quad \forall T_j \in \mathbb{T}_h, \quad j \in \{ 1, \ldots ,Ne\} \right\}.
\end{align*}
We define a pair of the standard dG spaces $(V_{h}^{m_1},Q_h^{m_2} )$ as
\begin{align*}
\displaystyle
V_{dc,h}^{m_1} := ( P_{dc,h}^{m_1})^d, \quad Q_{h}^{m_2} := P_{dc,h}^{m_2} \cap Q,
\end{align*}
for $m_1 \in \mathbb{N}$ and $m_2 \in \mathbb{N}_0$. We use the discrete space $(V_{dc,h}^{m_1},Q_h^{m_2} )$ in Section \ref{numerial=comp}. Let $(V_{dc,h}^{CR},Q_h^{0})$ be a pair of discontinuous finite element spaces defined by
\begin{align}
\displaystyle
V_{dc,h}^{CR} := (X_{dc,h}^{CR})^d, \quad Q_h^{0} = P_{dc,h}^{0} \cap Q, \label{CR7}
\end{align}
with norms
\begin{align*}
\displaystyle
|v_h|_{V_{dc,h}^{CR}} := \left ( \sum_{i=1}^d |v_{h,i}|^2_{h} \right)^{\frac{1}{2}}, \quad \| q_h \|_{Q_h^{0}} := \| q_h \|_{L^2(\Omega)}
\end{align*}
for any $v_h = (v_{h,1},\ldots,v_{h,d})^T  \in V_{dc,h}^{CR}$ and $q_h \in Q_h^{0}$. We {will} use the discrete space $(V_{dc,h}^{CR},Q_h^{0})$ in the WOPSIP method (Section \ref{sec=WOPSIP}). {Furthermore, we define a norm as
\begin{align*}
\displaystyle
|v_h|_{V_{dc,h*}^{CR}} := \left ( \sum_{i=1}^d |v_{h,i}|^2_{h*} \right)^{\frac{1}{2}}.
\end{align*}
For any $v_h \in V_{dc,h}^{CR}$, $|v_h|_{V_{dc,h*}^{CR}} \leq |v_h|_{V_{dc,h}^{CR}}$ for $h \leq 1$. The norm $|\cdot|_{V_{dc,h*}^{CR}}$ will be used in Theorem \ref{thr7}.
}

Let $I_{T_j}^{CR} : H^{1}(T_j) \to \mathbb{P}^1(T_j)$ be the CR interpolation operator such that for any $\varphi \in H^{1}(T_j)$, 
\begin{align}
\displaystyle
I_{T_j}^{CR}: H^{1}(T_j) \ni \varphi \mapsto I_{T_j}^{CR} \varphi := \sum_{i=1}^{d+1} \left(  \frac{1}{| {F}_{T_j,i} |_{d-1}} \int_{{F}_{T_j,i}} {\varphi} d{s} \right) \theta_{T_j,i}^{CR} \in \mathbb{P}^1(T_j). \label{CR3}
\end{align}

We then present estimates of the anisotropic CR interpolation error. Obtaining the estimate is not a novel concept \cite{ApeNicSch01}. However, the proof provided here differs from those given in prior works. We here present a proof using the error estimate of the $L^2$-projection $\Pi_h^0$ in Theorem \ref{thr1}.

\begin{thr} \label{thr2}
For $j \in \{ 1, \ldots , Ne\}$,
\begin{align}
\displaystyle
|I_{T_j}^{CR} \varphi - \varphi |_{H^{1}({T}_j)} &\leq c \sum_{i=1}^d {h}_i \left\| \frac{\partial }{\partial r_i} \nabla \varphi \right \|_{L^2(T_j)^d} \quad \forall {\varphi} \in H^{2}({T}_j). \label{CR4}
\end{align}
\end{thr}

\begin{pf*}
Let ${\varphi} \in H^{2}({T}_j)$, $j \in \{ 1, \ldots , Ne\}$. Because $I_{T_j}^{CR} \varphi  \in \mathbb{P}^1$, Green's formula and the definition of the CR interpolation imply that
\begin{align*}
\displaystyle
\frac{\partial }{\partial {x}_k} (I_{T_j}^{CR} \varphi)
&= \frac{1}{|T_j|_d} \int_{T_j}  \frac{\partial }{\partial {x}_k} (I_{T_j}^{CR} \varphi) dx = \frac{1}{|T_j|_d} \sum_{i=1}^{d+1} n_{{T}_j}^{(k)} \int_{F_i} I_{T_j}^{CR} \varphi ds \\
&= \frac{1}{|T_j|_d} \sum_{i=1}^{d+1} n_{{T}_j}^{(k)} \int_{F_i} \varphi ds 
= \frac{1}{|T_j|_d} \int_{T_j}  \frac{\partial \varphi}{\partial {x}_k} dx = \Pi_{T_j}^0 \left(  \frac{\partial \varphi}{\partial {x}_k} \right),
\end{align*}
for $k \in \{ 1, \ldots , d\}$, where $ n_{{T}_j}^{(k)}$ denotes the $k$th component of the outer unit normal vector $n_{{T}_j}$. The inequality \eqref{L2ortho} yields
\begin{align*}
\displaystyle
| I_{T_j}^{CR} \varphi - \varphi |_{H^{1}(T_j)}^2
&= \sum_{k=1}^d \left \|  \frac{\partial}{\partial x_k} ( I_{T_j}^{CR} \varphi - \varphi)  \right\|^2_{L^2(T_j)} \\
&\hspace{-1cm} = \sum_{k=1}^d \left \| \Pi_{T_j}^0 \left(  \frac{\partial \varphi}{\partial {x}_k} \right) - \left(  \frac{\partial \varphi}{\partial {x}_k} \right) \right\|^2_{L^2(T_j)} 
\leq c \sum_{i,k=1}^d h_i^2 \left \| \frac{\partial^2 \varphi}{ \partial r_i \partial x_k} \right \|_{L^{2}(T_j)}^2,
\end{align*}
which leads to \eqref{CR4} together with the Jensen's inequality.
\qed
\end{pf*}
The vector-valued local interpolation operator 
\begin{align}
\displaystyle
\mathcal{I}_{T_j}^{CR}: H^{1}(T_j)^d \to \mathbb{P}^1(T_j)^d \label{vecCRint}
\end{align}
is defined component-wise, that is,
\begin{align*}
\displaystyle
\mathcal{I}_{T_j}^{CR} v := ( I_{T_j}^{CR} v_1,\ldots,  I_{T_j}^{CR} v_d)^T \quad \forall v = (v_1,\ldots,v_d)^T \in H^{1}(T_j)^d.
\end{align*}
We define a global interpolation operator $\mathcal{I}_{h}^{CR}: H^{1}(\Omega)^d \to V_{dc,h}^{CR}$ as
\begin{align}
\displaystyle
(\mathcal{I}_{h}^{CR} v )|_{T_j} = \mathcal{I}_{T_j}^{CR} (v |_{T_j}), \quad j \in \{ 1, \ldots , Ne\}, \quad \forall v \in H^{1}(\Omega)^d.  \label{CR8}
\end{align}	

\subsubsection{Discontinuous RT finite element space and an associated interpolation operator} \label{RTsp}
We introduce a discontinuous RT finite element space and an associated interpolation operator.

For $T_j \in \mathbb{T}_h$, $j \in \{ 1, \ldots ,Ne \}$, we define the local RT polynomial space as follows.
\begin{align}
\displaystyle
\mathbb{RT}^0(T_j) := \mathbb{P}^0(T_j)^d + x \mathbb{P}^0(T_j), \quad x \in \mathbb{R}^d. \label{RT1}
\end{align}
For $p \in \mathbb{RT}^0(T_j)$, the local degrees of freedom are defined as
\begin{align}
\displaystyle
{\chi}^{RT}_{T_j,i}({p}) := \int_{{F}_{T_j,i}} {p} \cdot n_{T_j,i} d{s} \quad \forall i \in \{ 1, \ldots ,d+1 \}, \label{RT2}
\end{align}
where $n_{T_j,i}$ is the outward normal to ${F}_{T_j,i}$. Setting $\Sigma_{T_j}^{RT} := \{ {\chi}^{RT}_{T_j,i} \}_{1 \leq i \leq d+1}$, the triple $\{ T_j ,  \mathbb{RT}^0 , \Sigma_{T_j}^{RT} \}$ a finite element. The local shape functions are 
\begin{align}
\displaystyle
\theta_{T_j,i}^{RT}(x) := \frac{\iota_{{F_{T_j,i}},{T}_j}}{d |T_j|_d} (x - P_{T_j,i}) \quad \forall i \in \{ 1, \ldots , d+1 \}, \label{RT3}
\end{align}
where $\iota_{{F_{T_j,i}},{T}_j} := 1$ if ${n}_{T_{j,i}}$ points outwards, and $ - 1$ otherwise \cite{ErnGue21a}. We define a discontinuous RT finite element space as follows.
\begin{align}
\displaystyle
V_{dc,h}^{RT} &:= \{ v_h \in L^1(\Omega)^d : \ v_h|_{T} \in \mathbb{RT}^0(T) \quad \forall T \in \mathbb{T}_h \}. \label{RT7}
\end{align}

Let $\mathcal{I}_{T_j}^{RT}: H^{1}(T_j)^d \to \mathbb{RT}^0(T_j)$ be the RT interpolation operator such that for any $v \in H^{1}(T_j)^d$,
\begin{align}
\displaystyle
\mathcal{I}_{T_j}^{RT}: H^{1}(T_j)^d \ni v \mapsto \mathcal{I}_{T_j}^{RT} v := \sum_{i=1}^{d+1} \left(  \int_{{F}_{T_j,i}} {v} \cdot n_{T_j,i} d{s} \right) \theta_{T_j,i}^{RT} \in \mathbb{RT}^0(T_j). \label{RT4}
\end{align}
The following two theorems are divided into the element of (Type \roman{sone}) or the element of (Type \roman{stwo}) in Section \ref{element=cond} when $d=3$.

\begin{thr} \label{thr3}
Let $T_j$ be the element with Conditions \ref{cond1} or \ref{cond2} and satisfy (Type \roman{sone}) in Section \ref{element=cond} when $d=3$. For any $\hat{v} \in H^{1}(\widehat{T})^d$ with ${v} = ({v}_1,\ldots,{v}_d)^T := {\Psi} \hat{v}$ and $j \in \{ 1, \ldots ,Ne\}$,
\begin{align}
\displaystyle
\| \mathcal{I}_{T_j}^{RT} v - v \|_{L^2(T_j)^d} 
&\leq  c \left( \frac{H_{T_j}}{h_{T_j}} \sum_{i=1}^d h_i \left \|  \frac{\partial v}{\partial r_i} \right \|_{L^2(T_j)^d} +  h_{T_j} \| \div {v} \|_{L^{2}({T}_j)} \right). \label{RT5}
\end{align}
\end{thr}

\begin{pf*}
A proof is provided in \cite[Theorem 2]{Ish21}.
\qed
\end{pf*}

\begin{thr} \label{thr4}
Let $d=3$. Let $T_j$ be an element with Condition \ref{cond2} that satisfies (Type \roman{stwo}) in Section \ref{element=cond}. For $\hat{v} \in H^{1}(\widehat{T})^3$ with ${v} = ({v}_1,v_2,{v}_3)^T := {\Psi} \hat{v}$ and $j \in \{ 1, \ldots ,Ne\}$,
\begin{align}
\displaystyle
&\| \mathcal{I}_{T_j}^{RT} v - v \|_{L^2(T_j)^3} 
\leq c \frac{H_{T_j}}{h_{T_j}} \Biggl(  h_{T_j} |v|_{H^1(T_j)^3} \Biggr). \label{RT6}
\end{align}
\end{thr}

\begin{pf*}
A proof is provided in \cite[Theorem 3]{Ish21}.
\qed
\end{pf*}

We define a global  interpolation operator $\mathcal{I}_{h}^{RT}: H^{1}(\Omega)^d \cup V_{dc,h}^{CR}  \to V_{dc,h}^{RT}$ by
\begin{align}
\displaystyle
(\mathcal{I}_{h}^{RT} v )|_{T_j} = \mathcal{I}_{T_j}^{RT} (v |_{T_j}), \quad j \in \{ 1, \ldots , Ne\}, \quad \forall v \in H^{1}(\Omega)^d \cup V_{dc,h}^{CR}.  \label{RT8}
\end{align}
We also define the global interpolation $\Pi_h^0$ to the space $P_{dc,h}^{0}$ as
\begin{align*}
\displaystyle
(\Pi_h^0 \varphi)|_{T_j} := \Pi_{\textcolor{blue}{T_j}}^0 (\varphi|_{T_j}) \ \forall T_j \in \mathbb{T}_h, \ j \in \{ 1, \ldots ,Ne\}, \ \forall \varphi \in L^2(\Omega).
\end{align*}
Between the RT interpolation $\mathcal{I}_{h}^{RT}$ and the $L^2$-projection $\Pi_h^0$, the following relation holds.
\begin{lem} \label{lem2}
For $j \in \{ 1, \ldots , Ne\}$,
\begin{align}
\displaystyle
\div (\mathcal{I}_{T_j}^{RT} v) = \Pi_{T_j}^0 (\div v) \quad \forall v \in H^{1}(T_j)^d. \label{RT9a}
\end{align}
By combining \eqref{RT9a}, for any $v \in H^1(\Omega)^d$
\begin{align}
\displaystyle
\div (\mathcal{I}_{h}^{RT} v) = \Pi_h^0 (\div v). \label{RT9b}
\end{align}
\end{lem}

\begin{pf*}
A proof is provided in \cite[Lemma 16.2]{ErnGue21a}.
\qed
\end{pf*}

\section{Discrete Poincar\'e inequality}
{In this section, we present the discrete Poincar\'e inequality on anisotropic meshes. In \cite{Bre03,PieErn12}, the inequality was proven under the shape regularity condition. However, it is difficult to derive the inequality on anisotropic meshes. Here, we present a proof using a dual problem. Therefore, we impose that $\Omega$ is convex to establish the aforementioned inequality.}

The following relation plays an important role in the discontinuous Galerkin finite element analysis on anisotropic meshes.

\begin{lem} \label{lem3}
For any $w \in H^1(\Omega)^d$ and $\psi_h \in P_{dc,h}^{1}$,
\begin{align}
\displaystyle
&\int_{\Omega} \left( \mathcal{I}_h^{RT} w \cdot \nabla_h \psi_{h} + \div \mathcal{I}_h^{RT} w  \psi_{h} \right) dx\notag\\
&\quad =  \sum_{F \in \mathcal{F}_h^i} \int_{F}  \{ \! \{ w \} \!\}_{\omega,F} \cdot n_F \Pi_F^0 [\![ \psi_{h} ]\!]_F ds + \sum_{F \in \mathcal{F}_h^{\partial}} \int_{F} ( w \cdot n_F) \Pi_F^0 \psi_{h} ds.  \label{wop=3}
\end{align}
\end{lem}

\begin{pf*}
For any $w \in H^1(\Omega)^d$ and $\psi_h \in P_{dc,h}^{1}$, using Green’s formula and the fact $\mathcal{I}_h^{RT} w \cdot n_F \in \mathbb{P}^{0}(F)$ for any $F \in \mathcal{F}_h$, we derive
\begin{align*}
\displaystyle
&\int_{\Omega} \left( \mathcal{I}_h^{RT} w \cdot \nabla_h \psi_{h} + \div \mathcal{I}_h^{RT} w  \psi_{h} \right) dx
= \sum_{T \in \mathbb{T}_h} \int_{\partial T} (\mathcal{I}_h^{RT} w \cdot n_T) \psi_{h} ds  \notag \\
&\quad = \sum_{F \in \mathcal{F}_h^i} \int_{F} \left(  [\![ \mathcal{I}_h^{RT} w \cdot n]\!]_F \{\! \{ \psi_{h} \} \! \}_{\overline{\omega},F} + \{\! \{ \mathcal{I}_h^{RT} w  \}\! \}_{\omega,F} \cdot n_F [\![ \psi_{h} ]\!]_F \right) ds  \notag \\
&\quad + \sum_{F \in \mathcal{F}_h^{\partial}} \int_{F} (\mathcal{I}_h^{RT} w \cdot n_F) \psi_{h} ds.
\end{align*}
{Recall that the weighted and skew-weighted averages $ \{\! \{ \cdot \} \! \}_{{\omega},F}$ and  $ \{\! \{ \cdot \} \! \}_{\overline{\omega},F}$ were defined in Section \ref{regularmesh}}.

By the midpoint rule, we have
\begin{align*}
\displaystyle
[\![ \mathcal{I}_{h}^{RT} w \cdot n ]\!]_F (x_F)
&= \frac{1}{|F|_{d-1}} \int_F [\![ \mathcal{I}_{h}^{RT} w \cdot n ]\!]_F (x) ds \\
&= \frac{1}{|F|_{d-1}} \int_F [\![  w \cdot n ]\!]_F (x) ds = 0,
\end{align*}
which leads to
\begin{align}
\displaystyle
[\![ \mathcal{I}_{h}^{RT} w \cdot n ]\!]_F (x) = 0 \quad \forall x \in F. \label{wop=4}
\end{align}
Using \eqref{wop=4} and the properties of the projection $\Pi_F^0$ and $\mathcal{I}_h^{RT}$ yields
\begin{align*}
\displaystyle
&\int_{\Omega} \left( \mathcal{I}_h^{RT} w \cdot \nabla_h \psi_{h} + \div \mathcal{I}_h^{RT} w  \psi_{h} \right) dx \\
&\quad =  \sum_{F \in \mathcal{F}_h^i} \int_{F}  \{ \! \{ \mathcal{I}_h^{RT} w \} \!\}_{\omega,F} \cdot n_F \Pi_F^0 [\![ \psi_{h} ]\!]_F ds + \sum_{F \in \mathcal{F}_h^{\partial}} \int_{F} (\mathcal{I}_h^{RT} w \cdot n_F) \Pi_F^0 \psi_{h} ds \\
&\quad =  \sum_{F \in \mathcal{F}_h^i} \int_{F}  \{ \! \{  w \} \!\}_{\omega,F} \cdot n_F \Pi_F^0 [\![ \psi_{h} ]\!]_F ds + \sum_{F \in \mathcal{F}_h^{\partial}} \int_{F} ( w \cdot n_F) \Pi_F^0 \psi_{h} ds,
\end{align*}
which is the desired equality.
\qed
\end{pf*}

The right-hand terms in \eqref{wop=3} are estimated as follows.

\begin{lem} \label{lem4}
For any $w \in H^1(\Omega)^d$ and $\psi_h \in P_{dc,h}^{1}$,
\begin{align}
\displaystyle
&\left| \sum_{F \in \mathcal{F}_h^i} \int_{F}  \{ \! \{ w \} \!\}_{\omega,F} \cdot n_F \Pi_F^0 [\![ \psi_{h} ]\!]_F ds \right|  \notag\\
&\quad \leq c |\psi_{h}|_J \left(  h \| w \|_{L^2(\Omega)^d} + h^{\frac{3}{2}} \| w \|_{L^2(\Omega)^d}^{\frac{1}{2}} | w |_{H^1(\Omega)^d}^{\frac{1}{2}} \right), \label{wop=5} \\
&\left| \sum_{F \in \mathcal{F}_h^{\partial}} \int_{F} ( w \cdot n_F) \Pi_F^0 \psi_{h} ds \right| \notag \\
&\quad \leq c |\psi_{h}|_J \left(  h \| w \|_{L^2(\Omega)^d} + h^{\frac{3}{2}} \| w \|_{L^2(\Omega)^d}^{\frac{1}{2}} | w |_{H^1(\Omega)^d}^{\frac{1}{2}} \right). \label{wop=6}
\end{align}
\end{lem}

\begin{pf*}
Using the H\"older's inequality, the weighted average and the trace inequality \eqref{trace} yields
\begin{align*}
\displaystyle
&\int_{F} \left| \{ \! \{ w \} \!\}_{\omega,F} \cdot n_F \Pi_F^0 [\![ \psi_{h} ]\!]_F \right | ds \\
&\quad \leq c \Bigg \{ \left(  \| w_1 \|_{L^2(T_1)^d} + h_{T_1}^{\frac{1}{2}}  \| w_1 \|_{L^2(T_1)^d}^{\frac{1}{2}} | w_1 |_{H^1(T_1)^d}^{\frac{1}{2}} \right)^2 \\
&\quad \quad  + \left(  \| w_2 \|_{L^2(T_2)^d} + h_{T_2}^{\frac{1}{2}}  \| w_2 \|_{L^2(T_2)^d}^{\frac{1}{2}} | w_2 |_{H^1(T_2)^d}^{\frac{1}{2}}  \right)^{2} \Biggr\}^{\frac{1}{2}} \\
&\quad \quad \times \left( \omega_{T_1,F}^2 { \ell_{T_1,F}}^{-1} + \omega_{T_2,F}^2 {\ell_{T_2,F}}^{-1} \right)^{\frac{1}{2}} \| \Pi_F^{0} [\![ \psi_{h} ]\!] \|_{L^2(F)}.
\end{align*}
Using the Cauchy--Schwarz inequality,
\begin{align*}
\displaystyle
&\left| \sum_{F \in \mathcal{F}_h^i} \int_{F}  \{ \! \{ w \} \!\}_{\omega,F} \cdot n_F \Pi_F^0 [\![ \psi_{h} ]\!]_F ds \right| \\
&\quad \leq  c \sum_{F \in \mathcal{F}_h^i}  h^{-1} \left( \omega_{T_1,F}^2 { \ell_{T_1,F}}^{-1} + \omega_{T_2,F}^2 {\ell_{T_2,F}}^{-1} \right)^{\frac{1}{2}} \| \Pi_F^{0} [\![ \psi_{h} ]\!] \|_{L^2(F)} \\
&\quad \quad \times \sum_{T \in \mathbb{T}_F} h  \left(  \| w \|_{L^2(T)^d} + h_{T}^{\frac{1}{2}} \| w \|_{L^2(T)^d}^{\frac{1}{2}} | w |_{H^1(T)^d}^{\frac{1}{2}} \right) \\
&\quad \leq  c \left( \sum_{F \in \mathcal{F}_h^i}  h^{-2} \left( \omega_{T_1,F}^2 { \ell_{T_1,F}}^{-1} + \omega_{T_2,F}^2 {\ell_{T_2,F}}^{-1} \right) \| \Pi_F^{0} [\![ \psi_{h} ]\!] \|_{L^2(F)}^2 \right)^{\frac{1}{2}} \\
&\quad \quad \times \left ( \sum_{F \in \mathcal{F}_h^i}  \sum_{T \in \mathbb{T}_F} h^2  \left(  \| w \|_{L^2(T)^d} + h_{T}^{\frac{1}{2}} \| w \|_{L^2(T)^d}^{\frac{1}{2}} | w |_{H^1(T)^d}^{\frac{1}{2}} \right)^2 \right)^{\frac{1}{2}}
\end{align*}
which leads to the inequality \eqref{wop=5} together with the weight \eqref{weight} and the Cauchy--Schwarz and Jensen inequalities. {Here, $\mathbb{T}_F$ denotes the set of the simplices in $\mathbb{T}_h$ that share $F$ as a common face.}

By an analogous argument, the estimate \eqref{wop=6} holds.
\qed
\end{pf*}

\begin{lem} \label{lem5}
Let $h \leq 1$. Thus, for any $w \in H^1(\Omega)^d$ and $\psi_h \in P_{dc,h}^{1}$,
\begin{align}
\displaystyle
\left| \sum_{F \in \mathcal{F}_h^i} \int_{F}  \{ \! \{ w \} \!\}_{\omega,F} \cdot n_F \Pi_F^0 [\![ \psi_{h} ]\!]_F ds \right|
&\leq c |\psi_{h}|_{J*} \| w \|_{H^1(\Omega)^d}, \label{wop=7} \\
\left| \sum_{F \in \mathcal{F}_h^{\partial}} \int_{F} ( w \cdot n_F) \Pi_F^0 \psi_{h} ds \right|
&\leq  c |\psi_{h}|_{J*}  \| w \|_{H^1(\Omega)^d}. \label{wop=8}
\end{align}
\end{lem}

\begin{pf*}
By analogous proof with Lemma \ref{lem4}, we can obtain the target inequalities.
\qed
\end{pf*}

The following lemma provides a discrete Poincar\'e inequality. 

\begin{lem}[Discrete Poincar\'e inequality] \label{lem6}
Assume that $\Omega$ is convex. Let $\{ \mathbb{T}_h\}$ be a family of meshes with the semi-regular property (Assumption \ref{ass1}) and $h \leq 1$. Then, there exists a positive constant $C_{dc}^{P}$ independent of $h$ {but dependent on the maximum angle} such that
\begin{align}
\displaystyle
 \| \psi_{h} \|_{L^2(\Omega)} \leq C_{dc}^{P} | \psi_{h} |_{h*} \quad \forall \psi_h \in P_{dc,h}^{1}. \label{wop=9}
\end{align}
\end{lem}

\begin{pf*}
Let $\psi_h \in P_{dc,h}^{1}$. We consider the following problem. Find $z \in H^{2}(\Omega) \cap H_0^1(\Omega)$ such that
\begin{align*}
\displaystyle
- \varDelta z = \psi_{h} \quad \text{in $\Omega$}, \quad z = 0 \quad \text{on $\partial \Omega$}.
\end{align*}
We then have a priori estimates $|z|_{H^1(\Omega)} \leq C_P \| \psi_{h} \|_{L^2(\Omega)}$ and $ |z|_{H^2(\Omega)} \leq \| \psi_h \|_{L^2(\Omega)}$, where $C_P$ is the Poincar\'e constant.

We provide the following {equality} for analysis.
\begin{align*}
\displaystyle
 \int_{\Omega} ( \div \mathcal{I}_{h}^{RT}  (\nabla z) ) \psi_h {dx}
 &= \sum_{T \in \mathbb{T}_h} \int_{\partial T} n_T \cdot  \mathcal{I}_{h}^{RT} (\nabla z) \psi_h ds - \int_{\Omega} \mathcal{I}_{h}^{RT} (\nabla z) \cdot \nabla_h \psi_h dx \\
 &=  \int_{\Omega} ( \nabla z -  \mathcal{I}_{h}^{RT}  (\nabla z) ) \cdot \nabla_h \psi_h dx  - \int_{\Omega}  \nabla z \cdot \nabla_h \psi_h dx \\
 &\quad + \sum_{T \in \mathbb{T}_h} \int_{\partial T} \mathcal{I}_{h}^{RT}  (\nabla z) \cdot n_T \psi_h ds \\
 &= \int_{\Omega} ( \nabla z -  \mathcal{I}_{h}^{RT}  (\nabla z) ) \cdot \nabla_h \psi_h dx  - \int_{\Omega}  \nabla z \cdot \nabla_h \psi_h dx \\
 &\quad +  \sum_{F \in \mathcal{F}_h^i} \int_{F}  \{ \! \{  \nabla z \} \!\}_{\omega,F} \cdot n_F \Pi_F^0 [\![ \psi_{h} ]\!]_F ds \\
 &\quad + \sum_{F \in \mathcal{F}_h^{\partial}} \int_{F} (  \nabla z \cdot n_F) \Pi_F^0 \psi_{h} ds,
\end{align*}
where we apply the analogous argument with Lemma \ref{lem3} for the last equality. This {equality} yields
\begin{align*}
\displaystyle
&\| \psi_h \|^2_{L^2(\Omega)} = \int_{\Omega} \psi_h^2 dx = \int_{\Omega} - \varDelta z \psi_h dx = - \int_{\Omega} \div (\nabla z)  \psi_h dx \\
&\quad =  \int_{\Omega} ( \Pi_h^{0}  \div (\nabla z) - \div (\nabla z) ) \psi_h dx - \int_{\Omega} ( \Pi_h^{0}  \div (\nabla z) ) \psi_h dx\\
&\quad =  \int_{\Omega} ( \Pi_h^{0}  \div (\nabla z) - \div (\nabla z) ) ( \psi_h - \Pi_h^{0} \psi_h) dx - \int_{\Omega} ( \div \mathcal{I}_{h}^{RT}  (\nabla z) ) \psi_h dx\\
&\quad =  - \int_{\Omega} \div (\nabla z) \left( \psi_h - \Pi_h^{0} \psi_h \right) dx \\
&\quad \quad  - \int_{\Omega} ( \nabla z -  \mathcal{I}_{h}^{RT}  (\nabla z) ) \cdot \nabla_h \psi_h dx  + \int_{\Omega}  \nabla z \cdot \nabla_h \psi_h dx \\
&\quad \quad - \sum_{F \in \mathcal{F}_h^i} \int_{F}  \{ \! \{  \nabla z \} \!\}_{\omega,F} \cdot n_F \Pi_F^0 [\![ \psi_{h} ]\!]_F ds - \sum_{F \in \mathcal{F}_h^{\partial}} \int_{F} (  \nabla z \cdot n_F) \Pi_F^0 \psi_{h} ds.
\end{align*}
Using the H\"older's inequality, the error estimates \eqref{L2ortho}, \eqref{RT5}, \eqref{RT6}, {Lemma \ref{lem4} and Lemma \ref{lem5}}, we have
\begin{align*}
\displaystyle
\| \psi_h \|^2_{L^2(\Omega)}
&\leq \| \varDelta z \|_{L^2(\Omega)} \| \psi_h - \Pi_h^{0} \psi_h \|_{L^2(\Omega)} + \|  \nabla z -  \mathcal{I}_{h}^{RT}  (\nabla z) \|_{L^2(\Omega)} |\psi_h|_{H^1(\mathbb{T}_h)} \\
&\quad + |z|_{H^1(\Omega)} |\psi_h|_{H^1(\mathbb{T}_h)} + c \| \nabla z \|_{H^1(\Omega)^d} |\psi_h|_{J*} \\
&\leq c( h + 1 ) \| \psi_h \|_{L^2(\Omega)} |\psi_h|_{h*},
\end{align*}
which leads to the target inequality if $h \leq 1$.
\qed
\end{pf*}

\section{WOPSIP method for the Stokes equation} \label{sec=WOPSIP}
This section provides an analysis of the WOPSIP method for the Stokes equations on anisotropic meshes.

\subsection{WOPSIP method} \label{sec3=1}
We consider the WOPSIP method for the Stokes equation \eqref{stokes1} as follows. We aim to find $(u_h,p_h) \in V_{dc,h}^{CR} \times Q_h^{0}$ such that
\begin{subequations} \label{wop=1}
\begin{align}
\displaystyle
\nu a_h^{wop}(u_h,v_h) + b_h(v_h , p_h) &= \int_{\Omega} f \cdot v_h dx \quad \forall v_h \in V_{dc,h}^{CR}, \label{wop=1a} \\
b_h(u_h , q_h) &= 0 \quad \forall q_h \in Q_h^{0}, \label{wop=1b}
\end{align}
\end{subequations}
where $a_h^{wop}: (V + V_{dc,h}^{CR}) \times (V + V_{dc,h}^{CR}) \to \mathbb{R}$ and $b_h: (V+V_{dc,h}^{CR}) \times Q_h^{0} \to \mathbb{R}$ respectively denote bilinear forms defined by
\begin{align*}
\displaystyle
a_h^{wop}(u_h,v_h) &:= \sum_{i=1}^d a_h^{i}(u_{h,i},v_{h,i}), \\
a_h^{i}(u_{h,i},v_{h,i}) &:= \int_{\Omega} \nabla_h u_{h,i} \cdot \nabla_h v_{h,i} dx + \sum_{F \in \mathcal{F}_h} \kappa_F  \int_F \Pi_F^{0} [\![ u_{h,i}]\!] \Pi_F^{0} [\![ v_{h,i}]\!] ds, \\
b_h(v_h , q_h) &:= - \int_{\Omega} \divh v_h q_h dx,
\end{align*}
where $\{ u_{h,i}\}_{i=1}^d$ and $\{ v_{h,i}\}_{i=1}^d$ respectively denote the Cartesian components of $u_h$ and $v_h$. Recall that the parameter $\kappa_F$ is defined in \eqref{penalty0}. Using the H\"older's inequality, we obtain
\begin{align}
\displaystyle
|a_h^{wop}(u(h),v_h) | &\leq c |u(h)|_{V_{dc,h}^{CR}}  |v_h|_{V_{dc,h}^{CR}} \quad \forall u(h) \in V + V_{dc,h}^{CR}, \  \forall v_h \in V_{dc,h}^{CR}, \label{ah=sta} \\
|b_h(u(h),q_h) | &\leq c |u(h)|_{V_{dc,h}^{CR}}  \| q_h \|_{Q_{h}^0} \quad \forall u(h) \in V + V_{dc,h}^{CR}, \  \forall q_h \in Q_{h}^0. \label{bh=sta}
\end{align}
We set 
\begin{align*}
\displaystyle
V_{dc,h,\div}^{CR} := \{ v_h \in V_{dc,h}^{CR}  ; \ b_h(v_h,q_h) = 0 \ \forall q_h \in Q_h^0 \}.
\end{align*}
We then observe that
\begin{align*}
\displaystyle
\nu a_h^{wop}(v_h,v_h) = \nu {| v_h |^2_{V_{dc,h}^{CR}}} \quad \forall v_h \in V_{dc,h,\div}^{CR}.
\end{align*}

\begin{rem} \label{rem3}
Let $x_F$ be the barycentre point of $F \in \mathcal{F}_h$. Let $F \in \mathcal{F}_h^i$ with $F = T_{1} \cap T_{2}$, $T_{1},T_{2} \in \mathbb{T}_h$. We  define a finite element space as
\begin{align*}
\displaystyle
V_{h0}^{CR} &:= \biggl \{ v_h \in V_{dc,h}^{CR}; \  v_h|_{T_1 \cap F}(x_F) = v_h|_{T_2 \cap F}(x_F) \quad \forall F \in \mathcal{F}_h^i, \\
&\hspace{1cm} v_h|_F (x_F) = 0 \quad \forall F \in \mathcal{F}_h^{\partial} \biggr \} \subset V_{dc,h}^{CR},
\end{align*}
with a norm
\begin{align*}
\displaystyle
|v_h|_{V_{h0}^{CR}} := |v_h|_{H^1(\mathbb{T}_h)^d}.
\end{align*}
We note that the space $V_{h0}^{CR}$ is the classical CR finite element space on a conforming mesh.

For $v_h \in V_{h0}^{CR}$, by the midpoint rule, for $F \in \mathcal{F}_h^i$,
\begin{align}
\displaystyle
\Pi_F^{0} [\![ v_h ]\!] = \frac{1}{|F|_{d-1}} \int_F [\![ v_h ]\!] ds = [\![ v_h ]\!] (x_F) = v_h|_{T_1 \cap F}(x_F) - v_h|_{T_2 \cap F}(x_F) = 0. \label{wop=2}
\end{align} 
By {an} analogous argument, 
\begin{align}
\displaystyle
\Pi_F^{0} v_h|_F  = 0 \quad \forall F \in \mathcal{F}_h^{\partial}.  \label{wop=2b}
\end{align} 
\end{rem}

\subsection{Stability of the WOPSIP method}
We show the discrete inf-sup condition using Fortin's criterion \cite{For77}.

\begin{lem}[Inf-sup stability] \label{lem7}
The nonconforming Stokes element of type $V_{dc,h}^{CR} \times Q_h^{0}$ satisfies the uniform inf-sup stability condition
\begin{align}
\displaystyle
\inf_{q_h \in Q_h^{0}} \sup_{v_h \in V_{dc,h}^{CR}} \frac{b_h(v_h,q_h)}{|v_h|_{V_{dc,h}^{CR}} \| q_h \|_{Q_h^{0}}} \geq \beta_* := \beta. \label{wop=10}
\end{align}
\end{lem}

\begin{pf*}
Let $v \in V$ and $q_h \in  Q_h^{0}$. Using $\nabla q_h \equiv 0$ on $T$ and the definition of $\mathcal{I}_{h}^{CR}$, we obtain
\begin{align}
\displaystyle
 \textcolor{blue}{b}(v , q_h)
 &=  \textcolor{blue}{-} \sum_{T \in \mathbb{T}_h} \int_T \div v q_h dx 
 =  \textcolor{blue}{-} \sum_{T \in \mathbb{T}_h} \int_{\partial T} (v \cdot n_T) q_h ds \nonumber\\
  &=   \textcolor{blue}{-}  \sum_{T \in \mathbb{T}_h} \int_{\partial T} (\mathcal{I}_{h}^{CR} v \cdot n_T) q_h ds =  \textcolor{blue}{-} \sum_{T \in \mathbb{T}_h} \int_T \div ( \mathcal{I}_{h}^{CR} v ) q_h dx \nonumber \\ 
&=  b_h(\mathcal{I}_{h}^{CR} v , q_h). \label{wop=11}
\end{align}
Given that $\varDelta ({I}_{T}^{CR} v_i) \equiv 0$ on $T \in \mathbb{T}_h$ and $n_T \cdot \nabla ({I}_{T}^{CR} v_i) \in \mathcal{P}^{0}$ on faces of $T$, and using the definition of ${I}_{T}^{CR}$ and the H\"older's inequality, we have
\begin{align*}
\displaystyle
| \mathcal{I}_{h}^{CR} v |_{H^1(\mathbb{T}_h)^d}^2
&= \sum_{T \in \mathbb{T}_h} \sum_{i=1}^d | {I}_{T}^{CR}v_i |_{H^1(T)}^2   = \sum_{T \in \mathbb{T}_h} \sum_{i=1}^d \int_{\partial T} n_T \cdot \nabla ({I}_{T}^{CR} v_i) {I}_{T}^{CR}v_i ds  \\
&=  \sum_{T \in \mathbb{T}_h} \sum_{i=1}^d \int_T \nabla ({I}_{T}^{CR} v_i) \cdot \nabla v_i dx \leq | \mathcal{I}_{h}^{CR} v |_{H^1(\mathbb{T}_h)^d} |v|_{H^1(\Omega)^d},
\end{align*}
which leads to
\begin{align}
\displaystyle
| \mathcal{I}_{h}^{CR} v |_{H^1(\mathbb{T}_h)^d} \leq  |v|_{H^1(\Omega)^d}. \label{wop=12}
\end{align}
Using the definition of the $L^2$-projection $\Pi_F^{0}$ and the definition of $\mathcal{I}_{h}^{CR} v$, we have for $i=1,\ldots,d$,
\begin{align}
\displaystyle
\| \Pi_F^0 [\![  (\mathcal{I}_{h}^{CR} v)_i ] \!] \|^2_{L^2(F)}
&= \int_{F} \Pi_F^0 [\![  (\mathcal{I}_{h}^{CR} v)_i ] \!] \Pi_F^0 [\![  (\mathcal{I}_{h}^{CR} v)_i ] \!] ds \notag \\
&= \int_{F}  [\![ v_i ] \!] \Pi_F^0 [\![  (\mathcal{I}_{h}^{CR} v)_i ] \!] ds = 0, \label{wop=13}
\end{align}
because $[\![ v_i ] \!]_F = 0$ for $v_i \in H_0^1(\Omega)$ and $i \in \{ 1 , \ldots , d \}$.

Using the continuous inf-sup condition \eqref{stokes3}, \eqref{wop=11}, \eqref{wop=12} and \eqref{wop=13}, we have
\begin{align*}
\displaystyle
\beta   \| q_h \|_{Q_h^{0}}
&\leq \sup_{v \in V} \frac{b(v , q_h)}{ | v |_{H^1(\Omega)^d} }
= \sup_{v \in V} \frac{{b_h}( \mathcal{I}_{h}^{CR} v , q_h)}{ | \mathcal{I}_{h}^{CR} v |_{V_{dc,h}^{CR}}} \frac{| \mathcal{I}_{h}^{CR} v |_{V_{dc,h}^{CR}}}{ | v |_{H^1(\Omega)^d}} \\
&\leq  \sup_{v \in  V} \frac{{b_h}(\mathcal{I}_{h}^{CR} v , q_h)}{ | \mathcal{I}_{h}^{CR} v |_{V_{dc,h}^{CR}}} 
\leq  \sup_{v_h \in  V_{dc,h}^{CR}} \frac{{b_h}( v_h , q_h)}{ | v_h |_{{V_{dc,h}^{CR}}}},
\end{align*}
which is the target inequality.
\qed
\end{pf*}

\begin{thr} [Stability of the WOPSIP method] \label{thr7}
Assume that $\Omega$ is convex. {Let $h \leq 1$.} For any $f \in L^2(\Omega)^d$, let $(u_h,p_h) \in V_{dc,h}^{CR} \times Q_h^0$ be the solution of \eqref{wop=1}. Then,
\begin{align}
\displaystyle
|u_h|_{V_{dc,h}^{CR}} \leq \frac{{{C_{dc}^P}}}{\nu} \| f \|_{L^2(\Omega)^d}, \quad \| p_h \|_{Q_h^{0}} \leq  \frac{{\textcolor{blue}{c}}}{\beta_*} \| f \|_{L^2(\Omega)^d}, \label{wop=14}
\end{align}
{where $C_{dc}^P$ is the constant from the discrete Poincar\'e inequality \eqref{wop=9}.}
\end{thr}

\begin{pf*}
{Let $(u_h,p_h) \in V_{dc,h}^{CR} \times Q_h^0$ be the solution of \eqref{wop=1}. Recall that $ | u_h |_{V_{dc,h*}^{CR}} \leq  | u_h |_{V_{dc,h}^{CR}}$ for $h \leq 1$.} Setting $v_h := u_h$ in \eqref{wop=1a} and $q_h := p_h$ in \eqref{wop=1b} and the discrete Poincar\'e inequality \eqref{wop=9} yields
\begin{align*}
\displaystyle
\nu | u_h |^2_{V_{dc,h}^{CR}}
&\leq \| f  \|_{L^2(\Omega)^d} \|  u_h \|_{L^2(\Omega)^d} \\
&\leq {C_{dc}^P} \|  f \|_{L^2(\Omega)^d}  | u_h |_{V_{dc,h*}^{CR}} \leq {C_{dc}^P} \|  f \|_{L^2(\Omega)^d}  | u_h |_{V_{dc,h}^{CR}},
\end{align*}
which leads to the first inequality of \eqref{wop=14}.

For the estimate of the pressure, using the inf-sup condition \eqref{wop=10}, the equation \eqref{wop=1a}, the H\"older's inequality, the discrete Poincar\'e inequality \eqref{wop=9}, and the first inequality of \eqref{wop=14} yields
\begin{align*}
\displaystyle
\beta_* \| p_h \|_{Q_h^0}
&\leq \sup_{v_h \in V_{dc,h}^{CR}} \frac{|b_h(v_h,p_h)|}{|v_h|_{{V_{dc,h}^{CR}}}} \\
&= \sup_{v_h \in V_{dc,h}^{CR}} \frac{|\int_{\Omega} f  \cdot  v_h dx - \nu a_h^{wop} (u_h,v_h)|}{|v_h|_{V_{dc,h}^{CR}}} \\
&\leq  \sup_{v_h \in V_{dc,h}^{CR}} \frac{ {C_{dc}^P} \| f \|_{L^2(\Omega)^d} |v_h|_{V_{dc,h}^{CR}} + {c} \nu |u_h|_{V_{dc,h}^{CR}} |v_h|_{V_{dc,h}^{CR}} }{|v_h|_{V_{dc,h}^{CR}}} \\
&=  {C_{dc}^P} \| f \|_{L^2(\Omega)^d} + {c} \nu |u_h|_{V_{dc,h}^{CR}}
\leq c \| f \|_{L^2(\Omega)^d},
\end{align*}
which leads to the second inequality of \eqref{wop=14}.
\qed
\end{pf*}

\subsection{Error estimates of the WOPSIP method}

\begin{lem} \label{lem8}
Assume that $\Omega$ is convex. Let $(u,p) \in V \times Q$ be the solution of \eqref{stokes1} and $(u_h,p_h) \in V_{dc,h}^{CR} \times Q_h^{0}$ be the solution of \eqref{wop=1}. It then holds that
\begin{align}
\displaystyle
&|u-u_h|_{V_{dc,h}^{CR}} + \| p - p_h \|_{Q_h^{0}} \notag \\
&\quad \leq C(\beta_*) \left\{ \inf_{v_h \in V_{dc,h}^{CR}}  |u-v_h|_{V_{dc,h}^{CR}} + \frac{1}{\nu} \inf_{q_h \in Q_h^{0}}  \| p - q_h \|_{Q_h^{0}} + \frac{1}{\nu} E_h(u,p) \right\}, \label{wop=15}
\end{align}
where
\begin{align}
\displaystyle
E_h(u,p) := \sup_{w_h \in V_{dc,h}^{CR}} \frac{ |\nu a_h^{wop}(u,w_h) - \int f \cdot w_h dx + b_h(w_h,p) |}{|w_h|_{V_{dc,h}^{CR}}}. \label{wop=16}
\end{align}
\end{lem}

\begin{pf*}
A proof can be similarly completed for the CR finite element approximation and standard; e.g., see \cite[Lemma 14.3.1]{Ish22a}.
\qed
\end{pf*}

The essential part for error estimates is the consistency error term \eqref{wop=16}.

\begin{lem}[Asymptotic Consistency]  \label{lem9}
Assume that $\Omega$ is convex. Let $(u,p) \in (V \cap H^{2}(\Omega)^d) \times (Q \cap H^{1}(\Omega))$ be the solution of the homogeneous Dirichlet Stokes problem \eqref{stokes1} with data $f \in L^2(\Omega)^d$. Let $\{ \mathbb{T}_h\}$ be a family of conformal meshes with the semi-regular property (Assumption \ref{ass1}). Let $T \in \mathbb{T}_h$ be the element with Conditions \ref{cond1} or \ref{cond2} and satisfy (Type \roman{sone}) in Section \ref{element=cond} when $d=3$. Then,
\begin{align}
\displaystyle
&E_h(u,p) \notag \\
&\quad \leq  c \nu \left( \sum_{i,j=1}^d \sum_{T \in \mathbb{T}_h} h_j^2 \left \| \frac{\partial}{\partial r_j} \nabla u_i \right \|_{L^2(T)^d}^2 \right)^{\frac{1}{2}} + c \nu h \| \varDelta u \|_{L^2(\Omega)^d} + c h \| f \|_{L^2(\Omega)^d} \notag \\
&\quad \quad +  c h  |p|_{H^1(\Omega)}  + c h^{\frac{3}{2}}  \| f \|_{L^2(\Omega)^d}^{\frac{1}{2}} \left( \nu \| \varDelta u \|_{L^2(\Omega)^d} + | p |_{H^1(\Omega)} \right)^{\frac{1}{2}}. \label{stokes=17}
\end{align}
Furthermore, let $d=3$ and let \textcolor{blue}{$T \in \mathbb{T}_h$} be the element with Condition \ref{cond2} and satisfy (Type \roman{stwo}) in Section \ref{element=cond}. It then holds that
\begin{align}
\displaystyle
E_h(u,p)
&\leq c \nu h \| \varDelta u \|_{L^2(\Omega)^3} + c h | p |_{H^1(\Omega)} + c h \| f \|_{L^2(\Omega)^3} \label{stokes=18}
\end{align}
if $|u|_{H^2(\Omega)^3} \leq c \| \varDelta u\|_{L^2(\Omega)^3}$ holds.
\end{lem}

\begin{pf*}
For  $i=1,\ldots,d$, we first have
\begin{align*}
\displaystyle
 \div ( \mathcal{I}_{h}^{RT} (\nu  \nabla u_i   - p e_i ) ) &= \Pi_h^{0} \div (\nu \nabla u_i - p e_i) =  \Pi_h^{0} \left( \nu \varDelta u_i - \frac{\partial p}{\partial x_i} \right)
 = - \Pi_h^{0} f_i,
\end{align*}
where $(e_1,\ldots,e_d)$ denotes the Cartesian basis of $\mathbb{R}^d$.

{Let $w_h \in V_{dc,h}^{CR}$.} Setting $w := \nu  \nabla u_i   - p e_i $ in \eqref{wop=3} yields
\begin{align*}
\displaystyle
&\nu a_h^{wop}(u,w_h) - \int_{\Omega} f \cdot  w_h dx + b_h(w_h,p) \\
&\quad =  \sum_{i=1}^d \Biggl \{ \nu \int_{\Omega} ( \nabla u_i - \mathcal{I}_{h}^{RT} \nabla u_i) \cdot \nabla_h w_{h,i} dx   \\
&\quad \quad + \sum_{F \in \mathcal{F}_h^i} \int_{F}  \{ \! \{ w \} \!\}_{\omega,F} \cdot n_F \Pi_F^0 [\![ w_{h,i} ]\!]_F ds + \sum_{F \in \mathcal{F}_h^{\partial}} \int_{F} ( w \cdot n_F) \Pi_F^0 w_{h,i} ds  \\
&\quad \quad - \int_{\Omega} \left( f_i -  \Pi_h^0 f_i \right) w_{h,i}  dx  + \int_{\Omega} \left( \mathcal{I}_{h}^{RT} (p e_i) - (p e_i) \right) \cdot \nabla_h w_{h,i} dx { \Biggr \}} \\
&\quad =: I_1 + I_2 + I_3 + I_4 + I_5. 
\end{align*}

Let $T \in \mathbb{T}_h$ be the element with Conditions \ref{cond1} or \ref{cond2} and satisfy (Type \roman{sone}) in Section \ref{element=cond} when $d=3$. Using the H\"older's inequality, the Cauchy--Schwarz inequality and the RT interpolation error \eqref{RT5}, the term $I_1$ is estimated as
\begin{align*}
\displaystyle
|I_1|
&\leq c \nu  \sum_{i=1}^d \sum_{T \in \mathbb{T}_h} \| \nabla u_i - \mathcal{I}_{h}^{RT} \nabla u_i \|_{L^2(T)} | w_{h,i} |_{H^1(T)} \\
&\leq c \nu  \sum_{i=1}^d \sum_{T \in \mathbb{T}_h}\left(  \sum_{j=1}^d h_j \left \| \frac{\partial}{\partial r_j} \nabla u_i \right \|_{L^2(T)^d} + h_T \| \varDelta u_i \|_{L^2(T)} \right) | w_{h,i} |_{H^1(T)} \\
&\leq c \nu \left\{ \left( \sum_{i,j=1}^d \sum_{T \in \mathbb{T}_h} h_j^2 \left \| \frac{\partial}{\partial r_j} \nabla u_i \right \|_{L^2(T)^d}^2 \right)^{\frac{1}{2}} + h \| \varDelta u \|_{{L^2(\Omega)^d}} \right\}  | w_{h} |_{V_{dc,h}^{CR}}.
\end{align*}
Using the Cauchy--Schwarz and Jensen inequalities for $a_1,a_2,a_3 \in \mathbb{R}_+$ and $b_1 , b_2 \in \mathbb{R}_+ \cup \{ 0\}$,
\begin{align*}
\displaystyle
&\sum_{i=1}^d a_1 (a_2+b_1)^{\frac{1}{2}}(a_3+b_2)^{\frac{1}{2}} \\
&\leq c \left( \sum_{i=1}^d a_1^2 \right)^{\frac{1}{2}} \left\{ \left( \sum_{i=1}^d a_2^2 \right)^{\frac{1}{2}} + \left( \sum_{i=1}^d b_1^2 \right)^{\frac{1}{2}} \right\}^{\frac{1}{2}} \left\{ \left( \sum_{i=1}^d a_3^2 \right)^{\frac{1}{2}} + \left( \sum_{i=1}^d b_2^2 \right)^{\frac{1}{2}} \right\}^{\frac{1}{2}}.
\end{align*}
According to the previous inequality, the triangle inequality, \eqref{stokes6}, \eqref{wop=5} and \eqref{wop=6}, the terms $I_2$ and $I_3$ are estimated as
\begin{align*}
\displaystyle
|I_2| 
&\leq c \sum_{i=1}^d  |w_{h,i}|_J \Biggl(  h \| \nu  \nabla u_i   - p e_i \|_{L^2(\Omega)^d}\\
&\quad + h^{\frac{3}{2}} \| \nu  \nabla u_i   - p e_i \|_{L^2(\Omega)^d}^{\frac{1}{2}} | \nu  \nabla u_i   - p e_i |_{H^1(\Omega)^d}^{\frac{1}{2}} \Biggr) \\
&\leq c \sum_{i=1}^d  |w_{h,i}|_J \Biggl\{ h (\nu |u_i|_{H^1(\Omega)} + \| p e_i \|_{L^2(\Omega)^d} )  \\
&\quad + h^{\frac{3}{2}}  (\nu |u_i|_{H^1(\Omega)} + \| p e_i \|_{L^2(\Omega)^d} )^{\frac{1}{2}} (\nu \| \varDelta u_i \|_{L^2(\Omega)} + | p |_{H^1(\Omega)} )^{\frac{1}{2}}  \Biggr\} \\
&\leq c | w_{h} |_{V_{dc,h}^{CR}} \Biggl\{ h (\nu |u|_{H^1(\Omega)^d} + \| p \|_{L^2(\Omega)}) \\
&\quad + h^{\frac{3}{2}} \left( \nu |u|_{H^1(\Omega)^d} + \| p \|_{L^2(\Omega)} \right)^{\frac{1}{2}} \left( \nu \| \varDelta u \|_{L^2(\Omega)^d} + | p |_{H^1(\Omega)} \right)^{\frac{1}{2}} \Biggr \} \\
&\leq c | w_{h} |_{V_{dc,h}^{CR}}  \Biggl\{ h \| f \|_{L^2(\Omega)^d} + h^{\frac{3}{2}}  \| f \|_{L^2(\Omega)^d}^{\frac{1}{2}} \left( \nu \| \varDelta u \|_{L^2(\Omega)^d} + | p |_{H^1(\Omega)} \right)^{\frac{1}{2}} \Biggr \}, \\
|I_3|
&\leq c | w_{h} |_{V_{dc,h}^{CR}}  \Biggl\{ h \| f \|_{L^2(\Omega)^d} + h^{\frac{3}{2}}  \| f \|_{L^2(\Omega)^d}^{\frac{1}{2}} \left( \nu \| \varDelta u \|_{L^2(\Omega)^d} + | p |_{H^1(\Omega)} \right)^{\frac{1}{2}} \Biggr \}.
\end{align*}
Using the H\"older's inequality, the Cauchy--Schwarz inequality, the stability of $\Pi_h^0$ and the estimate \eqref{L2ortho}, the term $I_4$ is estimated as
\begin{align*}
\displaystyle
|I_4|
&= \left| \sum_{i=1}^d \int_{\Omega} \left( f_i -  \Pi_h^0 f_i \right) \left(  w_{h,i} - \Pi_h^0 w_{h,i}  \right) dx \right| \\
&\leq \sum_{i=1}^d \sum_{T \in \mathbb{T}_h} \| f_i -  \Pi_h^0 f_i \|_{L^2(T)}  \| w_{h,i} - \Pi_h^0 w_{h,i} \|_{L^2(T)} \\
&\leq c h \| f \|_{L^2(\Omega)^d}  | w_{h} |_{V_{dc,h}^{CR}}.
\end{align*}
Using the H\"older's inequality, the Cauchy--Schwarz inequality and the RT interpolation error \eqref{RT5}, the term $I_5$ is estimated as
\begin{align*}
\displaystyle
|I_5|
&\leq \sum_{i=1}^d \sum_{T \in \mathbb{T}_h} \| \mathcal{I}_{h}^{RT} (p e_i) - (p e_i)  \|_{L^2(\Omega)}  | w_{h,i} |_{H^1(T)} \\
&\leq c  \sum_{i=1}^d \sum_{T \in \mathbb{T}_h} \left(  \sum_{j=1}^d h_j \left \| \frac{\partial (p e_i)}{\partial r_j} \right \|_{L^2(T)^d} + h_T \| \div (p e_i) \|_{L^2(T)} \right) | w_{h,i} |_{H^1(T)} \\
&\leq c \left\{ \left( \sum_{j=1}^d \sum_{T \in \mathbb{T}_h} h_j^2 \left \| \frac{\partial p}{\partial r_j} \right \|_{L^2(T)}^2  \right)^{\frac{1}{2}} + h |p|_{H^1(\Omega)} \right \} | w_{h} |_{V_{dc,h}^{CR}} \\
&\leq c h  |p|_{H^1(\Omega)}  | w_{h} |_{V_{dc,h}^{CR}}.
\end{align*}
Gathering the above inequalities yields the target inequality \eqref{stokes=17}.

Let $d=3$ and let {$T \in \mathbb{T}_h$} be the element with Condition \ref{cond2} and satisfy (Type \roman{stwo}) in Section \ref{element=cond}. Using the H\"older's inequality, the Cauchy--Schwarz inequality and the RT interpolation error \eqref{RT6}, the terms $I_1$ and $I_5$ are estimated as
\begin{align*}
\displaystyle
|I_1|
&\leq c \nu h  |  u |_{H^2(\Omega)^3}  | w_{h} |_{V_{dc,h}^{CR}}, \\
|I_5|
&\leq c h |p|_{H^1(\Omega )}  | w_{h} |_{V_{dc,h}^{CR}}, 
\end{align*}
which implies the target inequality \eqref{stokes=18}.
\qed
\end{pf*}

\begin{thr}[Error Estimate] \label{thr8}
Assume that $\Omega$ is convex. Let $(u,p) \in (V \cap H^{2}(\Omega)^d) \times (Q \cap H^{1}(\Omega))$ be the solution to the homogeneous Dirichlet problem for the Stokes equation \eqref{stokes1} with data $f \in L^2(\Omega)^d$. Let $\{ \mathbb{T}_h\}$ be a family of conformal meshes with the semi-regular property (Assumption \ref{ass1}). Let $T \in \mathbb{T}_h$ be the element with Conditions \ref{cond1} or \ref{cond2} and satisfy (Type \roman{sone}) in Section \ref{element=cond} when $d=3$. {Let $(u_h,p_h) \in V_{dc,h}^{CR} \times Q_h^{0}$ be the solution of \eqref{wop=1}.} Then,
\begin{align}
\displaystyle
&|u-u_h|_{V_{dc,h}^{CR}} + \| p - p_h \|_{Q_h^{0}} \notag \\
&\quad \leq  c  \left( \sum_{i,j=1}^d \sum_{T \in \mathbb{T}_h} h_j^2 \left \| \frac{\partial}{\partial r_j} \nabla u_i \right \|_{L^2(T)^d}^2 \right)^{\frac{1}{2}} + c h \| \varDelta u \|_{L^2(\Omega)^d} + c \frac{h}{\nu} \| f \|_{L^2(\Omega)^d}  \notag \\
&\quad \quad  +  c \frac{h}{\nu}  |p|_{H^1(\Omega)} + c \frac{h^{\frac{3}{2}}}{\nu}  \| f \|_{L^2(\Omega)^d}^{\frac{1}{2}} \left( \nu \| \varDelta u \|_{L^2(\Omega)^d} + | p |_{H^1(\Omega)} \right)^{\frac{1}{2}}. \label{stokes=19}
\end{align}
Furthermore, let $d=3$ and let {$T \in \mathbb{T}_h$} be the element with Condition \ref{cond2} and satisfy (Type \roman{stwo}) in Section \ref{element=cond}. Then
\begin{align}
\displaystyle
&|u-u_h|_{V_{dc,h}^{CR}} + \| p - p_h \|_{Q_h^{0}} \notag \\
&\quad \leq  c h \| \varDelta u \|_{L^2(\Omega)^3} + c \frac{h}{\nu} | p |_{H^1(\Omega)} + c \frac{h}{\nu} \| f \|_{L^2(\Omega)^3} \label{stokes=20}
\end{align}
if $|u|_{H^2(\Omega)^3} \leq c \| \varDelta u\|_{L^2(\Omega)^3}$ holds.
\end{thr}

\begin{pf*}
{Let $v \in V$.} Using the definition of the $L^2$-projection $\Pi_F^{0}$ and the definition of $\mathcal{I}_{h}^{CR} v$, for $i=1,\ldots,d$, we have
\begin{align}
\displaystyle
\| \Pi_F^0 [\![  (\mathcal{I}_{h}^{CR} v)_i - v_i] \!] \|^2_{L^2(F)}
&= \int_{F} \Pi_F^0 [\![  (\mathcal{I}_{h}^{CR} v)_i  - v_i] \!] \Pi_F^0 [\![  (\mathcal{I}_{h}^{CR} v)_i - v_i] \!] ds \notag \\
&= \int_{F}  [\![  (\mathcal{I}_{h}^{CR} v)_i  - v_i ] \!] \Pi_F^0 [\![  (\mathcal{I}_{h}^{CR} v)_i - v_i] \!] ds = 0. \label{stokes=21}
\end{align}
Therefore, using \eqref{CR4} and \eqref{stokes=21}, we obtain
\begin{align}
\displaystyle
 \inf_{v_h \in {V_{dc,h}^{CR}}}  |u-v_h|_{V_{dc,h}^{CR}}
 &\leq  |u-\mathcal{I}_{h}^{CR} u |_{V_{dc,h}^{CR}}  = |u-\mathcal{I}_{h}^{CR} u|_{H^1(\mathbb{T}_h)^d} \notag\\
 &\leq c  \left( \sum_{i,j=1}^d \sum_{T \in \mathbb{T}_h} h_j^2 \left \| \frac{\partial}{\partial r_j} \nabla u_i \right \|_{L^2(T)^d}^2 \right)^{\frac{1}{2}} \label{stokes=22} \\
 &\leq c h |u|_{H^2(\Omega)^d}.  \label{stokes=23}
\end{align}
The estimate \eqref{L2ortho} yields
\begin{align}
\displaystyle
  \inf_{q_h \in Q_h^{0}}  \| p - q_h \|_{Q_h^{0}} 
  &\leq  \| p - \Pi_h^0 p \|_{Q_h^{0}} \notag \\
&  \leq c \left( \sum_{j=1}^d \sum_{T \in \mathbb{T}_h} h_j^2 \left \| \frac{\partial p}{\partial r_j} \right \|_{L^2(T)}^2  \right)^{\frac{1}{2}} \leq c h |p|_{H^1(\Omega)}. \label{stokes=24}
\end{align}
From \eqref{wop=15}, \eqref{stokes=17}, \eqref{stokes=22} and \eqref{stokes=24}, we conclude by deriving the desired inequality \eqref{stokes=19}.

The estimate \eqref{stokes=20} is deduced from \eqref{wop=15}, \eqref{stokes=18}, \eqref{stokes=23} and \eqref{stokes=24}.
\qed
\end{pf*}

\section{Numerical experiments} \label{numerical=sec}
{Let $d=2$ and $\Omega := (0,1)^2$. } {We have conducted experiments which confirm that the WOPSIP method does NOT converge when we use the parameter $\kappa_{F*}$. Therefore, in all experiments, we used $\kappa_F$.}

\subsection{Comparison of calculations for two schemes on anisotropic meshes} \label{numerial=comp}
In this section, we present numerical tests for some schemes on anisotropic meshes. We set $\varphi (x_1,x_2) := x_1^2(x_1-1)^2x_2^2(x_2-1)^2$. The function $f$ of the Stokes equation \eqref{intro1} with $\nu = 1$,
\begin{align*}
\displaystyle
-  \varDelta u + \nabla p = f \quad \text{in $\Omega$}, \quad \div u  = 0 \quad \text{in $\Omega$}, \quad u = 0 \quad \text{on $\partial \Omega$}, 
\end{align*}
is given so that the exact solution is
{
\begin{align*}
\displaystyle
\begin{pmatrix}
 u_1  \\
  u_2
\end{pmatrix}
&=
\begin{pmatrix}
 \frac{\partial \varphi}{\partial x_2}  \\
  - \frac{\partial \varphi}{\partial x_1} 
\end{pmatrix}
=
\begin{pmatrix}
2 x_1^2 (x_1-1)^2 x_2 (x_2-1)^2 + 2 x_1^2 (x_1-1)^2 x_2^2 (x_2-1)\\
  - 2  x_1 (x_1-1)^2 x_2^2 (x_2-1)^2 - 2  x_1^2 (x_1-1) x_2^2 (x_2-1)^2 
\end{pmatrix}, \\
p &=  x_1^2 - x_2^2.
\end{align*}
}

\begin{figure}[htbp]
  \begin{minipage}[b]{0.45\linewidth}
    \centering
    \includegraphics[keepaspectratio, scale=0.15]{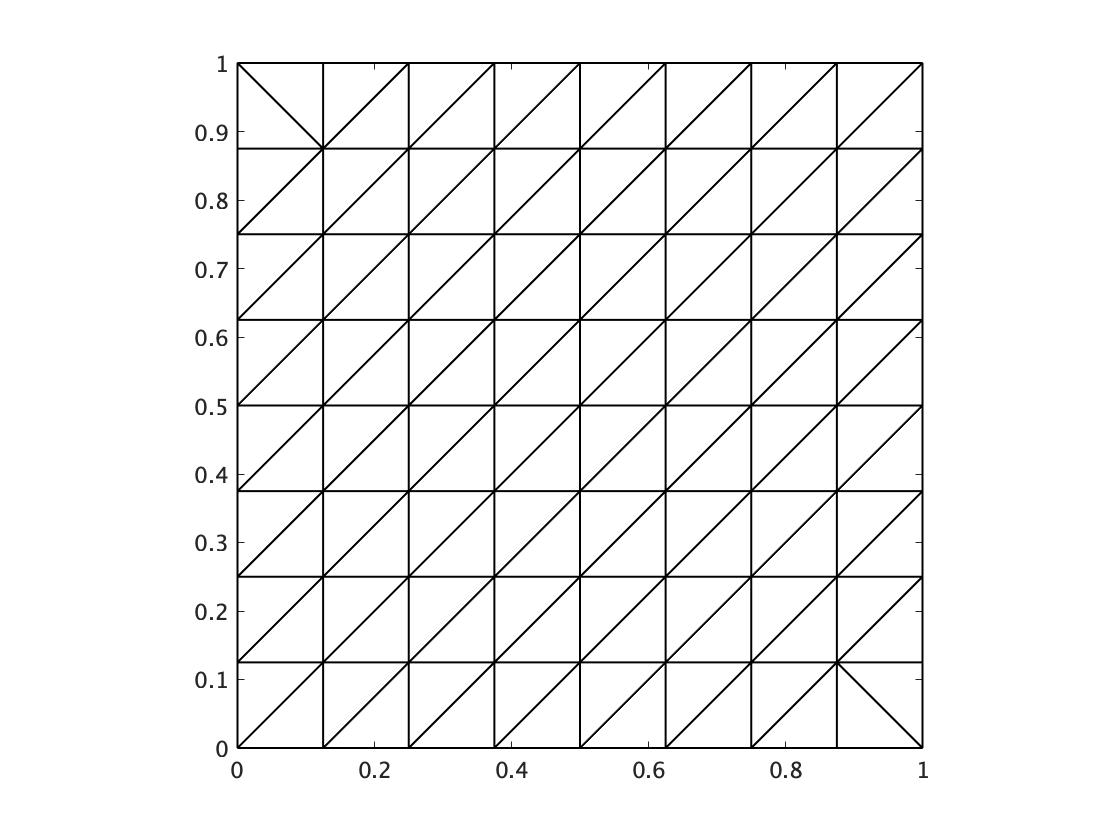}
    \caption{(\Roman{lone}) Standard mesh}
     \label{fig1}
  \end{minipage}
  \begin{minipage}[b]{0.45\linewidth}
    \centering
    \includegraphics[keepaspectratio, scale=0.15]{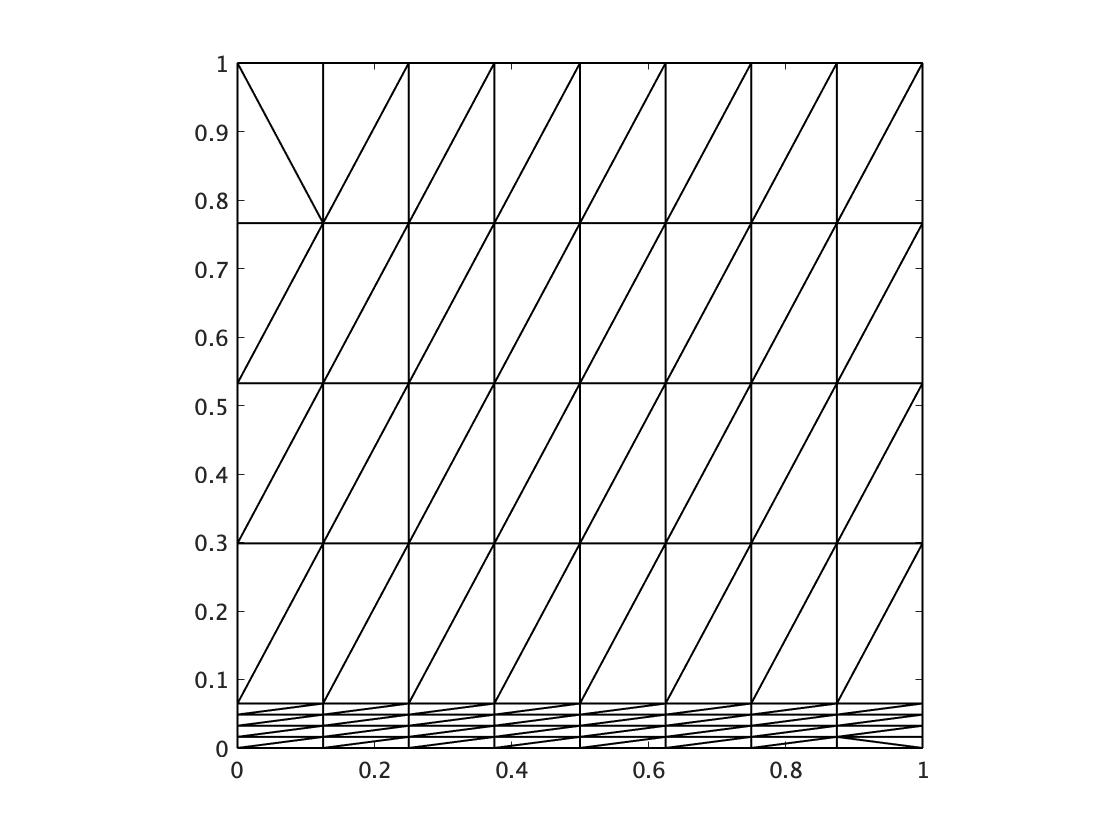}
    \caption{(\Roman{ltwo}) Shishkin mesh, $\delta = \frac{1}{128}$}
     \label{fig2}
  \end{minipage}
\end{figure}

\begin{figure}[htbp]
  \begin{minipage}[b]{0.45\linewidth}
    \centering
    \includegraphics[keepaspectratio, scale=0.15]{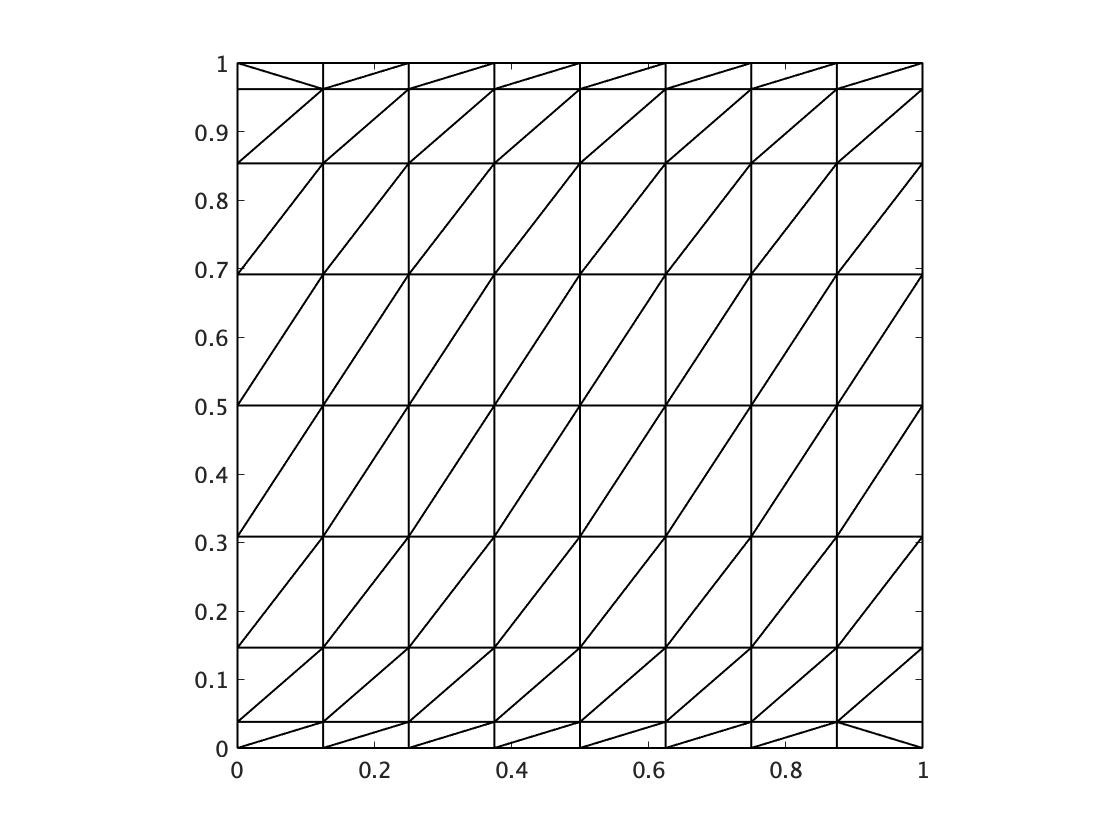}
    \caption{(\Roman{lthree}) Anisotropic mesh}
     \label{fig33}
  \end{minipage}
  \begin{minipage}[b]{0.45\linewidth}
    \centering
    \includegraphics[keepaspectratio, scale=0.15]{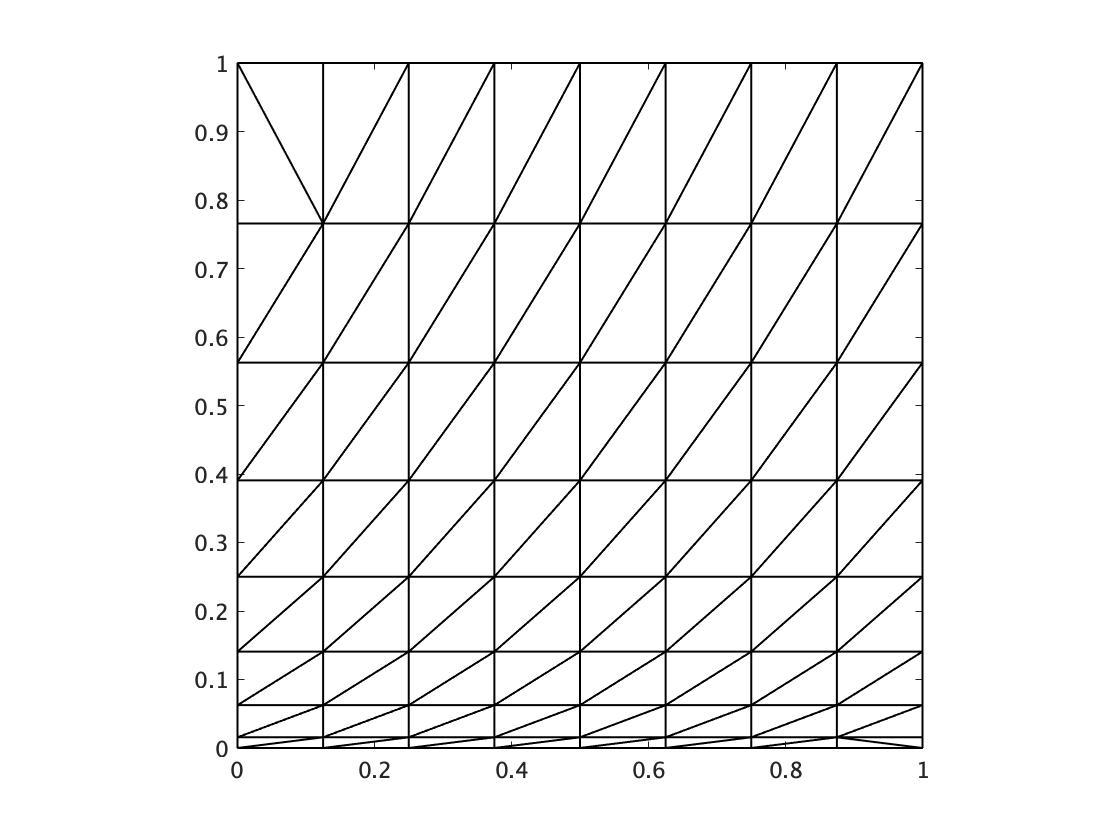}
    \caption{(\Roman{lfour}) Anisotropic mesh}
     \label{fig44}
  \end{minipage}
\end{figure}

Let {$N \in \{ 32,64 \}$} be the division number of each side of the bottom and the height edges of $\Omega$. We consider four types of mesh partitions. Let $(x_1^i, x_2^i)^T$ be grip points of triangulations $\mathbb{T}_h$ defined as follows. Let $i \in \mathbb{N}$.
\begin{description}
  \item[(\Roman{lone}) Standard mesh (Fig. \ref{fig1})] 
\begin{align*}
\displaystyle
x_1^i := \frac{i}{N}, \quad x_2^i := \frac{i}{N}, \quad  i \in \{0, \ldots , N \}.
\end{align*}
  \item[(\Roman{ltwo}) Shishkin mesh (Fig. \ref{fig2})]
\begin{align*}
\displaystyle
x_1^i &:= \frac{i}{N}, \quad  i \in \{0 , \ldots , N \}, \\
 x_2^i &:=
\begin{cases}
\tau \frac{2}{N} i, \quad  i \in \left\{0, \ldots , \frac{N}{2} \right \}, \\
\tau + (1 - \tau) \frac{2}{N} \left( i - \frac{N}{2} \right), \quad i \in \left\{ \frac{N}{2}+1 , \ldots , N \right\},
\end{cases}
\end{align*}
where $\tau := 4 \delta | \ln(N) |$ with $\delta = \frac{1}{128}$, see \cite[Section 2.1.2]{Lin10}.
\item[(\Roman{lthree}) Anisotropic mesh {from \cite{CheLiuQia10}} (Fig. \ref{fig33})]
\begin{align*}
\displaystyle
x_1^i := \frac{i}{N}, \quad x_2^i := \frac{1}{2}\left( 1 - \cos \left( \frac{i \pi}{N} \right) \right), \quad  i \in \{0, \ldots, N \}.
\end{align*}
 \item[(\Roman{lfour}) Anisotropic mesh (Fig. \ref{fig44})]
 \begin{align*}
\displaystyle
x_1^i := \frac{i}{N}, \quad x_2^i := \left ( \frac{i}{N} \right)^{2}, \quad  i \in \{0, \ldots, N \}.
\end{align*}
\end{description}

The shape-regularity condition is known: there exists a constant $\gamma_1 \> 0$ such that
\begin{align*}
\displaystyle
\rho_{T} \geq \gamma_1 h_{T} \quad \forall \mathbb{T}_h \in \{ \mathbb{T}_h \}, \quad \forall T \in \mathbb{T}_h,
\end{align*}
and is equivalent to the following condition. There exists a constant $\gamma_2 \> 0$ such that for any $\mathbb{T}_h \in \{ \mathbb{T}_h \}$ and simplex $T \in \mathbb{T}_h$, we have
\begin{align*}
\displaystyle
|T|_2 \geq \gamma_2 h_{T}^2. 
\end{align*}
A proof is provided in \cite[Theorem 1]{BraKorKri08}. The semiregularity mesh condition defined in \eqref{NewGeo} is equivalent to the maximum angle condition. The following parameters are computed.
\begin{align*}
\displaystyle
\textit{MinAngle} := \max_{T \in \mathbb{T}_h} \frac{|L_3|_1^2}{|T|_2}, \quad \textit{MaxAngle} := \max_{T \in \mathbb{T}_h} \frac{|L_1|_1 |L_2|_1}{|T|_2},
\end{align*}
where $L_i$, $i=1,2,3,$ denote the edges of the simplex $T \in \mathbb{T}_h$ with $|L_1|_1 \leq |L_2|_1 \leq |L_3|_1$. 

\begin{table}[h]
\caption{Mesh conditions}
\centering
\begin{tabular}{l|l|l|l} \hline
 Mesh &$N$ &  $\textit{MinAngle} $ & $\textit{MaxAngle}$   \\ \hline \hline
 \Roman{lone} & 32  & 4.00000 & 2.00000   \\
 & 64  &  4.00000  & 2.00000     \\
\hline
 \Roman{ltwo} & 32  & 9.66647   & 2.00000    \\
 & 64  &  8.21423  & 2.00000     \\
\hline
 \Roman{lthree}  & 32  & 2.61132e+01  & 2.00000   \\
& 64  &  5.19640e+01    &  2.00000    \\
\hline
 \Roman{lfour} & 32  & 6.40625e+01 & 2.00000   \\
 & 64  & 1.28031e+02   & 2.00000    \\
\hline
\end{tabular}
\label{table=mesh}
\end{table}

Notably, a sequence with meshes (\Roman{lone}) or (\Roman{ltwo}) satisfies the shape-regularity condition, but a sequence with meshes (\Roman{lthree}), or (\Roman{lfour}) does not satisfy the shape regularity condition.

We adopt the following schemes.
\begin{description}
  \item[(1)] Scheme \eqref{wop=1}.
 \item[(2)] Well-balanced CR (WBCR) method ({\cite{Lin14}}). We define the classical RT finite element space as
\begin{align*}
\displaystyle
V_{c,h}^{RT} &:= \{ v_h \in V_{dc,h}^{RT} : \  [\![ v_h \cdot n ]\!]_F = 0, \ \forall F \in \mathcal{F}_h \}.
\end{align*}
Let $\mathcal{I}_{h}^{RTc}: V_{h0}^{CR} \to V_{c,h}^{RT}$ be a classical RT interpolation operator.  Find $(u_h,p_h) \in V_{h0}^{CR} \times Q_h^0$ such that
\begin{subequations} \label{stokes51}
\begin{align}
\displaystyle
a_h^{CR}(u_h,v_h) + b_h^{CR}(v_h , p_h) &= \int_{\Omega} f \cdot \mathcal{I}_{h}^{RTc} v_h dx \quad \forall v_h \in V_{h0}^{CR}, \label{stokes51a} \\
b_h^{CR}(u_h , q_h) &= 0 \quad \forall q_h \in Q_h^0, \label{stokes51b}
\end{align}
\end{subequations}
where
\begin{align*}
\displaystyle
a_h^{CR}(u_h,v_h) &:= \sum_{i=1}^d \int_{\Omega} \nabla_h u_{h,i} \cdot \nabla_h v_{h,i} dx, \\
b_h^{CR}(v_h , q_h) &:= - \int_{\Omega} \divh v_h q_h dx,
\end{align*}
for $u_h,v_h \in V_{h0}^{CR}$ and $q_h \in Q_h^0$. Recall that space $V_{h0}^{CR}$ was defined in Remark \ref{rem3}.
\end{description}

If an exact solution $u$ is known, the error $e_h := u - u_h$ and $e_{h/2} := u - u_{h/2}$ are computed numerically for two mesh sizes $h$ and $h/2$. The convergence indicator $r$ is defined by
\begin{align*}
\displaystyle
r = \frac{1}{\log(2)} \log \left( \frac{\| e_h \|_X}{\| e_{h/2} \|_X} \right).
\end{align*}
We compute the convergence order with respect to norms defined by
{
\begin{align*}
\displaystyle
E_{u_h}^{(\text{Mesh No.})} &:= \frac{| u - u_h |_{V_{dc,h}^{CR}} }{|u|_{H^1(\Omega)^2}} \quad \text{for scheme (1)}, \\
E_{u_h}^{(\text{Mesh No.})} &:=  \frac{| u - u_h |_{H^1(\mathbb{T}_h)^2}}{|u|_{H^1(\Omega)^2} } \\
&\quad \quad \text{for schemes (2), (3), and (4)}, \\
 \quad E_{u_h,L^2}^{(\text{Mesh No.})} &:=  \frac{\| u - u_h \|_{L^2(\Omega)^2} }{\|u\|_{L^2(\Omega)^2} }, \quad
 E_{p_h}^{(\text{Mesh No.})} :=  \frac{ \| p - p_h \|_{L^2(\Omega)}}{ \| p \|_{L^2(\Omega)}}.
\end{align*}
}

For the computation of scheme (1), we used the CG method without preconditioners, and the quadrature of the five order for computation of the right-hand side in \eqref{wop=1a} and \eqref{stokes51a}, see \cite[p. 85, Table 30.1]{ErnGue21b}. 

\begin{table}[h]
\caption{Scheme (1)}
\centering
\begin{tabular}{l | l |  l | l | l |  l | l | l} \hline
Mesh & $N$ & $E_{u_h}^{(\text{Mesh No.})}$ & $r$ & $E_{u_h,L^2}^{(\text{Mesh No.})}$  & $r$ & $E_{p_h}^{(\text{Mesh No.})}$ & $r$ \\ \hline \hline
\Roman{lone} & 32 &  8.10569e-01       &    &  2.12630e-01   & &   3.61598e-02  &  \\
& 64 & 4.08981e-01       &  0.99   &  5.42357e-02  & 1.97 &   1.35562e-02  & 1.42  \\
 \hline
   \Roman{ltwo} & 32    &  1.15924   &     &   4.33629e-01    &  & 6.52059e-02      &    \\
  & 64    & 5.79411e-01  &1.00     & 1.08800e-01    & 1.99 & 2.22654e-02     & 1.55    \\
  \hline
   \Roman{lthree}  & 32   & 1.05163   &     &   3.60039e-01     &  & 5.24322e-02     &   \\
  & 64   & 5.34097e-01   &  0.98  &   9.31283e-02     & 1.95  & 1.76734e-02     & 1.57   \\
  \hline
   \Roman{lfour}  &  32  & 1.23942    &     &  4.97459e-01      &  & 7.17788e-02     &  \\
  &  64  & 6.36438e-01   &  0.96   &  1.31655e-01     & 1.92  & 2.44549e-02     &1.55  \\
\hline
\end{tabular}
\label{table01}
\end{table}

\begin{table}[h]
\caption{Scheme (2)}
\centering
\begin{tabular}{l | l |  l | l | l |  l | l | l } \hline
Mesh & $N$ & $E_{u_h}^{(\text{Mesh No.})}$ & $r$ & $E_{u_h,L^2}^{(\text{Mesh No.})}$  & $r$ & $E_{p_h}^{(\text{Mesh No.})}$ & $r$ \\ \hline \hline
 \Roman{lone} & 32 & 1.30431e-01   &    & 1.10175e-02  &  & 2.26926e-02  &   \\
   & 64   & 6.53265e-02 &  1.00  & 2.76911e-03   &1.99  &  1.13420e-02& 1.00  \\ \hline
  \Roman{ltwo} & 32  & 1.77909e-01 &    & 2.07770e-02 &  &3.42518e-02  &  \\
 &  64  & 8.70267e-02 & 1.03   & 5.01619e-03  & 2.05 & 1.67556e-02 & 1.03 \\ \hline
 \Roman{lthree}  & 32   & 1.48023e-01  &    & 1.40474e-02   &  &  2.53116e-02  &  \\
 & 64   & 7.42163e-02  &  1.00  & 3.54266e-03   &  1.99 & 1.26452e-02 & 1.00 \\ \hline
 \Roman{lfour}  &  32  & 1.59293e-01  &    & 1.85984e-02   &  & 3.39050e-02 &  \\
  &  64 & 7.99498e-02   & 0.99   &  4.71503e-03  &1.99 & 1.69444e-02 & 1.00  \\
\hline
\end{tabular}
\label{table=wbcr00}
\end{table}

The table \ref{table01} implies that the inf-sup condition is satisfied on {our} anisotropic meshes. {Furthermore, the table \ref{table=wbcr00} also implies that the inf-sup condition is satisfied on our anisotropic meshes. In the CR and WBCR methods, the inf-sup condition is proved on anisotropic meshes, see \cite{ApeKemLin21}. }


\subsection{Comparison of calculations for a problem with boundary layers}
\begin{figure}[htbp]
  \begin{minipage}[b]{0.3\linewidth}
    \centering
    \includegraphics[keepaspectratio, scale=0.12]{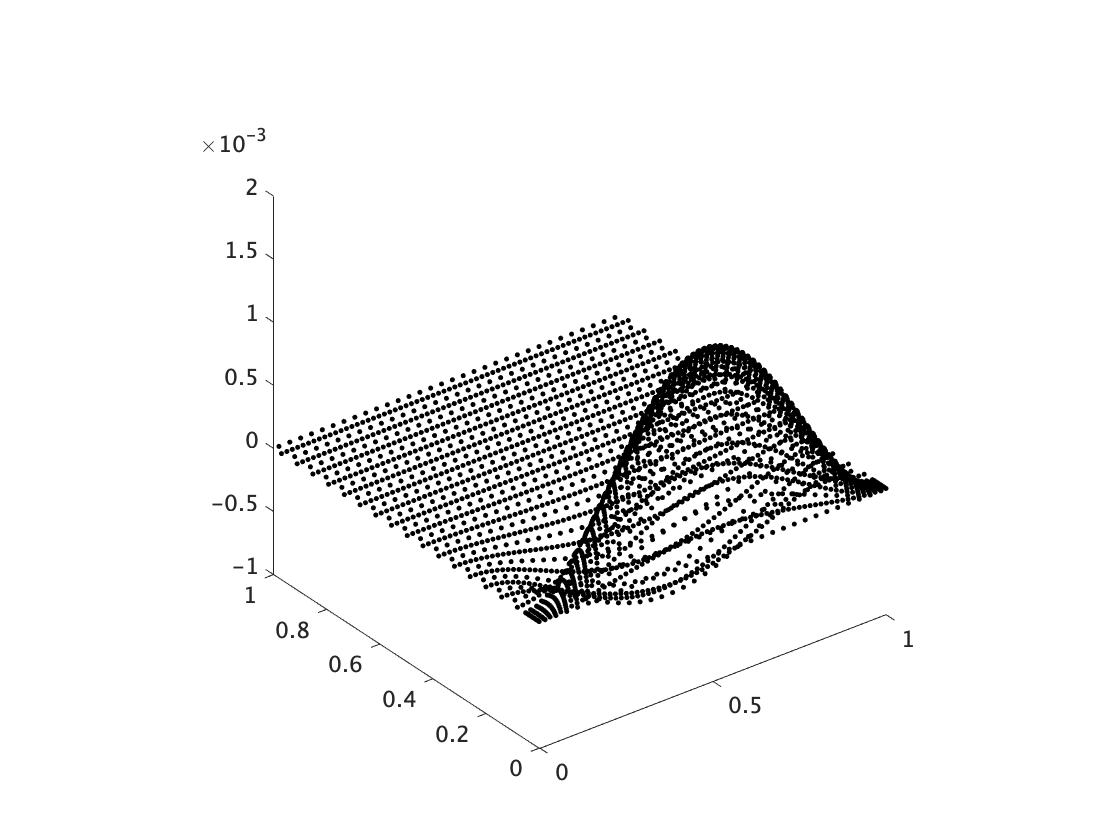}
    \subcaption{Plot of $u_1$}
     \label{fig25}
  \end{minipage}
  \begin{minipage}[b]{0.3\linewidth}
    \centering
   \includegraphics[keepaspectratio, scale=0.12]{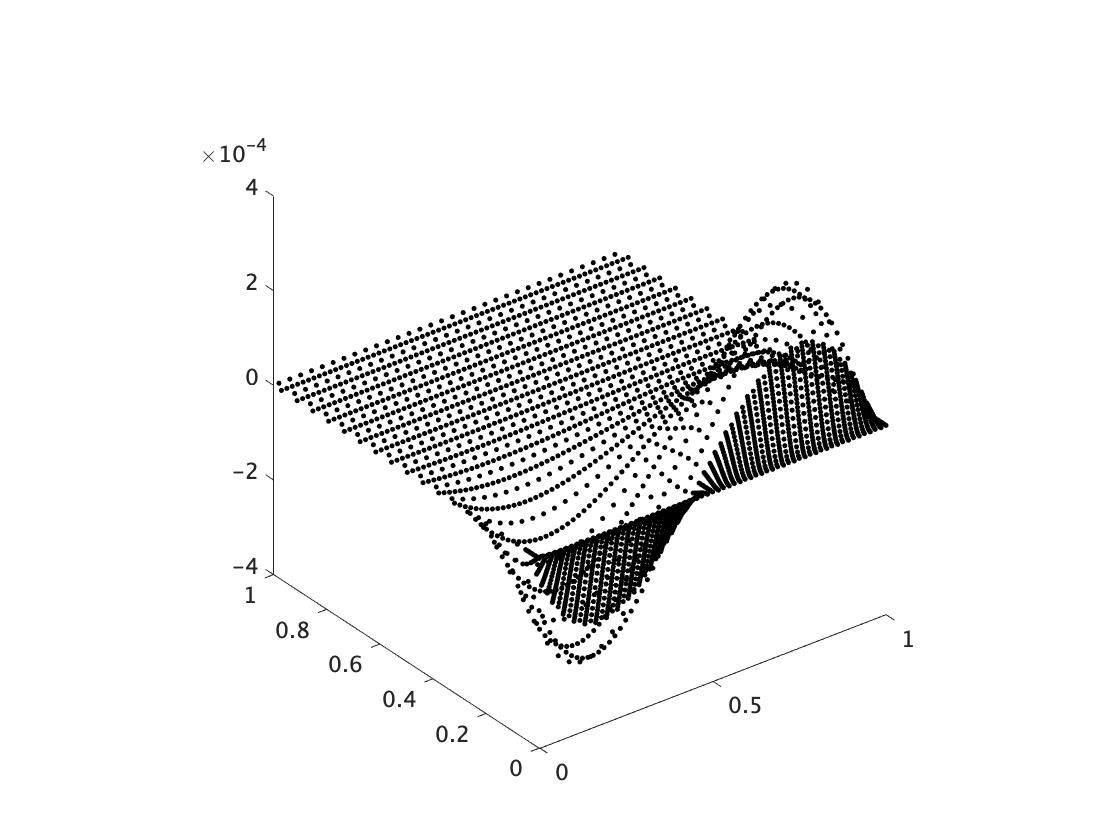}
    \subcaption{Plot of $u_2$}
     \label{fig26}
  \end{minipage}
  \begin{minipage}[b]{0.3\linewidth}
    \centering
  \includegraphics[keepaspectratio, scale=0.12]{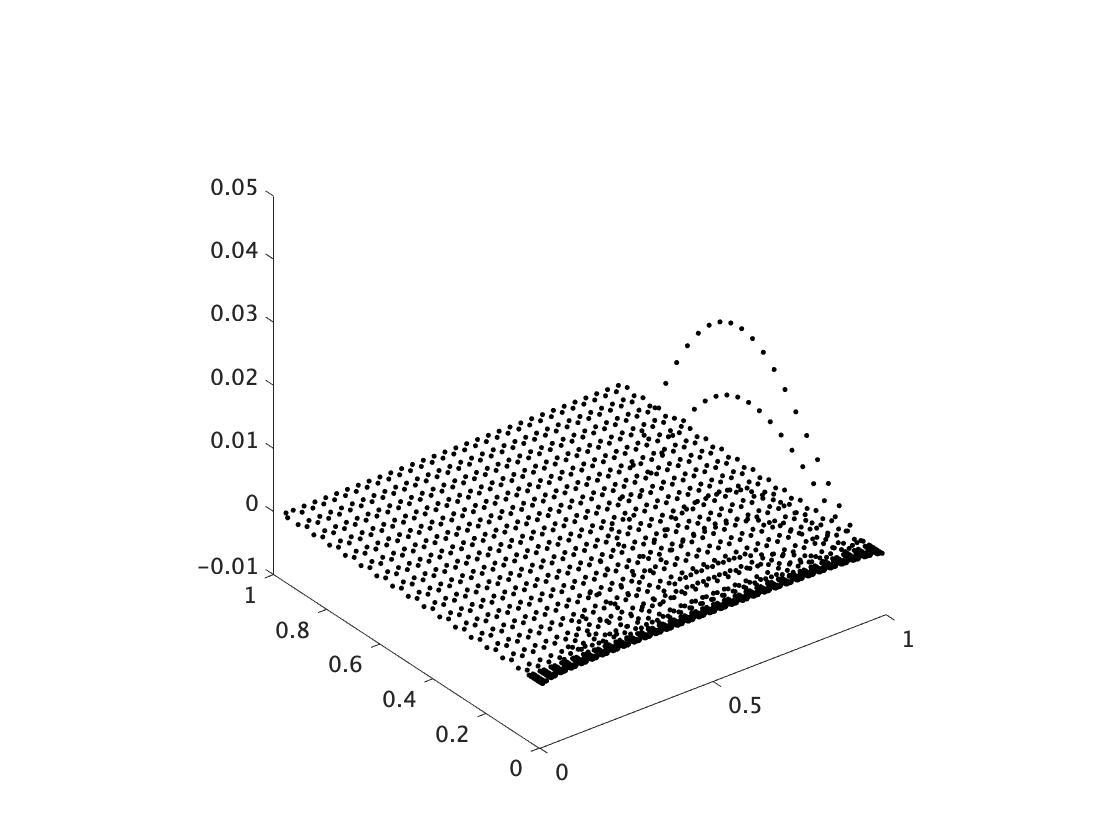}
    \subcaption{Plot of $p$}
     \label{fig27}
  \end{minipage}
  \caption{Mesh (\Roman{ltwo}) Plot of the exact solution, $\eta = \frac{1}{16}$, $\delta = \frac{1}{256}$}
  \label{fexfig=test2}
\end{figure}

 {
We set $\varphi (x_1,x_2) := x_1^2(x_1-1)^2x_2^2(x_2-1)^2 \exp \left( - \frac{x_2}{\eta} \right)$ with $\eta := \sqrt{\delta}$ and  $\delta \in \left\{ \frac{1}{64}, \frac{1}{128}, \frac{1}{256} \right \}$. The function $f$ of the Stokes equation \eqref{intro1} with $\nu = 1$ is given so that the exact solution is
\begin{align*}
\displaystyle
\begin{pmatrix}
 u_1  \\
  u_2
\end{pmatrix}
=
\begin{pmatrix}
 \frac{\partial \varphi}{\partial x_2}  \\
  - \frac{\partial \varphi}{\partial x_1} 
\end{pmatrix},
\quad p =  x_1^2(x_1-1)^2 \exp \left( - \frac{x_2}{ \delta} \right) - \frac{\delta}{30} + \frac{\delta}{30} \exp \left( - \frac{1}{\delta} \right).
\end{align*}
Note that $\int_{\Omega} p dx \approx 0$. See Fig. \ref{fexfig=test2}.
}
{
\begin{rem}
Through a simple calculation,
\begin{align*}
\displaystyle
\frac{\partial^2 u_1}{\partial x_2^2}
&=  \left( - \frac{1}{\eta^3} x_1^2(x_1-1)^2x_2^2(x_2-1)^2  + C_1\left( x_1 , x_2, \eta \right) \right) \exp \left( - \frac{x_2}{\eta} \right),\\
\frac{\partial^2 u_2}{\partial x_2^2}
&=  - \frac{2}{\eta^2} \left( x_1(x_1-1)^2  +  x_1^2(x_1-1)  \right) x_2^2 (x_2-1)^2 \exp \left( - \frac{x_2}{\eta} \right)\\
&\quad + C_2 \left( x_1 , x_2, \eta \right)  \exp \left( - \frac{x_2}{\eta} \right),\\
\frac{\partial p}{\partial x_2}
&= - \frac{1}{\delta} x_1^2(x_1-1)^2 \exp \left( - \frac{x_2}{ \delta} \right),
\end{align*}
where $C_1\left( x_1 , x_2, \eta \right)$ and $C_2\left( x_1,x_2, \eta \right)$ are functions. In this numerical experiment, we set  $\eta = \sqrt{\delta}$.
\end{rem}
}

In the computation, the following schemes were used:
\begin{description}
  \item[(1)] Scheme \eqref{wop=1}  (WOPSIP method).
  \item[(2)] Scheme \eqref{stokes51} (WBCR method).
\end{description}
Let $N \in \{ 16,32,64,128,256 \}$ be the division number for each side of the bottom and height edges of $\Omega$. We consider meshes (\Roman{lone}) and (\Roman{ltwo}) with $\delta \in \left\{ \frac{1}{64},\frac{1}{128}, \frac{1}{256} \right \}$. The notation $\# Np$ denotes the number of nodal points on $V_{dc,h}^{CR} \times Q_h^{0}$ or $V_{h0}^{CR} \times Q_h^{0}$ including nodal points on the boundary. We further use the CG method without preconditioners and the quadrature of the five orders for the computation of the right-hand side in \eqref{wop=1a} and \eqref{stokes51a}.

\begin{table}[htb]
\caption{Mesh (\Roman{lone}) WOPSIP method, $\eta = \frac{1}{8}$, $\delta = \frac{1}{64}$}
\centering
\begin{tabular}{l|l | l || l |l||l|l} \hline
$N$ & $\# Np$&  $h $ & $E_{u_h}^{(\Roman{lone})}$ & $r$ & $E_{p_h}^{(\Roman{lone})}$ & $r$  \\ \hline \hline
16& 3,584 &  8.84e-02 &6.84774e-01   &   &  8.83176e-01    &        \\
32& 14,336 &  4.42e-02  & 3.81183e-01  &  0.85 &  5.23969e-01    & 0.75  \\
64 &  57,344 &   2.21e-02  &1.98511e-01   & 0.94  & 2.72238e-01    &   0.94  \\
128 &  229,376 &   1.10e-02    & 1.00464e-01  &  0.98 &  1.37158e-01  &  0.99       \\
256 &  917,504 &   5.52e-03  & 5.03923e-02  & 1.00  & 6.86957e-02    &  1.00  \\
\hline
\end{tabular}
\label{newtable1_64=1}
\end{table}

\begin{table}[htb]
\caption{Mesh (\Roman{ltwo}) WOPSIP method, $\eta = \frac{1}{8}$, $\delta = \frac{1}{64}$}
\centering
\begin{tabular}{l|l | l || l |l||l|l} \hline
$N$ & $\# Np$&  $h $ & ${E_{u_h}^{(\Roman{ltwo})}}$ & $r$ & ${E_{p_h}^{(\Roman{ltwo})}}$ & $r$  \\ \hline \hline
16& 3,584 &  8.84e-02 & 4.99659e-01 &   &6.39849e-01     &        \\
32& 14,336 &  4.42e-02  &2.66696e-01   & 0.91   & 3.41427e-01   &  0.91 \\
64 &  57,344 &   2.21e-02  &1.42234e-01 & 0.91  & 1.83575e-01    &   0.90  \\
128 &  229,376 &   1.10e-02    & 7.58869e-02  & 0.91  &  9.90536e-02  &  0.89       \\
256 &  917,504 &   5.52e-03  & 4.04865e-02 & 0.91  &  5.34495e-02   &  0.89  \\
\hline
\end{tabular}
\label{newtable1_64=2}
\end{table}

\begin{table}[htb]
\caption{Mesh (\Roman{lone}) WOPSIP method, $\eta = \frac{1}{8 \sqrt{2}}$, $\delta = \frac{1}{128}$}
\centering
\begin{tabular}{l|l | l || l |l||l|l} \hline
$N$ & $\# Np$&  $h $ & $E_{u_h}^{(\Roman{lone})}$ & $r$ & $E_{p_h}^{(\Roman{lone})}$ & $r$  \\ \hline \hline
16& 3,584 &  8.84e-02 &8.23155e-01   &   &  1.05497   &        \\
32& 14,336 &  4.42e-02  & 4.89480e-01   &  0.75 &  8.08138e-01   &  0.38   \\
64 &  57,344 &   2.21e-02  & 2.70311e-01  &  0.86 &   4.90630e-01     &  0.72   \\
128 &  229,376 &   1.10e-02    & 1.41134e-01 &  0.94 &    2.57758e-01     &  0.93       \\
256 &  917,504 &   5.52e-03  & 7.15209e-02 & 0.98  & 1.30291e-01       &  0.98   \\
\hline
\end{tabular}
\label{newtable1}
\end{table}

\begin{table}[htb]
\caption{Mesh (\Roman{ltwo}) WOPSIP method, $\eta = \frac{1}{8 \sqrt{2}}$, $\delta = \frac{1}{128}$}
\centering
\begin{tabular}{l|l | l || l |l||l|l} \hline
$N$ & $\# Np$&  $h $ & $E_{u_h}^{(\Roman{ltwo})}$ & $r$ & $E_{p_h}^{(\Roman{ltwo})}$ & $r$  \\ \hline \hline
16& 3,584 &  1.30e-01 & 5.97426e-01 &   &    9.36325e-01    &          \\
32& 14,336 &  6.39e-02  & 3.13359e-01 & 0.93  &   4.83980e-01    & 0.95   \\
64 &  57,344 &   3.14e-02  & 1.62259e-01  &   0.95 &   2.51145e-01    &  0.95    \\
128 &  229,376 &   1.54e-02    & 8.35607e-02 &  0.96 & 1.30872e-01       &  0.94     \\
256 &  917,504 &   7.55e-03  & 4.30401e-02  &   0.96 &    6.83930e-02    &  0.94  \\
\hline
\end{tabular}
\label{newtable2}
\end{table}

\begin{table}[htb]
\caption{Mesh (\Roman{lone}) WOPSIP method, $\eta = \frac{1}{16}$, $\delta = \frac{1}{256}$}
\centering
\begin{tabular}{l|l | l || l |l||l|l} \hline
$N$ & $\# Np$&  $h $ & $E_{u_h}^{(\Roman{lone})}$ & $r$ & $E_{p_h}^{(\Roman{lone})}$ & $r$  \\ \hline \hline
16& 3,584 &  8.84e-02 & 9.50427e-01   &   & 1.12990      &      \\
32& 14,336 &  4.42e-02  & 6.04935e-01 &  0.65 &    9.94423e-01         & 0.18  \\
64 &  57,344 &   2.21e-02  & 3.49110e-01  &  0.79 &    7.71574e-01    &    0.37  \\
128 &  229,376 &   1.10e-02    & 1.94504e-01 &  0.84 &  4.72726e-01    &   0.71 \\
256 &  917,504 &   5.52e-03  &  1.02293e-01  &  0.93 &    2.49118e-01    &    0.92  \\
\hline
\end{tabular}
\label{newtable3}
\end{table}

\begin{table}[htb]
\caption{Mesh (\Roman{ltwo}) WOPSIP method, $\eta = \frac{1}{16}$, $\delta = \frac{1}{256}$}
\centering
\begin{tabular}{l|l | l || l |l||l|l} \hline
$N$ & $\# Np$&  $h $ & $E_{u_h}^{(\Roman{ltwo})}$ & $r$ & $E_{p_h}^{(\Roman{ltwo})}$ & $r$  \\ \hline \hline
16& 3,584 & 1.35e-01  &7.89295e-01    &   &1.46492        &        \\
32& 14,336 &   6.69e-02  & 4.10272e-01  & 0.94  &  7.48414e-01      & 0.97    \\
64 &  57,344 & 3.31e-02     & 2.11273e-01 &  0.96 & 3.81425e-01     &  0.97  \\
128 &  229,376 & 1.64e-02     & 1.07657e-01  &  0.97 &   1.94667e-01    &  0.97       \\
256 &  917,504 &  8.13e-03   & 5.45969e-02  &  0.98 &   9.95082e-02    & 0.97   \\
\hline
\end{tabular}
\label{newtable4}
\end{table}

\begin{table}[htb]
\caption{Mesh (\Roman{lone}) WBCR method, $\eta = \frac{1}{8}$, $\delta = \frac{1}{64}$}
\centering
\begin{tabular}{l|l | l || l |l||l|l} \hline
$N$ & $\# Np$&  $h $ & $E_{u_h}^{(\Roman{lone})}$ & $r$ & $E_{p_h}^{(\Roman{lone})}$ & $r$  \\ \hline \hline
16& 3,584 &  8.84e-02 &   9.81333e-01  &   &  1.38484      &        \\
32& 14,336 &  4.42e-02  &5.23575e-01   & 0.91  &  6.14362e-01       & 1.17   \\
64 &  57,344 &   2.21e-02  & 2.65825e-01&  0.98 &    2.84631e-01    &  1.11 \\
128 &  229,376 &   1.10e-02    &  1.33413e-01  & 0.99  &  1.38739e-01      & 1.04 \\
256 &  917,504 &   5.52e-03  & 6.67689e-02 &  1.00 &     6.88943e-02   & 1.01 \\
\hline
\end{tabular}
\label{newtable1_64=3}
\end{table}

\begin{table}[htb]
\caption{Mesh (\Roman{ltwo}) WBCR method, $\eta = \frac{1}{8}$, $\delta = \frac{1}{64}$}
\centering
\begin{tabular}{l|l | l || l |l||l|l} \hline
$N$ & $\# Np$&  $h $ & ${E_{u_h}^{(\Roman{ltwo})}}$ & $r$ & ${E_{p_h}^{(\Roman{ltwo})}}$ & $r$  \\ \hline \hline
16& 3,584 &  8.84e-02 & 7.58108e-01     &   &  7.22029e-01       &        \\
32& 14,336 &  4.42e-02  & 3.93515e-01 &  0.95 &      3.58341e-01   & 1.01  \\
64 &  57,344 &   2.21e-02  &2.04706e-01   &  0.94 &  1.86657e-01      & 0.94  \\
128 &  229,376 &   1.10e-02    & 1.06968e-01 &  0.94 &    9.95828e-02    & 0.91 \\
256 &  917,504 &   5.52e-03  & 5.60719e-02 &  0.93 &    5.35364e-02     & 0.90  \\
\hline
\end{tabular}
\label{newtable1_64=4}
\end{table}

\begin{table}[htb]
\caption{Mesh (\Roman{lone}) WBCR method, $\eta = \frac{1}{8 \sqrt{2}}$, $\delta = \frac{1}{128}$}
\centering
\begin{tabular}{l|l | l || l |l||l|l} \hline
$N$ & $\# Np$&  $h $ & $E_{u_h}^{(\Roman{lone})}$ & $r$ & $E_{p_h}^{(\Roman{lone})}$ & $r$  \\ \hline \hline
16& 3,584 &  8.84e-02 &  1.26704    &   &  1.75263      &        \\
32& 14,336 &  4.42e-02  & 7.05425e-01  &  0.84 &    9.54909e-01    &  0.88   \\
64 &  57,344 &   2.21e-02  & 3.61245e-01  &  0.97 & 5.12135e-01      &   0.90  \\
128 &  229,376 &   1.10e-02    &  1.81656e-01 & 0.99  &   2.60500e-01     &  0.98       \\
256 &  917,504 &   5.52e-03  & 9.09561e-02  & 1.00  &     1.30634e-01   &   1.00  \\
\hline
\end{tabular}
\label{newtable5}
\end{table}

\begin{table}[htb]
\caption{Mesh (\Roman{ltwo}) WBCR method, $\eta = \frac{1}{8 \sqrt{2}}$, $\delta = \frac{1}{128}$}
\centering
\begin{tabular}{l|l | l || l |l||l|l} \hline
$N$ & $\# Np$&  $h $ & $E_{u_h}^{(\Roman{ltwo})}$ & $r$ & $E_{p_h}^{(\Roman{ltwo})}$ & $r$  \\ \hline \hline
16& 3,584 &  1.30e-01 &9.89743e-01  &   &    1.02200    &        \\
32& 14,336 &  6.39e-02  & 5.02759e-01   &   0.98&   4.96331e-01      &  1.04\\
64 &  57,344 &   3.14e-02  &2.53831e-01   &  0.99 &2.53006e-01        &   0.97  \\
128 &  229,376 &   1.54e-02    & 1.28052e-01   & 0.99  &   1.31165e-01    &  0.95    \\
256 &  917,504 &   7.55e-03  & 6.47475e-02 &  0.98 & 6.84384e-02       & 0.94   \\
\hline
\end{tabular}
\label{newtable6}
\end{table}

\begin{table}[htb]
\caption{Mesh (\Roman{lone}) WBCR method, $\eta = \frac{1}{16}$, $\delta = \frac{1}{256}$}
\centering
\begin{tabular}{l|l | l || l |l||l|l} \hline
$N$ & $\# Np$&  $h $ & $E_{u_h}^{(\Roman{lone})}$ & $r$ & $E_{p_h}^{(\Roman{lone})}$ & $r$  \\ \hline \hline
16& 3,584 &  8.84e-02 &  1.56033 &   &      2.05430  &      \\
32& 14,336 &  4.42e-02  &  9.44351e-01  &  0.72  &  1.20348      & 0.77 \\
64 &  57,344 &   2.21e-02  & 4.91889e-01 & 0.94  &   8.06744e-01     &   0.58  \\
128 &  229,376 &   1.10e-02    &2.48251e-01   &  0.99 & 4.77567e-01       &  0.76  \\
256 &  917,504 &   5.52e-03  &  1.24408e-01 &  1.00 &    2.49727e-01   &   0.94  \\
\hline
\end{tabular}
\label{newtable7}
\end{table}

\begin{figure}[htbp]
  \begin{minipage}[b]{0.3\linewidth}
    \centering
  \includegraphics[keepaspectratio, scale=0.12]{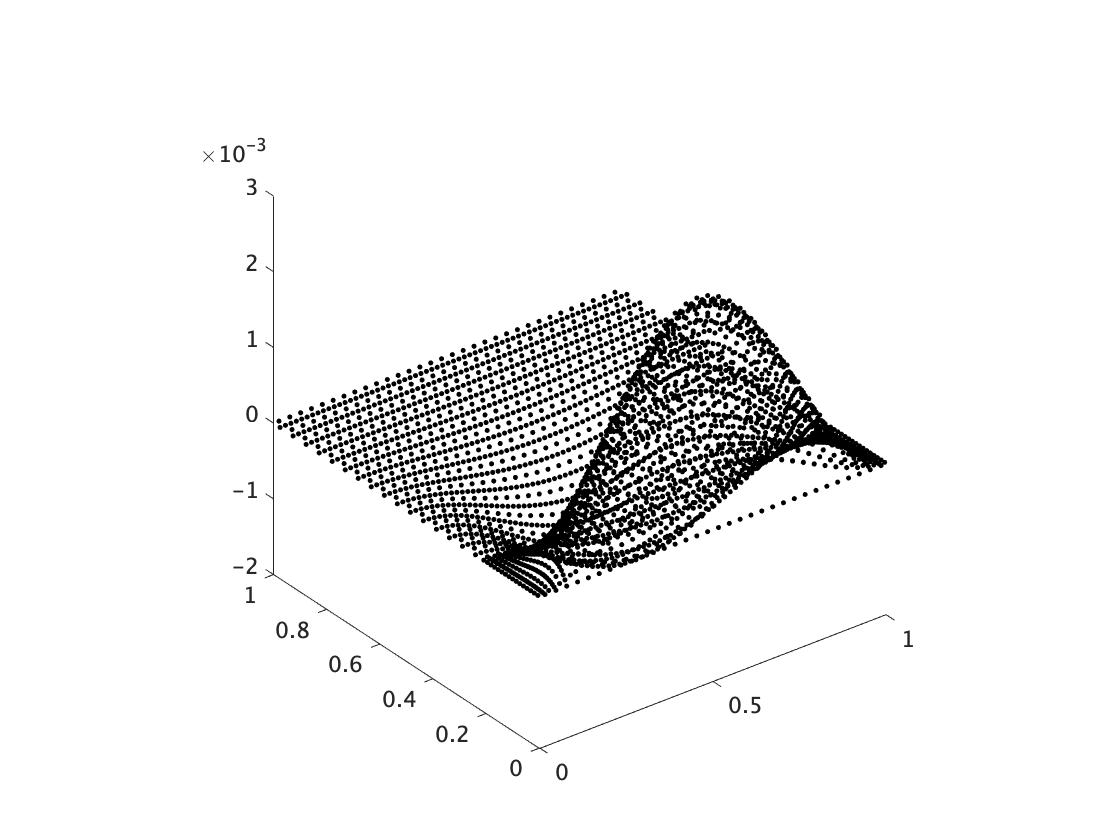}
    \subcaption{Plot of $u_{h,1}$}
     \label{fig25}
  \end{minipage}
  \begin{minipage}[b]{0.3\linewidth}
    \centering
   \includegraphics[keepaspectratio, scale=0.12]{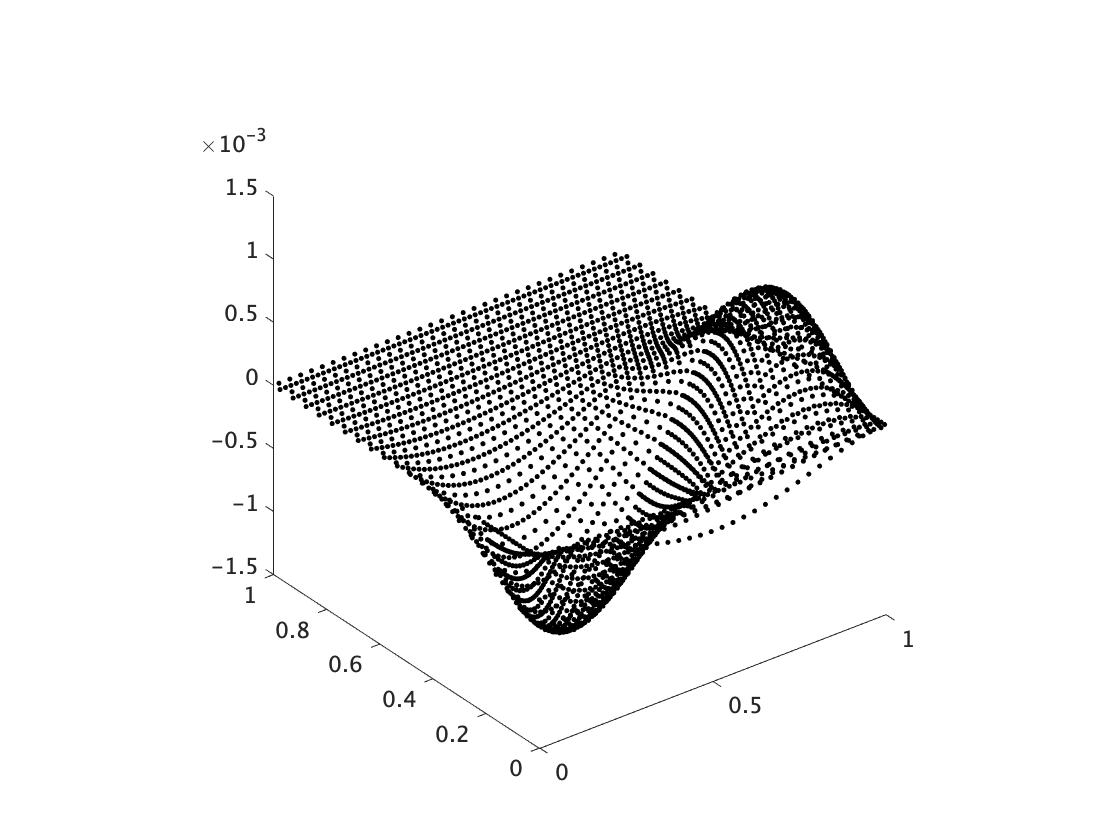}
    \subcaption{Plot of $u_{h,2}$}
     \label{fig26}
  \end{minipage}
  \begin{minipage}[b]{0.3\linewidth}
    \centering
  \includegraphics[keepaspectratio, scale=0.12]{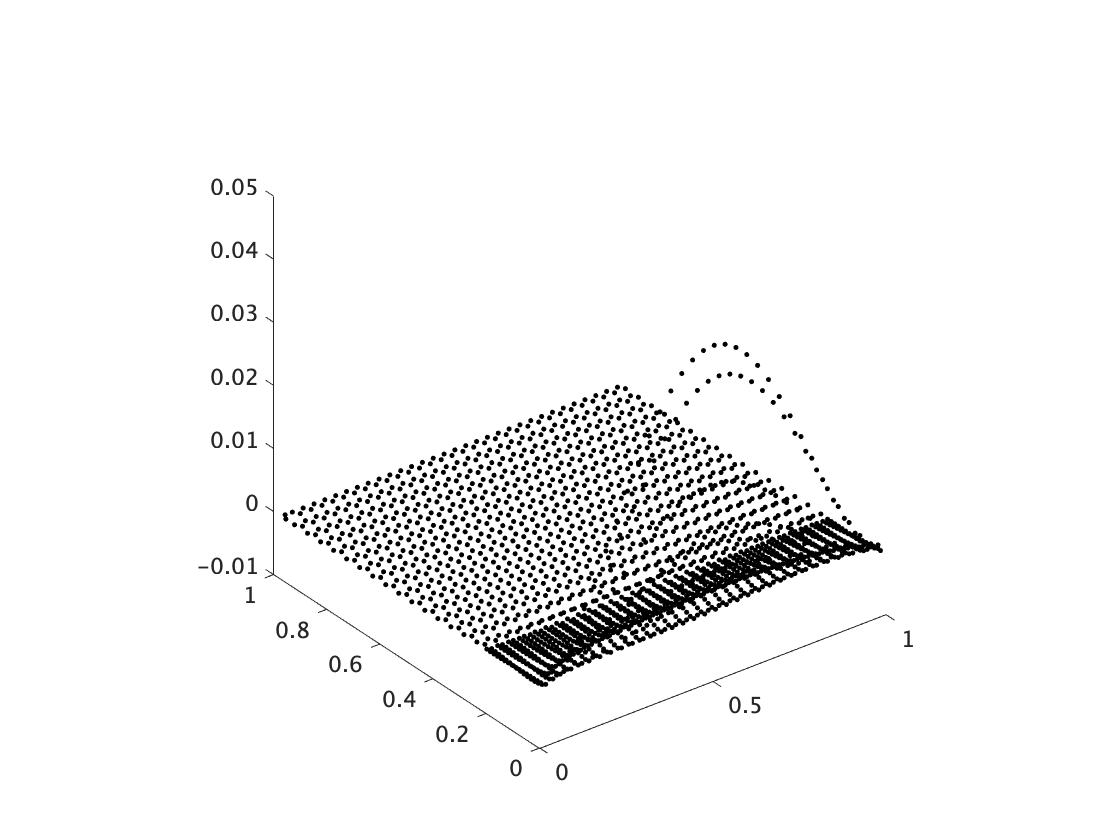}
    \subcaption{Plot of $p_h$}
     \label{fig27}
  \end{minipage}
  \caption{Mesh (\Roman{ltwo}) Plot of the WOPSIP solution, $\eta = \frac{1}{8}$, $\delta = \frac{1}{64}$}
  \label{test2=wop_1_64}
\end{figure}

\begin{figure}[htbp]
  \begin{minipage}[b]{0.3\linewidth}
    \centering
  \includegraphics[keepaspectratio, scale=0.12]{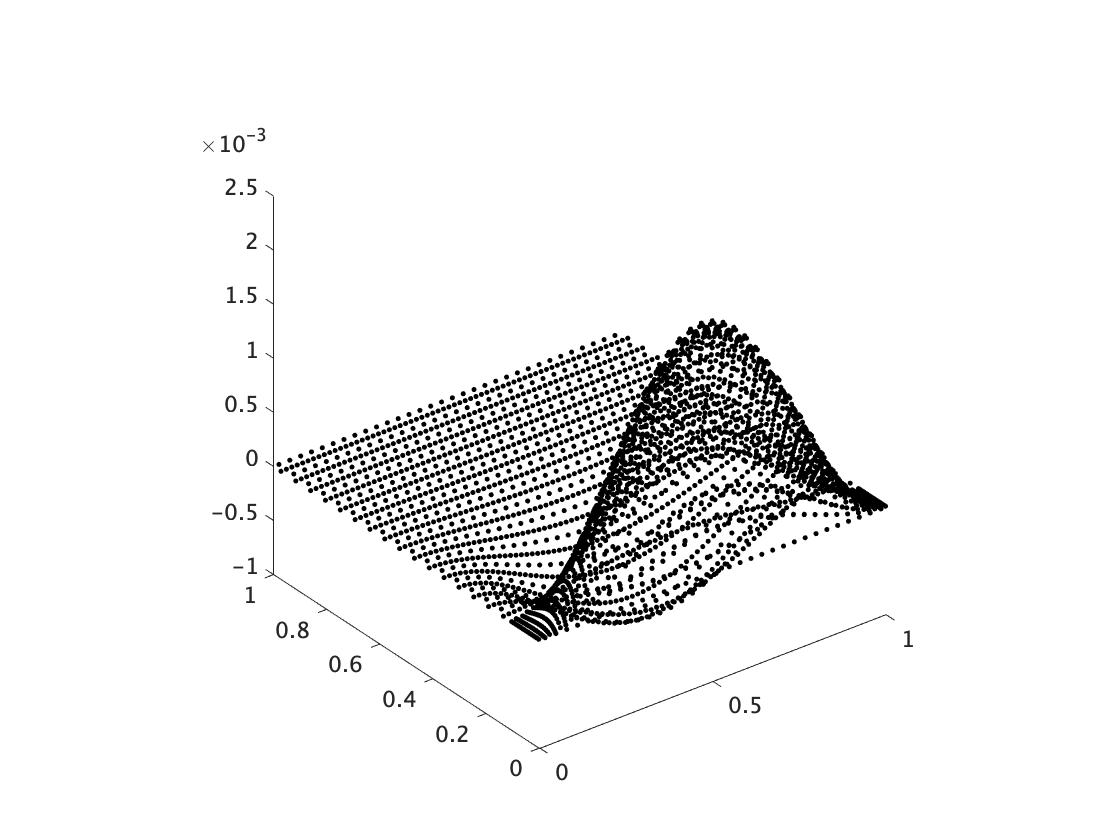}
    \subcaption{Plot of $u_{h,1}$}
     \label{fig25}
  \end{minipage}
  \begin{minipage}[b]{0.3\linewidth}
    \centering
   \includegraphics[keepaspectratio, scale=0.12]{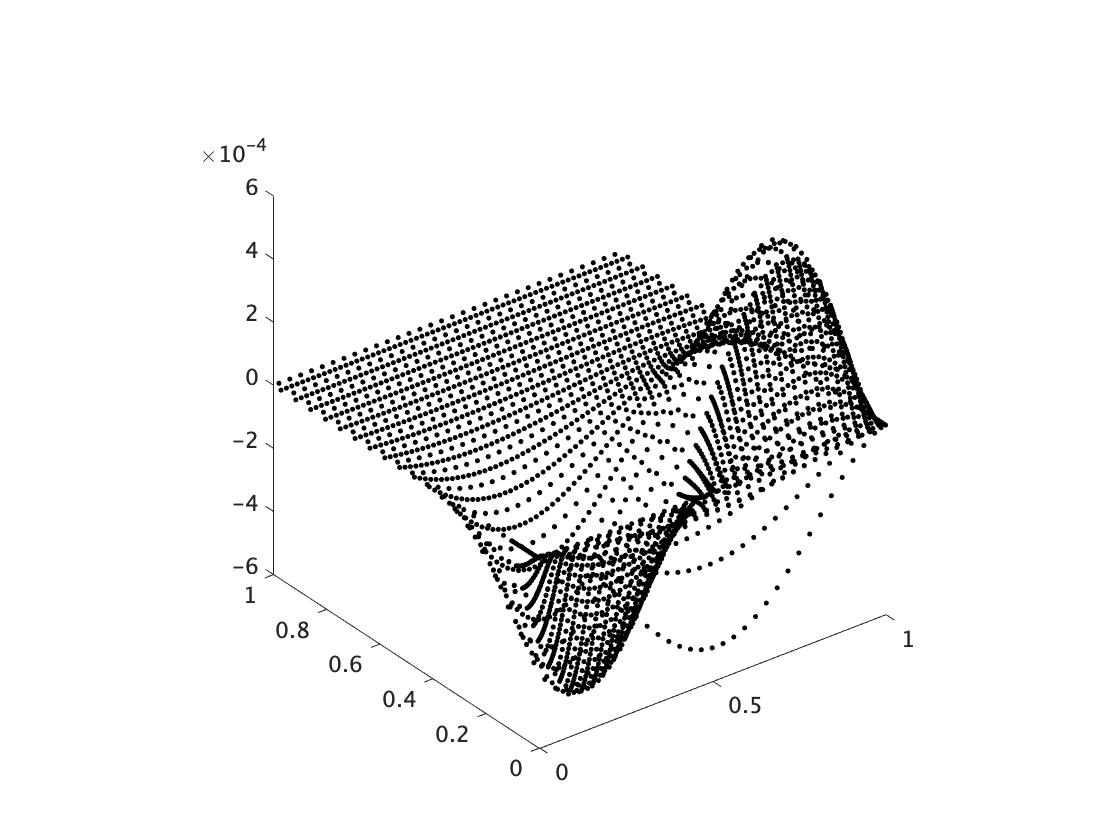}
    \subcaption{Plot of $u_{h,2}$}
     \label{fig26}
  \end{minipage}
  \begin{minipage}[b]{0.3\linewidth}
    \centering
  \includegraphics[keepaspectratio, scale=0.12]{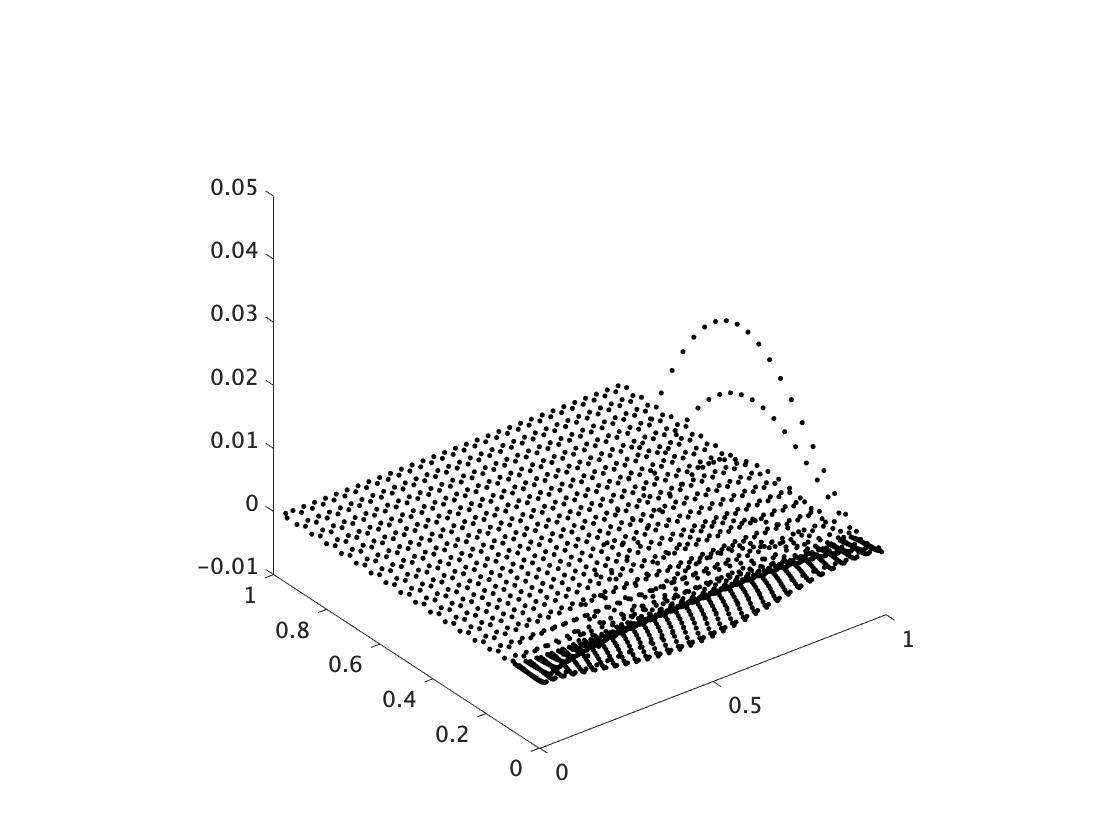}
    \subcaption{Plot of $p_h$}
     \label{fig27}
  \end{minipage}
  \caption{Mesh (\Roman{ltwo}) Plot of the WOPSIP solution, $\eta = \frac{1}{8 \sqrt{2}}$, $\delta = \frac{1}{128}$}
  \label{test2=wop_1_128}
\end{figure}

\begin{figure}[htbp]
  \begin{minipage}[b]{0.3\linewidth}
    \centering
    \includegraphics[keepaspectratio, scale=0.12]{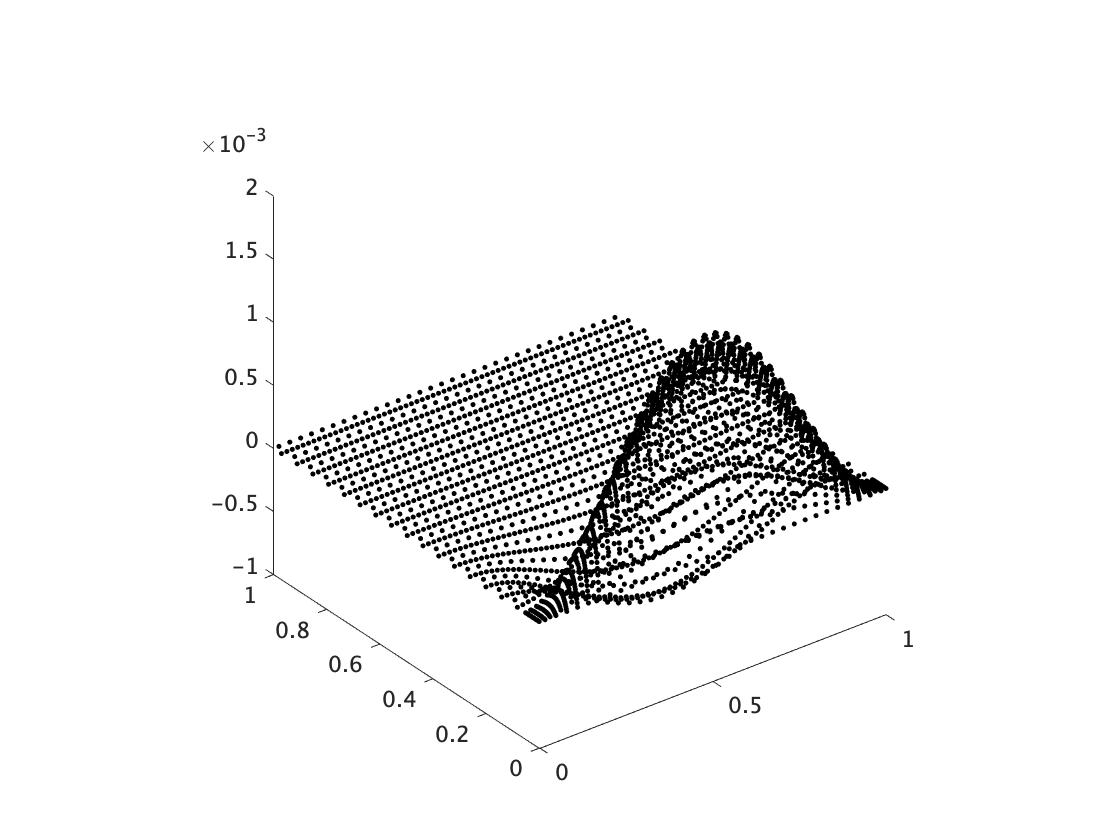}
    \subcaption{Plot of $u_{h,1}$}
     \label{fig25}
  \end{minipage}
  \begin{minipage}[b]{0.3\linewidth}
    \centering
   \includegraphics[keepaspectratio, scale=0.12]{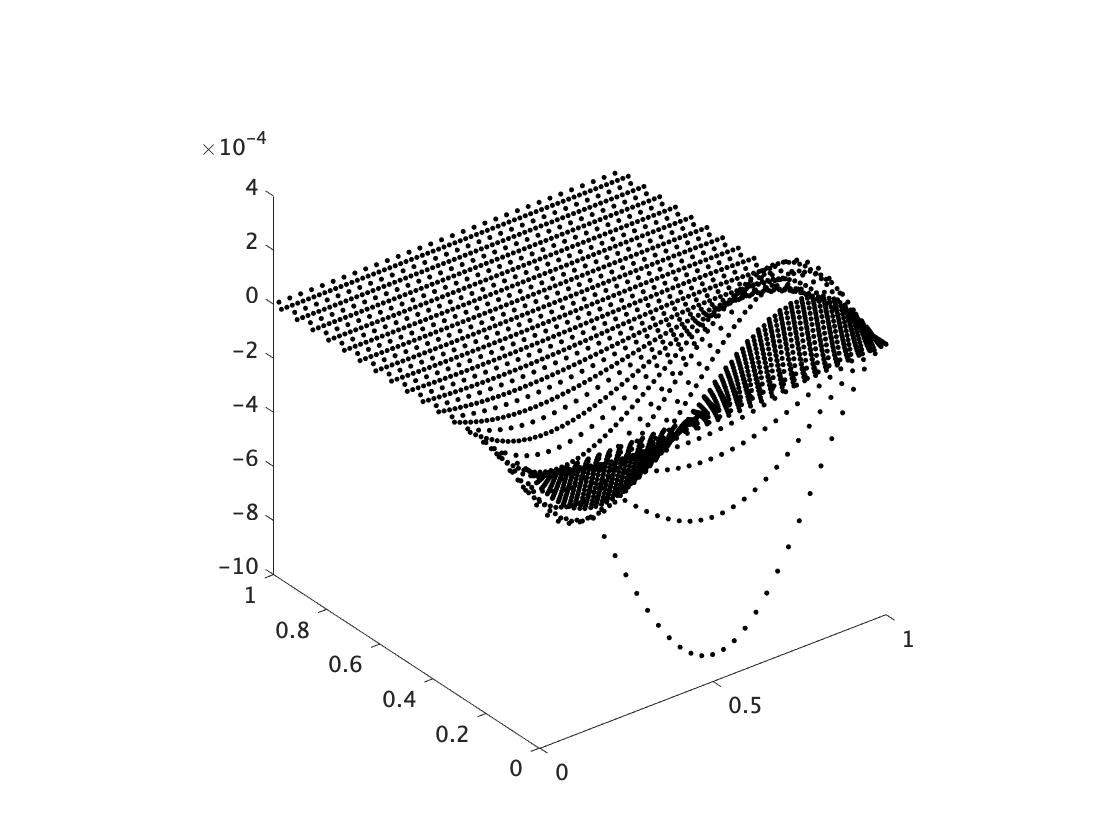}
    \subcaption{Plot of $u_{h,2}$}
     \label{fig26}
  \end{minipage}
  \begin{minipage}[b]{0.3\linewidth}
    \centering
  \includegraphics[keepaspectratio, scale=0.12]{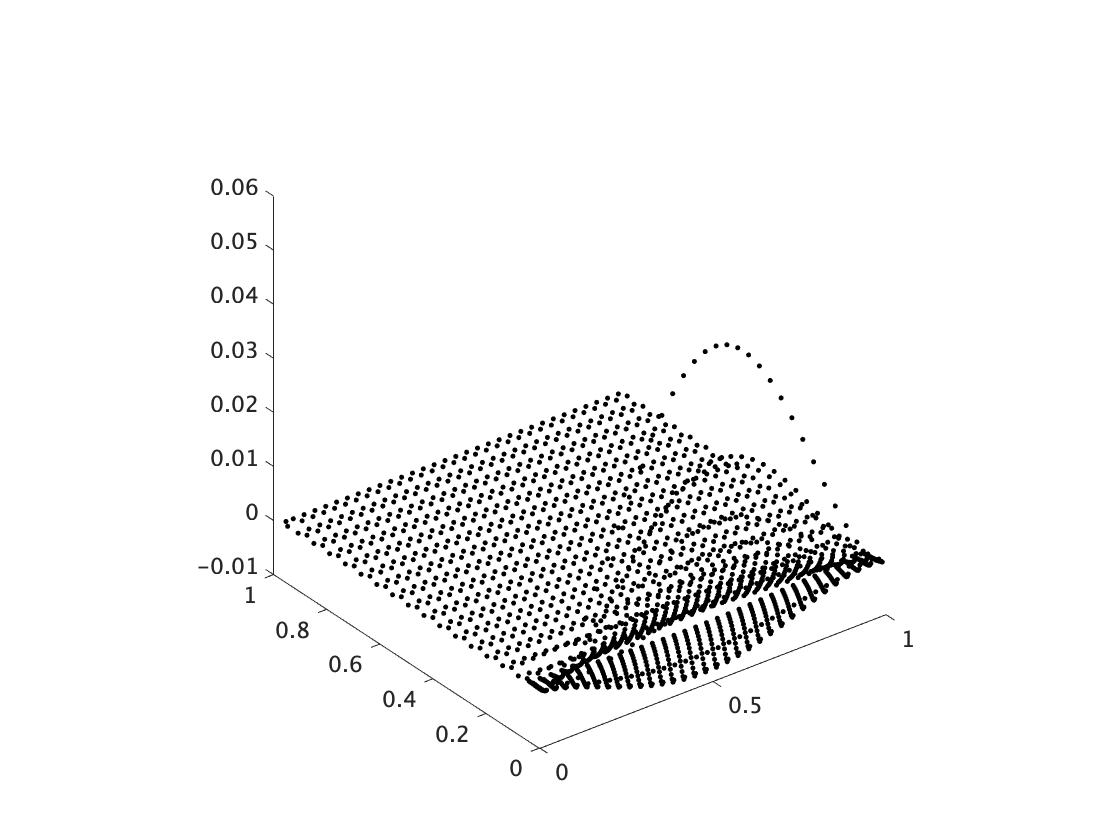}
    \subcaption{Plot of $p_h$}
     \label{fig27}
  \end{minipage}
  \caption{Mesh (\Roman{ltwo}) Plot of the WOPSIP solution, $\eta = \frac{1}{16}$, $\delta = \frac{1}{256}$}
  \label{test2=wop_1_256}
\end{figure}

\begin{figure}[htbp]
  \begin{minipage}[b]{0.3\linewidth}
    \centering
  \includegraphics[keepaspectratio, scale=0.12]{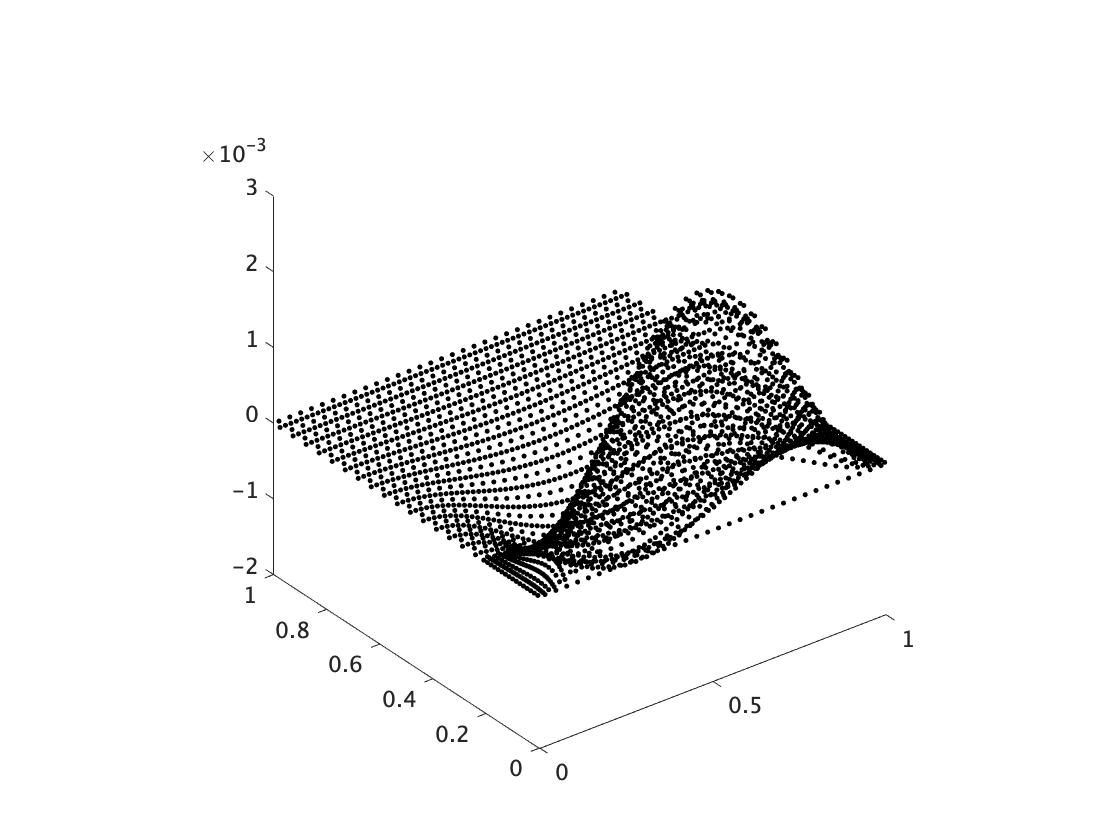}
    \subcaption{Plot of $u_{h,1}$}
     \label{fig25}
  \end{minipage}
  \begin{minipage}[b]{0.3\linewidth}
    \centering
   \includegraphics[keepaspectratio, scale=0.12]{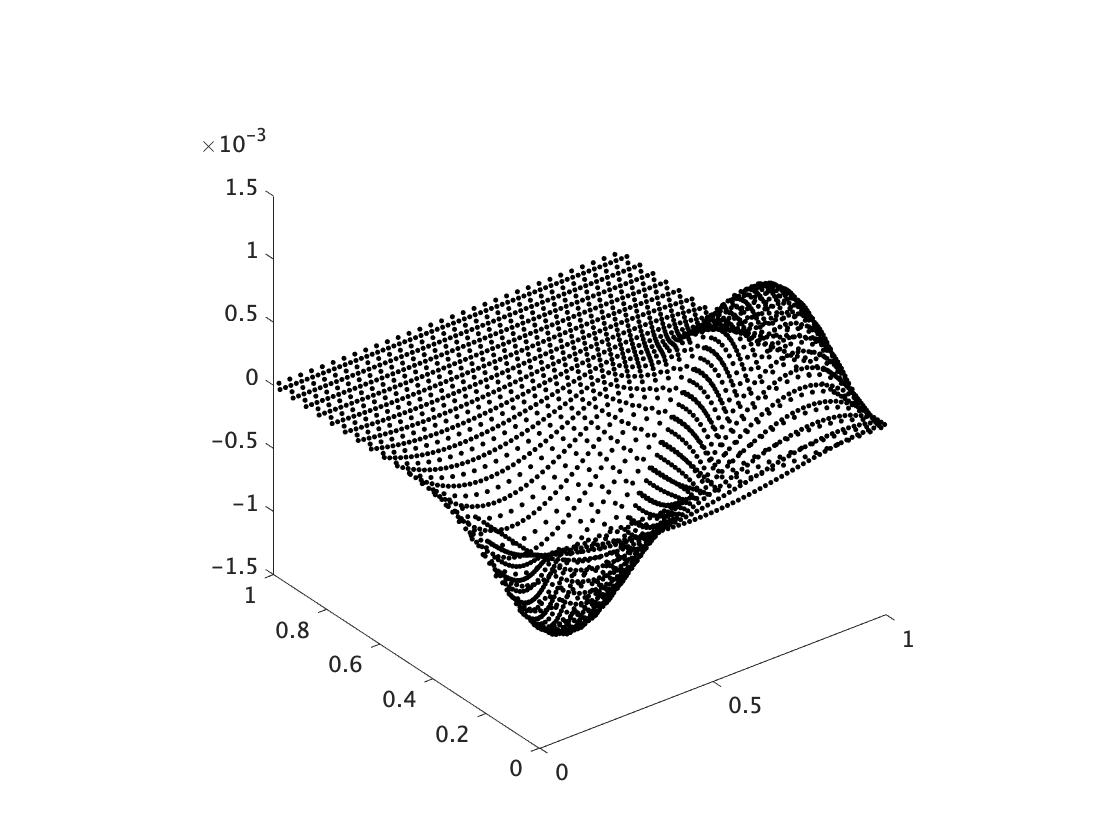}
    \subcaption{Plot of $u_{h,2}$}
     \label{fig26}
  \end{minipage}
  \begin{minipage}[b]{0.3\linewidth}
    \centering
  \includegraphics[keepaspectratio, scale=0.12]{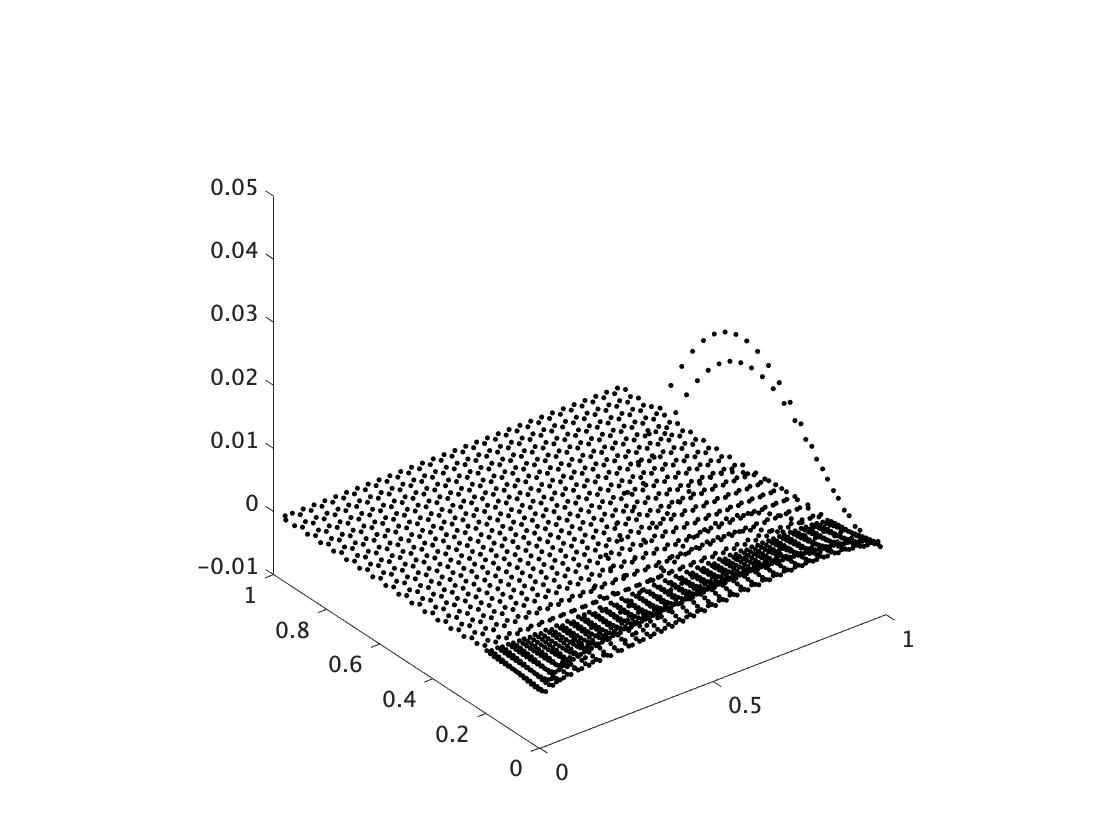}
    \subcaption{Plot of $p_h$}
     \label{fig27}
  \end{minipage}
  \caption{Mesh (\Roman{ltwo}) Plot of the WBCR solution, $\eta = \frac{1}{8}$, $\delta = \frac{1}{64}$}
  \label{test2=wbcr_1_64}
\end{figure}

\begin{figure}[htbp]
  \begin{minipage}[b]{0.3\linewidth}
    \centering
    \includegraphics[keepaspectratio, scale=0.12]{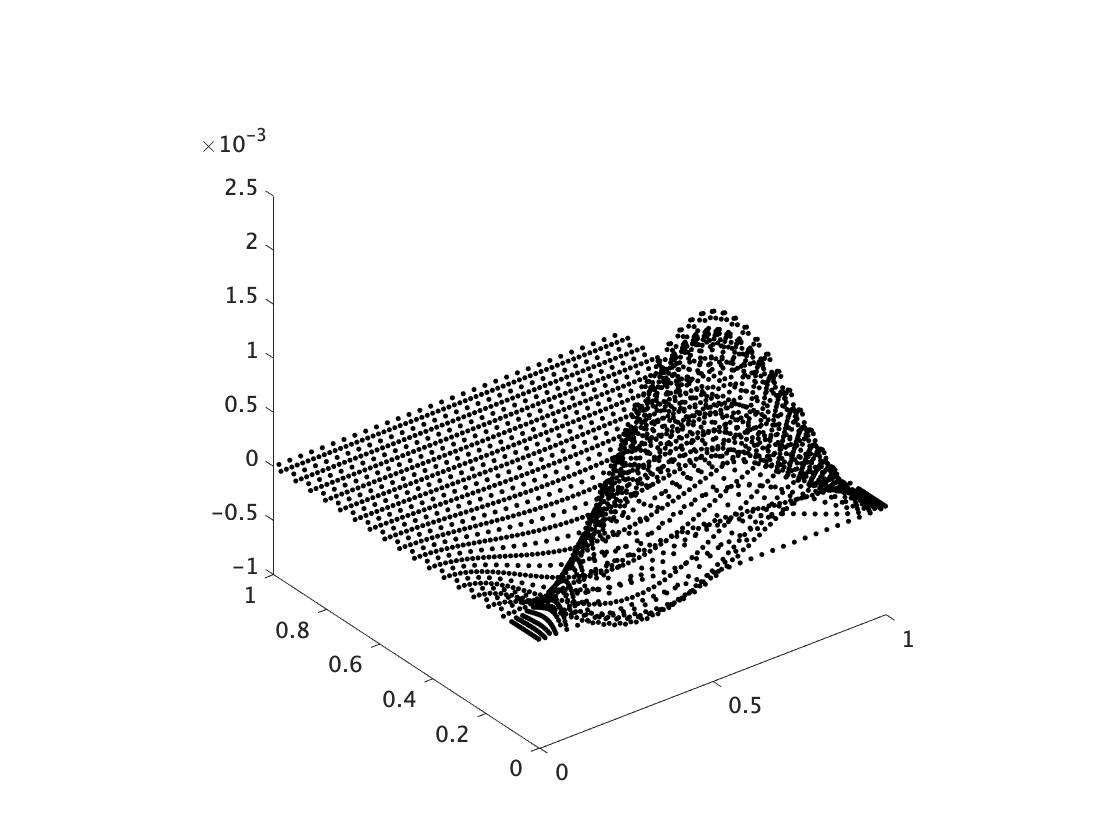}
    \subcaption{Plot of $u_{h,1}$}
     \label{fig25}
  \end{minipage}
  \begin{minipage}[b]{0.3\linewidth}
    \centering
   \includegraphics[keepaspectratio, scale=0.12]{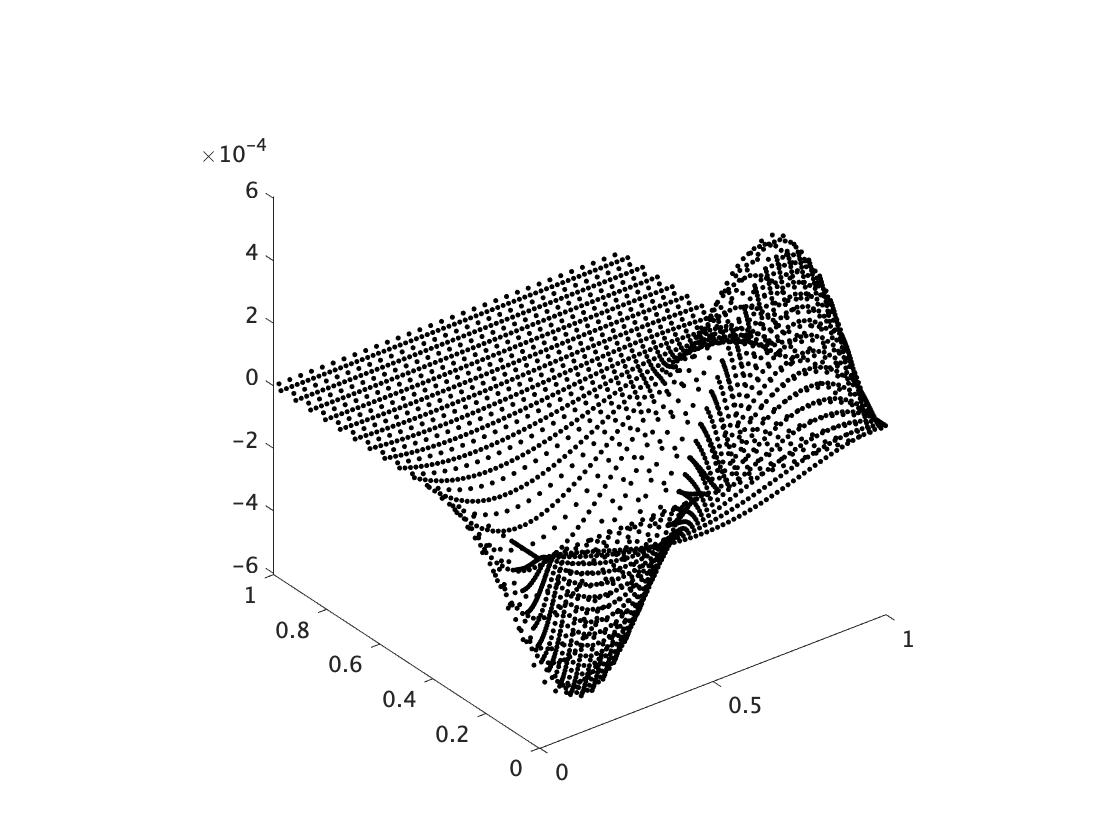}
    \subcaption{Plot of $u_{h,2}$}
     \label{fig26}
  \end{minipage}
  \begin{minipage}[b]{0.3\linewidth}
    \centering
  \includegraphics[keepaspectratio, scale=0.12]{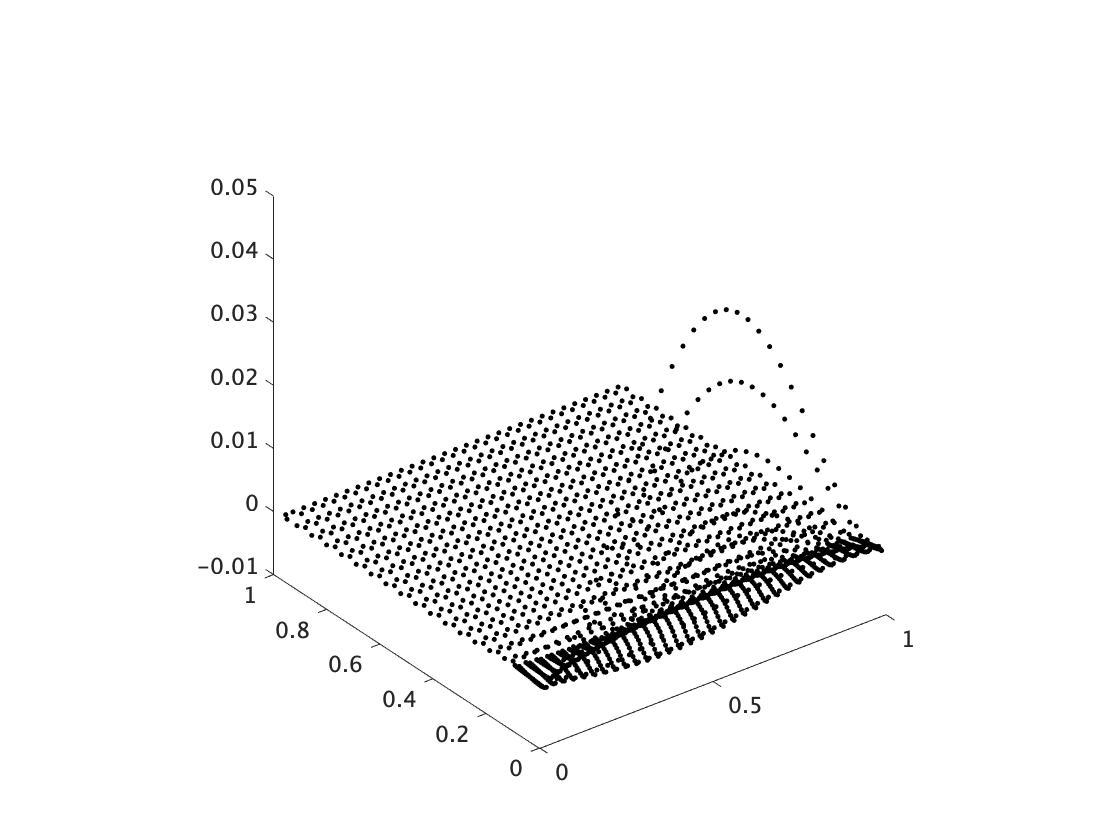}
    \subcaption{Plot of $p_h$}
     \label{fig27}
  \end{minipage}
  \caption{Mesh (\Roman{ltwo}) Plot of the WBCR solution, $\eta = \frac{1}{8 \sqrt{2}}$, $\delta = \frac{1}{128}$}
  \label{test2=wbcr_1_128}
\end{figure}

\begin{figure}[htbp]
  \begin{minipage}[b]{0.3\linewidth}
    \centering
    \includegraphics[keepaspectratio, scale=0.12]{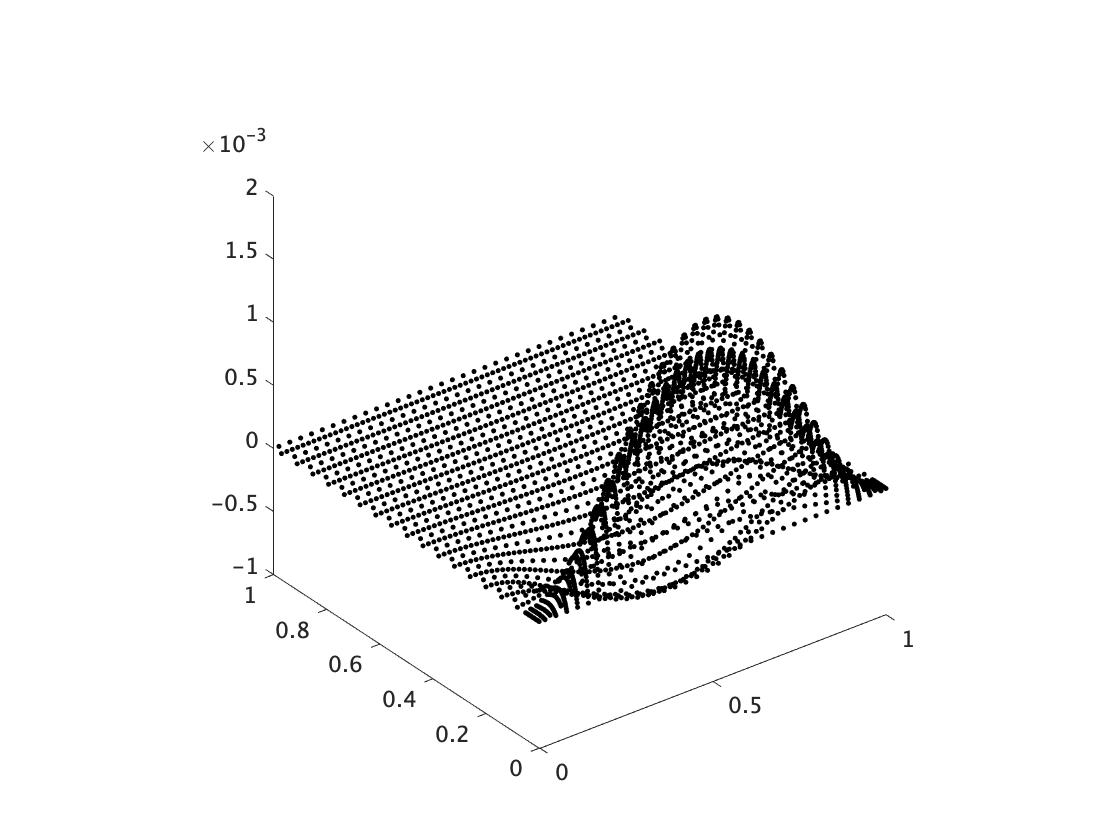}
    \subcaption{Plot of $u_{h,1}$}
     \label{fig25}
  \end{minipage}
  \begin{minipage}[b]{0.3\linewidth}
    \centering
   \includegraphics[keepaspectratio, scale=0.12]{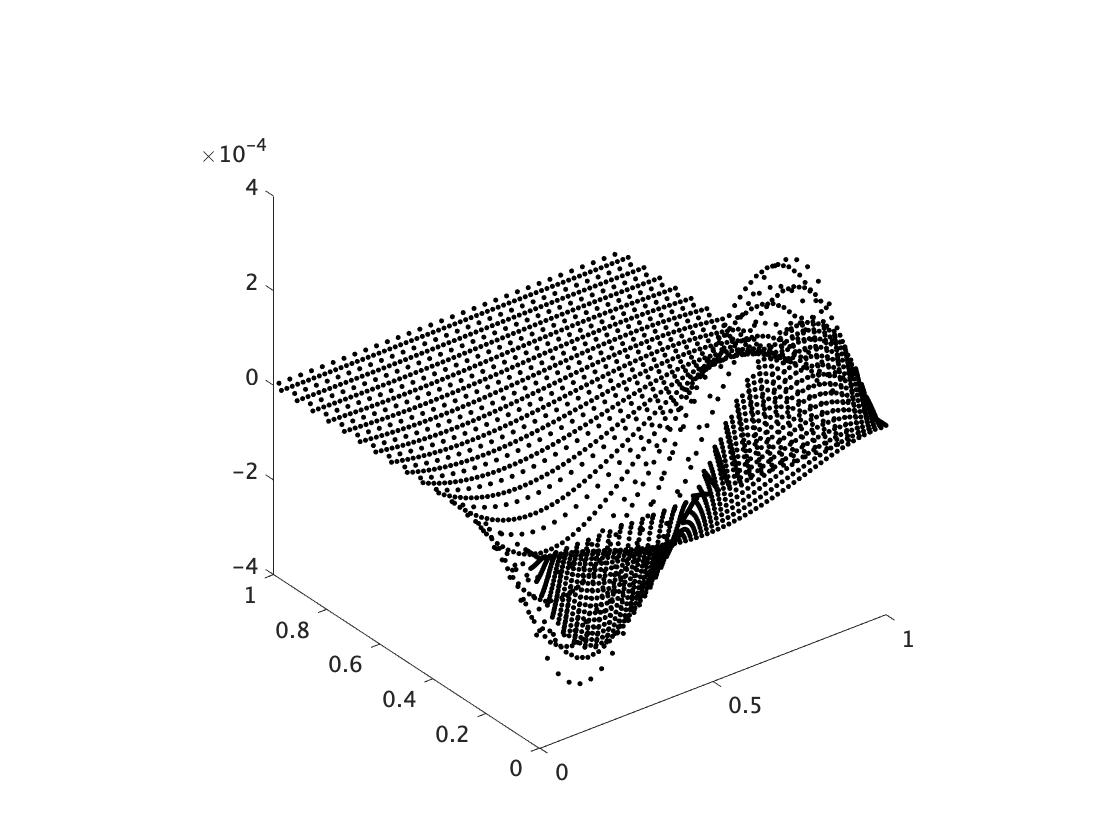}
    \subcaption{Plot of $u_{h,2}$}
     \label{fig26}
  \end{minipage}
  \begin{minipage}[b]{0.3\linewidth}
    \centering
  \includegraphics[keepaspectratio, scale=0.12]{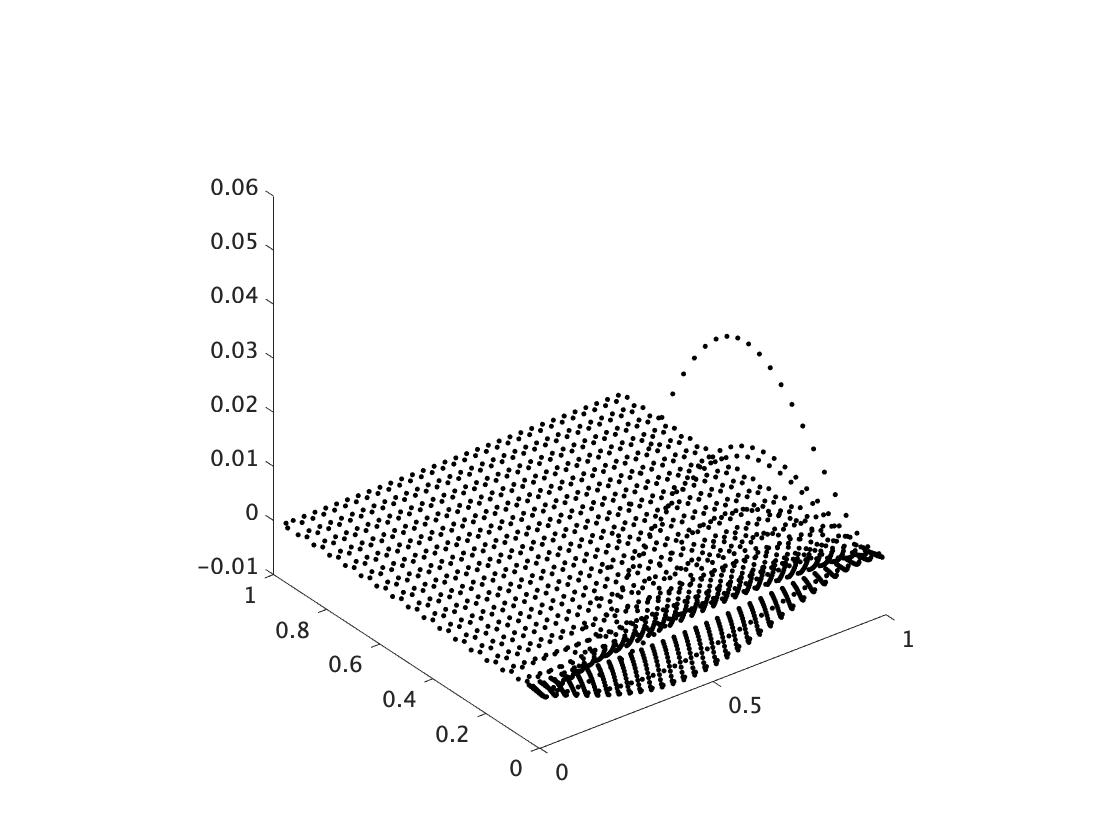}
    \subcaption{Plot of $p_h$}
     \label{fig27}
  \end{minipage}
  \caption{Mesh (\Roman{ltwo}) Plot of the WBCR solution, $\eta = \frac{1}{16}$, $\delta = \frac{1}{256}$}
  \label{test2=wbcr_1_256}
\end{figure}

\begin{table}[thb]
\caption{Mesh (\Roman{ltwo}) WBCR method, $\eta = \frac{1}{16}$, $\delta = \frac{1}{256}$}
\centering
\begin{tabular}{l|l | l || l |l||l|l} \hline
$N$ & $\# Np$&  $h $ & $E_{u_h}^{(\Roman{ltwo})}$ & $r$ & $E_{p_h}^{(\Roman{ltwo})}$ & $r$  \\ \hline \hline
16& 3,584 & 1.35e-01    & 1.32981  &   & 1.78771       &        \\
32& 14,336 &6.69e-02   &  6.72574e-01  &0.98  &   7.85551e-01     &   1.19 \\
64 &  57,344 &  3.31e-02    &  3.38546e-01 &  0.99 &   3.83474e-01      &   1.03 \\
128 &  229,376 & 1.64e-02    & 1.69928e-01 & 0.99  &   1.93935e-01     & 0.98        \\
256 &  917,504 &  8.13e-03   & 8.51702e-02  &  1.00 &    9.87999e-02    &   0.97 \\
\hline
\end{tabular}
\label{newtable8}
\end{table}
{
Tables \ref{newtable1_64=1} - \ref{newtable4} represent numerical results of the WOPSIP method for $\delta = \frac{1}{64}, \frac{1}{128}, \frac{1}{256}$, respectively. Tables \ref{newtable1_64=3} - \ref{newtable8} represent numerical results of the WBCR method for $\delta = \frac{1}{64}, \frac{1}{128}, \frac{1}{256}$, respectively. The WBCR method showed a trend similar to that of the WOPSIP method. We observe that the use of the Shishkin mesh (\Roman{ltwo}) with the parameter $\tau = 4 \delta \ln(N)$ near the bottom is likely to achieve the optimal convergence order, whereas Tables \ref{newtable3} and \ref{newtable7} show that the mesh needs to be split sufficiently to obtain the optimum convergence order on the standard mesh (\Roman{lone}). These demonstrate the effectiveness of the anisotropic meshes. \\
\indent
We illustrate approximate solutions for the boundary layer problem can be found below  (Figure \ref{test2=wop_1_64} - Figure \ref{test2=wbcr_1_256}) when Mesh (\Roman{ltwo}), $N = 32$ and $\delta = \frac{1}{64}, \frac{1}{128}, \frac{1}{256}$.  Even from the perspective of approximate solutions, it can be observed that numerical oscillations are occurring in $u_{h,2}$ and $p_h$ in both methods. However, the oscillations of the WBCR method seem to be small. Furthermore, it can be observed that the oscillations do not spread throughout the domain. Numerical oscillations seem inevitable even if the Shishkin mesh is used. These results imply that other techniques may be needed to reduce unnatural plots near the boundary layer for approximate solutions in Navier--Stokes equations (e.g., \cite[Section 3.5]{ErnGue04}, \cite{RooStyTob08}).
}

\section{Concluding remarks} \label{sec=conc}
In conclusion, we identify several topics related to the results described in this work. 


\subsection{Discrete Poincar\'e and Sobolev inequalities}
The discrete Poincar\'e and Sobolev inequalities play important roles in finite element analysis. However, to the best of our knowledge, the derivation of those inequalities remains an open question in nonconforming finite element methods on anisotropic meshes, although those inequalities are derived in \cite{PieErn12} on shape-regular meshes. In this article, we have derived {a} discrete Poincar\'e inequality by imposing that $\Omega$ is convex. Because we used the regularity of the solution of the dual problem, the convexity can be not violated in our method. In \cite{Bre03}, the discrete Poincar\'e inequalities for piecewise $H^1$ functions are proposed. However, the inverse, trace inequalities and the local quasi-uniformity for meshes under the shape-regular condition are used for the proof. Therefore, careful consideration of the results used in \cite{Bre03} may be necessary to remove the assumption that $\Omega$ is convex. For example, it may be necessary to assume the minimum angle condition for simplices within macro elements. We leave further investigation of this problem as a topic for future work.

\subsection{Advantages of the WOPSIP method}
The dG scheme based on the WOPSIP method has advantages on anisotropic mesh partitions. One key advantage of the WOPSIP method is that it does not require any penalty parameter tuning. The Stokes element proposed in \eqref{CR7} satisfies the inf-sup condition. Another advantage of the method is that the error analysis of the technique is studied on more general meshes (\cite{BreOweSun12,BreOweSun08,Bre15,BarBre14}) than conformal meshes. This enables meshes with hanging nodes, whereas handling those meshes in the classical CR nonconforming finite element might be difficult. Here, we briefly treat the WOPSIP method on nonconforming meshes in our framework.

Let $d=2$. According to \cite[Section 2.7]{Bre15}, a nonconforming $\mathbb{T}_h$ on $\Omega$ is defined as satisfying the following condition:
\begin{itemize}
 \item If an edge of $T \in \mathbb{T}_h$ contains hanging nodes, it is the union of the edges of other triangles in $\mathbb{T}_h$.
\end{itemize}
The definition of the set $\mathcal{F}_h^i$ in Section \ref{regularmesh} is modified as follows: An (open) interior edge of a triangle in $\mathbb{T}_h$ belongs to $\mathcal{F}_h^i$ if and only if
\begin{description}
  \item[(\Roman{lone})] it is the common edge of exactly two triangles, or
  \item[(\Roman{ltwo})] it contains a hanging node.
\end{description}
An excellent feature of the WOPSIP method for nonconforming is that the CR interpolation defined in \eqref{CR3} is unaffected by the hanging nodes; that is, the relations \eqref{wop=13} and \eqref{stokes=21} remain valid. Therefore, if the relation \eqref{wop=4} holds, then the error estimates in {Theorem \ref{thr8}} may provide coverage as $h \to 0$ for nonconforming meshes. We leave further investigation of the hanging nodes for future work.

\subsection{Penalty parameters}
One disadvantage of the WOPSIP method is that it increases the number of conditions. In particular, the penalty parameter $\kappa_F$ increases as $h \to 0$. On the mesh (\Roman{ltwo}) in Section \ref{numerical=sec}, we compute the following qualities (Table \ref{table4}).
\begin{align*}
\displaystyle
&\tau_{\max}^{(f)} := \max_{F \in \mathcal{F}_h^i} \frac{1}{h_F}, \quad \tau_{\max}^{(ave)} :=  \max_{F \in \mathcal{F}_h^i} \frac{1}{4} \left( \frac{1}{{ \ell_{T_1,F}}} + \frac{1}{{ \ell_{T_2,F}}} \right),\\
&\tau_{\max}^{(dg)} :=  \max_{F \in \mathcal{F}_h^i} \frac{2}{ \left( \sqrt{ \ell_{T_1,F}} + \sqrt{\ell_{T_2,F}} \right)^{2}}, \quad \tau_{\max}^{(wop)} :=  \max_{F \in \mathcal{F}_h^i} \frac{2}{h^2 \left( \sqrt{ \ell_{T_1,F}} + \sqrt{\ell_{T_2,F}} \right)^{2}}
\end{align*}

\begin{table}[h]
\caption{Comparison of penalty parameters: Mesh (\Roman{ltwo}): $\delta = \frac{1}{1024}$}
\centering
\begin{tabular}{l|l|l|l|l|l} \hline
$N $&  $1/h $ & $\tau_{\max}^{(f)}$  & $\tau_{\max}^{(ave)}$ & $\tau_{\max}^{(dg)}$  & $\tau_{\max}^{(wop)}$  \\ \hline \hline
16&  7.2179e+00 & 7.3866e+02 & 3.6942e+02    &3.6942e+02  & 1.9246e+04   \\
32& 1.4467e+01    & 1.1819e+03 &   5.9114e+02 & 5.9114e+02 & 1.2373e+05 \\
64 & 2.8998e+01 & 1.9698e+03 &   9.8540e+02  &  9.8540e+02    & 8.2860e+05  \\
128 & 5.8123e+01   & 3.3767e+03 & 1.6896e+03    & 1.6896e+03 & 5.7079e+06   \\
256 &1.1650e+02  &  5.9093e+03  & 2.9574e+03     & 2.9574e+03   & 4.0139e+07  \\
\hline
\end{tabular}
\label{table4}
\end{table}

Table \ref{table4} shows that the penalty parameter of the WOPSIP method is considerably larger than the others. The problem of the ill-conditioning on anisotropic meshes {is left for future research.} The improvement of the ill-conditioning for the Poisson equation is stated in \cite{BreOweSun08}.

\begin{acknowledgements}
We also thank the anonymous referees for their valuable comments and suggestions.
\end{acknowledgements}

 \section*{Data Availibility}
 Data supporting the findings of this study are available from the corresponding author upon request. 

%
 \section*{Conflict of interest}
The authors declare no conflicts of interest associated with this manuscript.



\appendix
\section*{Additional numerical results in Section \ref{numerial=comp}}
In this section, we adopt the following schemes.

\begin{description}
   \item[(3)] Hybrid dG (HdG) method with a penalty term (c.f., \cite{EggWal13}). We define a discrete space as ${\Lambda}_h^1 := \{ {\lambda}_h \in L^2(\mathcal{F}_h)^2: \  {\lambda}_h|_F \in \mathbb{P}^1(F) ^2 \ \forall F \in \mathcal{F}_h, \ {\lambda}_h|_{\partial \Omega} = 0\}$. Find $(u_h, \lambda_h , p_h) \in V_{dc,h}^1 \times {\Lambda}_h^1 \times P_{dc,h}^{0}$ such that
\begin{subequations} \label{HdG=1}
\begin{align}
\displaystyle
a_h^{HdG}(u_h,\lambda_h;v_h , \mu_h) + b_h^{HdG}(v_h , \mu_h;p_h) &= \int_{\Omega} f \cdot v_h dx \quad \forall (v_h , \mu_h) \in V_{dc,h}^1 \times {\Lambda}_h^1, \label{HdG=1a} \\
b_h^{HdG}(u_h , \lambda_h ; q_h) - 10^{-10} \int_{\Omega} p_h q_h dx &= 0 \quad \forall q_h \in P_{dc,h}^{0}, \label{HdG=1b}
\end{align}
\end{subequations}
where
\begin{align*}
\displaystyle
a_h^{HdG}(u_h,\lambda_h;v_h , \mu_h) &:= \int_{\Omega} \nabla_h u_{h} :  \nabla_h v_{h} dx \\
&\hspace{-1.0cm} + \sum_{T \in \mathbb{T}_h} \sum_{F \in \mathcal{F}_{T}} \frac{{\eta}}{h_F} \int_{F} (\lambda_h - u_{h} ) \cdot ( \mu_h - v_{h} ) ds \\
&\hspace{-1.0cm} + \sum_{T \in \mathbb{T}_h} \int_{\partial T} \left(  \nabla_h u_{h}  n_F \cdot ( \mu_h -  v_{h} )+ ( \lambda_h - u_{h} )\cdot \nabla_h v_{h} n_F \right) ds, \\
 b_h^{HdG}(v_h , \mu_h; q_h) &:=  - \int_{\Omega} \divh v_h q_h dx - \sum_{T \in \mathbb{T}_h} \int_{\partial T} ( \mu_h - v_h ) \cdot n_F  q_h ds, 
\end{align*}
for any $(u_h,\lambda_h)  \in V_{dc,h}^1 \times {\Lambda}_h^1$, $(v_h,\mu_h) \in V_{dc,h}^1 \times {\Lambda}_h^1$, and $q_h \in P_{dc,h}^{0}$. Here, ${\eta}$ is a positive real number.
 \item[(4)] Conforming finite element method with a penalty parameter. Let $k \in \mathbb{N}$. We define discrete spaces as $P_{c,h}^{k} := \{ p_h \in H^1(\Omega); \ p_h|_T \in \mathbb{P}^k(T) \ \forall T \in \mathbb{T}_h \}$, and $V_{c,h}^2 := (P_{c,h}^{2} \cap H_0^1(\Omega))^2$. Find $(u_h,p_h) \in V_{c,h}^{2} \times P_{c,h}^{1}$ such that
 \begin{subequations} \label{dG=00}
\begin{align}
\displaystyle
a_h^{c}(u_h,v_h) + b_h^{c}(v_h , p_h) &= \int_{\Omega} f \cdot v_h dx \quad \forall v_h \in V_{c,h}^{2}, \label{dG=00a} \\
b_h^{c}(u_h , q_h)  - 10^{-10} \int_{\Omega} p_h q_h dx &= 0 \quad \forall q_h \in P_{c,h}^{1}, \label{dG=00b}
\end{align}
\end{subequations}
where 
\begin{align*}
\displaystyle
a_h^{c}(u_{h} , v_{h}) &:= \int_{\Omega} \nabla u_{h} :  \nabla v_{h} dx, \quad b_h^{c} (v_h,q_h) :=  - \int_{\Omega} \div v_h q_h dx, 
\end{align*}
for any $u_h  \in V_{c,h}^2$, $v_h \in V_{c,h}^2$, and $q_h \in P_{c,h}^{1}$.
\end{description}

For the computation of the schemes (3), and (4), we used the FreeFEM software tool \cite{HecPioMorHyaOht12,Hec} based on code provided in the prior work \cite{OikKik17} and used {UMFPACK}.

\begin{table}[h]
\caption{Scheme (3)}
\centering
\begin{tabular}{l | l |  l | l | l |  l | l | l | l} \hline
Mesh & $N$ & ${\eta}$ & $E_{u_h}^{(\text{Mesh No.})}$ & $r$ & $E_{u_h,L^2}^{(\text{Mesh No.})}$  & $r$ & $E_{p_h}^{(\text{Mesh No.})}$ & $r$ \\ \hline \hline
 \Roman{lone} & 32 & 20 &  7.65762e-02  &    &  6.54741e-03 &  &  2.35152e-02 &  \\
   & 64  & 20 & 3.74884e-02  & 1.03  &  1.60004e-03   & 2.03  & 1.14994e-02  & 1.03  \\ \hline
  \Roman{ltwo} & 32  & 30  &  1.08934e-01  &    & 1.27498e-02  &  &   3.78265e-02&     \\
 &  64 & 30  & 4.50132e-02  &  1.28  & 3.03009e-03 &  2.07&  1.84420e-02 & 1.04 \\ \hline
 \Roman{lthree}  & 32  & 40   & 8.22514e-02  &    & 1.01627e-02  &  & 2.92376e-02  &  \\
 & 64  & 40  & 3.88132e-02  &  1.08  & 2.50863e-03  &  2.02 &1.44096e-02  & 1.02 \\ \hline
 \Roman{lfour}  &  32 & 40  &  9.87191e-02  &    & 1.34999e-02  &  & 3.86220e-02 &   \\
  &  64 & 40  & 4.60525e-02  & 1.10    &3.34145e-03   & 2.01 & 1.90213e-02 & 1.02   \\
\hline
\end{tabular}
\label{table04}
\end{table}

\begin{table}[h]
\caption{Scheme (4)}
\centering
\begin{tabular}{l | l |  l | l | l |  l | l | l } \hline
Mesh & $N$ & $E_{u_h}^{(\text{Mesh No.})}$ & $r$ & $E_{u_h,L^2}^{(\text{Mesh No.})}$  & $r$ & $E_{p_h}^{(\text{Mesh No.})}$ & $r$ \\ \hline \hline
 \Roman{lone} & 32 &  2.86851e-03 &    & 8.51503e-05 &  &  2.44204e-04 &   \\
   & 64   & 7.18451e-04  & 2.00   & 1.06518e-05 & 3.00 & 6.09198e-05  & 2.00   \\ \hline
  \Roman{ltwo} & 32  & 5.11492e-03  &    &  2.31144e-04 &  & 5.25401e-04  &   \\
 &  64  & 1.19256e-03  &  2.10  & 2.66944e-05 & 3.11  & 1.10624e-04  & 2.25 \\ \hline
 \Roman{lthree}  & 32   & 3.35828e-03  &    & 1.19972e-04 &  & 3.56077e-04  &  \\
 & 64   & 8.41989e-04  & 2.00   & 1.50243e-05 & 3.00 &   8.89802e-05 & 2.00 \\ \hline
 \Roman{lfour}  &  32  & 4.60495e-03  &    & 2.02053e-04 &  &  4.35519e-04 &  \\
  &  64 &  1.16069e-03  &  1.99  & 2.52328e-05 & 3.00 &  1.08824e-04 & 2.00  \\
\hline
\end{tabular}
\label{tablep2p1}
\end{table}

The table \ref{table04} implies that the inf-sup condition is satisfied on {our} anisotropic meshes. However, to the best of our knowledge, a theoretical analysis of the HdG method on anisotropic meshes has not been considered in the relevant literature. Furthermore, one needs to tune up the parameter $\eta$. The table \ref{tablep2p1} also implies that the inf-sup condition is satisfied on {our} anisotropic meshes. The inf-sup condition in two-dimensional cases is proven in the previous work \cite{BarWac19}.

\end{document}